\documentclass[12pt]{article} 
\usepackage{amsmath,stmaryrd,pb-diagram}
\usepackage{amssymb} 
\usepackage{amscd}
\usepackage{amsthm}
\usepackage[utf8]{inputenc}
\usepackage[french]{babel}
\usepackage{hyperref}
\def\rot{\mathrm{rot}}
\def\grad{\mathrm{grad}}
\def\vol{\mathrm{vol}}

\def\dom{\mathrm{dom}}
\def\ad{\mathrm{ad}}

\def\Im{{\mbox{Im}}}

\def\dim{{\mbox{dim}}}

\def\det{{\mbox{det}}}
\def\GL{{\mbox{GL}}}
\def\O{{\mbox{O}}}

\def\End{{\mbox{End}}}

\def\Spin{{\mbox{Spin}}}
\def\cali{{\cal I}}
 
\def\cala{{\cal A}} 
\def\calb{{\cal B}}
\def\cald{{\cal D}} 
 
\def\call{{\cal L}}

\def\calo{{\cal O}} 
 
\def\calh{{\cal H}}

\def\calV{{\cal V}}

\def\pfraco{{\mathfrak o}}

 \def\fracg{{\mathfrak g}}

\def\fracu{{\mathfrak u}}
\def\fracS{{\mathfrak {S}}}

\def\fracco{{\mathfrak {co}}}  
 \def\fracgl{{\mathfrak {gl}}}

\def\bbbone{\mbox{\rm 1\hspace {-.6em} l}}

\def\Tr{{\mbox{Tr}}}

\numberwithin{equation}{section}

\begin{document}
\enlargethispage{3cm}
\thispagestyle{empty}
\begin{center}
{\bf NOTES SUR LES VARI\'ET\'ES DIFF\'ERENTIABLES,}
\end{center}
\begin{center}
 {\bf STRUCTURES COMPLEXES ET QUATERNIONIQUES}
\end{center} 
\begin{center}
{\bf ET  APPLICATIONS}
\end{center}
   
\vspace{2cm}

\begin{center}
Michel DUBOIS-VIOLETTE
\footnote{Laboratoire de Physique Th\'eorique, UMR 8627\\
CNRS et Universit\'e Paris Sud 11,
B\^atiment 210\\ F-91 405 Orsay Cedex\\
Michel.Dubois-Violette$@$u-psud.fr} 
\end{center}
\vfill

\newpage

\section*{Introduction}
Ces notes constituent une introduction à la géométrie différentielle, à la géométrie  des variétés complexes et à celle des variétés quaternioniques en vue d'applications en physique théorique. On a supposé une certaine familiarité avec les notions classiques de géométrie différen\-tielle telles que les notions de fibré tangent et de connexion sur un fibré principal. Les exemples et la manière de les aborder sont souvent issus de problèmes liés à la physique théorique. On trouvera sous la mention ``exercices et compléments" des développements et des éclairages différents de certains aspects et notions introduites.\\

Afin d'avoir un cadre unifié, nous sommes parti des notions de pseudogroupe et de $G$-structure qui sont décrites dans le chapitre 1 où nous donnons une introduction descriptive des conditions d'intégrabilité des $G$-structures.\\

Le chapitre 2 est consacré à des théorèmes classiques d'intégrabilité dans $\mathbb C^n$. On y trouvera une démonstration d'un théorème de Malgrange qui est la contre-partie non-abélienne du lemme de Grotendieck-Dolbault dont nous donnons également une démonstration.\\

Dans le chapitre 3 nous introduisons les variétés presque complexes et les variétés complexes ainsi que les notions d'application $(\pm)$ holomorphe et d'application harmonique. On trouvera dans les exemples une description en termes de projecteurs hermitiens de la géométrie et de la structure complexe des espaces projectifs et des grassmaniennes complexes.\\

Dans le chapitre 4, nous introduisons les fibrés holomorphes et nous en donnons une description en termes de fibrés différentiables munis de connexions particulières à l'aide d'un théorème de Koszul et Malgrange et d'un théorème de Singer que nous démontrons.\\

Les chapitres 5, 6, 7, 8, 9 et 10 sont essentiellement une reproduction du texte d'un exposé fait au séminaire de mathématique de l'E.N.S. durant l'hiver 1981 sur les strutures complexes au-dessus des variétés. On y trouvera une démonstration du théorème d'Atiyah et Ward concernant la transformation de Penrose pour les instantons. On développe, dans ce cadre des structures complexes au-dessus des variétés, l'approche géométrique d'Elie Cartan à la théorie des spineurs ainsi que la théorie des structures spinorielles. Le chapitre 10 contient en particulier une version de la transformation de Penrose pour le noyau de l'opérateur de Dirac sur une variété riemannienne compacte à spin de dimension paire qui est conforme plate (Corollaire 10.9).\\

Le chapitre 11 est consacré aux variétés quaternioniques. On y trouvera des exemples reliés à la théorie de la gravitation (instantons gravitationnels, etc.).\\

Le chapitre 12 est consacré aux équations de Yang et Mills.  On en donne différentes formulations dont certaines sont peu connues et on décrit divers types de solutions de ces équations. Ce chapitre contient notamment le rappel d'un résultat de Harnad, Tafel et Shnider généralisant les résultats de Trautman et Novakowski relatifs aux connexions canoniques des fibrés de Stieffel. On y trouvera également quelques éléments complémentaires (par rapport à ceux du chapitre 8) sur les instantons de Yang et Mills. Finalement, on donne une interprétation des équations de Yang et Mills comme lois de conservation sur le fibré principal et également une interprétation comme conditions d'intégrabilité de systèmes linéaires.\\

Le chapitre 13 porte sur les équations d'Einstein et l'approche d'Einstein-Cartan. On développe une interprétation des équations d'Einstein comme conditions d'intégrabilité de systèmes linéaires, ce qui permet de les identifier à des conditions de courbure nulle pour des connexions appropriées. Comme pour les équations de Yang et Mills (dans le cas semi-simple), on montre que les équations d'Einstein sont équivalentes à une loi de conservation sur le fibré des repères. Cette loi de conservation est apparentée à la conservation de l'impulsion-énergie.\\

Les notations utilisées ici sont standard. Sauf mention particulière la convention d'Einstein de sommation des indices répétés en haut et en bas (contravariants et covariants) est utilisée dans ces notes.

\newpage
\tableofcontents
\newpage

\section{Pseudogroupes, variétés et $G$-structures}

\subsection{Pseudogroupes de transformations} 

Soit $X$ un espace topologique; un {\sl pseudogroupe de transformations sur $X$} est un ensemble $\Gamma$ tel que \\
a) $f\in \Gamma$ est un homéomorphisme d'un ouvert de $X$, $\dom(f)$ , sur un autre ouvert de $X$, $\Im(f)$;\\

\noindent b) pour toute famille $(\calo_\alpha)_{\alpha\in I}$ d'ouvert de $X$, un homoméomorphisme $f$ de $\cup_{\alpha\in I} \calo_\alpha$ sur un autre ouvert de $X$ est un élément de $\Gamma$ si et seulement si $f\restriction \calo_\alpha\in \Gamma$, $\forall\alpha\in I$.\\

\noindent c) $f_1,f_2\in \Gamma$ et $\dom(f_1)=\Im(f_2)\Rightarrow f_1\circ f_2\in \Gamma$\\

\noindent d) $Id_X\in \Gamma$\\

\noindent e) $f\in \Gamma\Rightarrow f^{-1}\in \Gamma$

\subsection{Exemples}

a) $X=\mathbb R^n$, $\Gamma(\mathbb R^n)$ est le pseudo-groupe des difféomorphismes locaux de classe $C^r$ ainsi que leurs inverses; ($C^0$= continu, $C^\omega$=analytique)

b) $X=\mathbb R^n, G\subset\GL(n,\mathbb R)$ un sous-groupe et $r\geq 1$
\[
\Gamma^r_G(\mathbb R^n)=\{f\in \Gamma^r(\mathbb R^n)\vert df(x)\in G,\>\>\> \forall x\in \dom(f)\}
\]
Si $1\leq r\leq s$ et $H\subset G$ un sous-groupe, on a $\Gamma^s_H(\mathbb R^n)\subset \Gamma^r_G(\mathbb R^n)$. {\sl Nous poserons}
$\Gamma(\mathbb R^n)=\Gamma^\infty(\mathbb R^n),\Gamma_G(\mathbb R^n)=\Gamma^\infty_G(\mathbb R^n)$

c) $X=\mathbb C^n, \Gamma(\mathbb C^n)$ est le pseudogroupe des transformations holomorphes. On notera que, en considérant $\GL(n,\mathbb C)\subset \GL(2n,\mathbb R)$, on a canoniquement :
\[
\Gamma(\mathbb C^n)=\Gamma^r_{\GL(n,\mathbb C)}(\mathbb R^{2n}), \>\>\> \forall r\geq 1.
\]

d) Comme on le voit sur l'exemple c) l'équation de b) $df(x)\in G$ implique souvent une régularité pour $f$ plus forte que celle que l'on a mise au départ. Il peut même arriver que les restrictions deviennent injectives par exemple, $\Gamma^r_{O(n)}(\mathbb R^n)$ est l'ensemble des déplacements euclidiens.

\subsection{Atlas et variétés} 

Soit $\Gamma$ un pseudogroupe de transformation sur l'espace topologique $X$ et $V$ un espace topologique. Un {\sl atlas de $V$ compatible avec $\Gamma$} est une famille de paires $(\calo_\alpha, \varphi_\alpha)_{\alpha\in A}$ où $(\calo_\alpha)_{\alpha\in A}$ est un recouvrement ouvert de $V,\varphi_\alpha$ est un homéomorphisme de $\calo_\alpha$ sur un ouvert de $X$ et où l'on a : $\forall \alpha, \beta\in A$
\[
\varphi_\alpha\circ \varphi^{-1}_\beta \restriction \varphi_\beta(\calo_\alpha\cap \calo_\beta)\in\Gamma
\]
Les $(\calo_\alpha, \varphi_\alpha)$ $(\alpha\in A)$ sont les {\sl cartes de l'atlas}. Une paire $(\calo,\varphi)$ où $\calo$ est un ouvert de $V$ et $\varphi$ un homoméomorphisme de $\calo$ sur un ouvert
de $X$ est une {\sl carte compatible avec l'atlas} ($\calo_\alpha,\varphi_\alpha)_{\alpha\in A}$ si $((\calo_\alpha,\varphi_\alpha)_{\alpha\in A},(\calo,\varphi))$ est encore un atlas de $V$ compatible avec $\Gamma$. En adjoignant à un atlas de $V$ compatible avec $\Gamma$ toutes les cartes compatibles on obtient un {\sl atlas complet} sur $V$ compatible avec $\Gamma$ (i.e. contenu dans aucun autre atlas sur $V$ compatible ave $\Gamma$) contenant l'atlas  initial.

Soit $n$ un entier; {\sl une variété de classe $C^r$ et de dimension $n$} est un espace topologique séparé muni d'un atlas complet compatible avec $\Gamma^r(\mathbb R^n)$. Nous parlerons de {\sl variétés différentiables} pour désigner les variétés de classe $C^\infty$ et de {\sl variétés analytiques réelles} pour les variétés $C^\omega$.\\
Si $G$ est un sous-groupe de $\GL(n,\mathbb R)$, une {\sl variété de type $G$} sera un espace topologique séparé muni d'un atlas complet compatible avec $\Gamma_G(\mathbb R^n)$, (on peut définir de manière analogue les {\sl variétés de type $G$ de classe $C^r$} pour $r\geq 1$). Une variété de type $G$ est canoniquement une variété différentiable et, plus généralement si $H$ est un sous-groupe de $G$, une variété de type $H$ est canoniquement une variété de type $G$.

\subsection{Exemples}

a) Considérons $\GL(n,\mathbb C)$ comme sous-groupe de $\GL(2n,\mathbb R)$. Une variété de type $\GL(n,\mathbb C)$ est un espace topologique muni d'un atlas complet compatible avec $\Gamma (\mathbb C^n)$; une telle variété est appelée {\sl variété complexe}. C'est canoniquement une variété analytique. Nous reviendrons sur cet exemple.

b) Soit $Sp(n,\mathbb R)$ le sous-groupe de $\GL(2n,\mathbb R)$ laissant invariant la forme bilinéaire antisymmétrique $\omega_0=\sum^n_{k=1} dx^k\wedge dx^{n+k}$ (groupe symplectique). La notion de variété de type $Sp(n,\mathbb R)$ est identique à la notion de {\sl variété symplectique} de dimension $2n$, i.e. de variété différentiable de dimension $2n$ munie d'une 2-forme extérieure fermée partout non dégénérée. En effet si $V$ est une variété de type $Sp(n,\mathbb R)$ et si $(\calo,\varphi)$ appartient à l'atlas correspondant la 2-forme $\varphi^\ast(\omega_0)=\omega$ sur $\calo$ est fermé et sa valeur en un point $x$ ne dépend pas de la carte choisie $(\calo, \varphi)$ telle que $x\in \calo$ en vertu de la compatibilité avec $\Gamma_{Sp(n,\mathbb R)}(\mathbb R)$; on a donc une 2-forme fermée $\omega$ globalement définie sur $V$ de rang maximum. Inversement, le théorème de Darboux assure l'existence sur toute variété symplectique de dimension $2n$ d'un atlas compatible avec $\Gamma_{Sp(n,\mathbb R)}$ tel que la forme symplectique soit obtenue de la manière précédente.

c) Soit $O(n)\subset \GL(n,\mathbb R)$; une variété de type $O(n)$ n'est autre qu'une variété riemannienne plate de dimension $n$ comme on le voit facilement. On voit, à travers ces exemples (nous en verrons d'autres) que l'intérêt de cette notion de variété de type $G(\subset \GL(n,\mathbb R))$ est en quelque sorte proportionnel à la taille de $\Gamma_G(\mathbb R^n)$. Nous allons maintenant introduire une notion plus faible qui reste très riche même lorsque $\Gamma_G(\mathbb R^n)$ est petit.

\subsection{$G$-structures}

Soit $V$ une variété différentiable de dimension $n$, soit $L(V)$ le fibré des repères tangents à $M$. $L(V)$ est un fibré principal de groupe structural $\GL(n,\mathbb R)$. Soit $G$ un sous-groupe de $\GL(n,\mathbb R)$, une $G$-{\sl structure sur} $V$ est un fibré principal $P=P(V,G)$ sur $V$ de groupe structural $G$ qui est un sous-fibré différentiable de $L(V)$. Une variété différentiable munie d'une $G$-structure sera appelée {\sl variété presque de type $G$.}

\subsection{Exemples}

a) On considère à nouveau $\GL(n,\mathbb C)$ comme sous-groupe de $\GL(2n,\mathbb R)$; explicitement $\mathbb C^n\leftrightarrow \mathbb R^{2n}$ $(z^k)\leftrightarrow (x^k,x^{n+k})$ avec $z^k=x^k+ix^{n+k}$ et $\GL(n,\mathbb C)$ s'identifie au sous-groupe de $\GL(2n,\mathbb R)$ des matrices inversibles commutant avec $J_0=\left(\begin{array}{cc}
0 & -1_n\\
1_n & 0
\end{array}
\right)$ qui représente la multiplication par $i$. Alors, une $\GL(n,\mathbb C)$-structure sur une variété différentielle $V$ de dimension $2n$ n'est autre qu'une {\sl structure presque complexe sur $V$}, i.e. un champ (différentiable) $x\mapsto J_x$ ($x\in V$) d'endomorphismes des espaces tangents $(J_x\in \End T_x(V))$ satisfaisant $J^2_x=-\bbbone_x$ ($\Rightarrow$ {\sl presque complexe}= presque de type $\GL(n,\mathbb C))$.

b) Plus généralement, si $\tau_0$ est un tenseur sur $\mathbb R^n$ et si $G\subset \GL(n,\mathbb R)$ est un sous-groupe laissant $\tau_0$ invariante, une $G$-structure $P$ sur $V$ permet de construire sur $V$ un champ tensoriel $x\mapsto\tau_x$ du type précédent ; (si $u\in P$ est un repère de $P$ en $x$, $\tau_x$ est le tenseur de composantes $\tau_0$ dans $u$). Si $G$ est exactement le sous-groupe $G_{\tau_0}$ de $\GL(n,\mathbb R)$ laissant $\tau_0$ invariant, le champ $x\mapsto \tau_x$ permet de reconstruire la $G$-structure $P$ (= ensemble des repères relativement auxquels les composantes du champ $\tau_x$ sont $\tau_0$).

c) Une $O(n)$-structure sur $V$ est une {\sl structure riemannienne sur $V$}, (voir b) pour retrouver la définition usuelle). Une variété presque de type $O(n)$ n'est donc autre qu'une variété riemannienne.

d) Soit $G\subset \GL(n,\mathbb R)$ un sous-groupe et soit $V$ une variété de type $G$; l'ensemble $P$ des repères naturels correspondant aux cartes de l'atlas complet de $V$ compatible avec $\Gamma_G(\mathbb R^n)$ est manifestement une $G$-structure sur $V$. Une variété de type $G$ est donc (canoniquement) presque de type $G$. La $G$-structure que nous venons d'introduire est d'un type bien particulier.

\subsection{Intégrabilité}
Une $G$-structure $P$ sur $V$ est dite {\sl intégrable} si il existe un atlas de $V$ (compatible avec $\Gamma(\mathbb R^n)$, $n=\dim(V))$ dont les repères naturels (ceux correspondant à ses cartes) sont tous des éléments de $P$ (chaque carte définit une section de $P$ sur son domaine). L'intégrabilité est donc équivalente à l'existence d'un atlas compatible avec $\Gamma_G(\mathbb R^n)$ contenu dans l'atlas de $V$ et donc à l'existence d'un atlas complet de $V$ compatible avec $\Gamma_G(\mathbb R^n)$ et contenu dans l'atlas de $V$ (compatible avec $\Gamma(\mathbb R^n)$). Il est clair que cet atlas complet compatible avec $\Gamma_G(\mathbb R^n)$ est alors unique. Autrement dit, une variété de type $G$ n'est autre qu'une variété munie d'une $G$-structure intégrable.

\subsection{Proposition}
{\sl Soit $P$ une $G$-structure intégrable sur $V$; $P$ admet une connexion sans torsion.}

\noindent \underbar{Démonstration}. Soit $(\calo_\alpha)$ un recouvrement ouvert localement fini de $V$ par des domaines de cartes $(\calo_\alpha,\varphi_\alpha)$ de l'atlas complet de $V$ compatible avec $\Gamma_G(\mathbb R^n)$ correspondant à la $G$-structure intégrable $P$ (comme dans 1.7), soit $(\chi_\alpha)$ une partition de l'unité subordonnée à $(\calo_\alpha)$ et soit $\omega_\alpha$ la forme de connexion plate sur $P\restriction \calo_\alpha$ pour laquelle le champ des repères naturels de la carte $(\calo_\alpha, \varphi_\alpha)$ est horizontal. Posons $\omega=\sum_\alpha \pi^\ast (\chi_\alpha)\omega_\alpha$, où $\pi:P\rightarrow V$ est la projection; $\omega$ est une forme de connexion sans torsion sur $P$.~$\square$

\subsection{Obstructions et intégrabilité formelle}

Etant donné une $G$-structure $P$ sur $V$, il y a un certain nombre d'obstructions ``formelles" à l'intégrabilité de $P$; dans les ``bons cas", l'annulation de ces obstructions (= {\sl intégrabilité formelle}) entraîne l'intégrabilité. Par exemple, une $O(n)$-structure sur $V$ est, comme on l'a vu plus haut, une structure riemannienne sur $V$; la courbure riemannienne correspondant  (i.e. la courbure de la connexion sans torsion sur la $O(n)$-structure donnée sur $V$) est une obstruction à l'intégrabilité, (c'est la seule dans ce cas), et son annulation implique l'intégrabilité. Dans le cas général la proposition 1.8 est la clef pour construire ces obstructions. L'idée de l'intégrabilité d'une $G$-structure $P$ correspond à l'existence, au voisinage de chaque point $x\in V$, de coordonnées locales satisfaisant un certain système différentiel $S^p$ (compatibilité de la carte avec $\Gamma_G(\mathbb R^n)$). On fait, en $x\in V$ un développement limité de ces coordonnées locales, en identifiant dans le système précédent les termes d'ordre $k$, on obtient en $x$ un système algébrique $S^P_k(x)$ qui est, en général surdéterminé; on a donc des conditions de compatibilité que l'on peut mettre sous la forme $c_{k-1}(x,P)=0$ où les $c_k$ sont des sections de fibrés associés à $P$ dont nous allons donner une description. Pour faire les développements limités, on est amené à étudier le comportement local des {\sl automorphismes infinitésimaux} d'une $G$ structure $P$ sur $V$. Un tel automorphisme infinitésimal est un champ de vecteur $X$ sur $V$ engendrant un groupe local à 1 paramètre laissant $P(\subset L(V))$ invariant. Puisque l'on ne s'intéresse qu'au comportement local, on peut supposer $V=\mathbb R^n$ et $P=\mathbb R^n\times G$, $X=X^\mu \frac{\partial}{\partial x^\mu}$ (où $x^\mu$ est le système de coordonnées naturel de $\mathbb R^n$) est un automorphisme infinitésimal de $P$ si et seulement si la matrice $(\partial X^\mu/\partial x^\nu=(\partial _\nu X^\mu)$ appartient à l'algèbre de Lie $\fracg$ de $G$. Faisons un développement limité de $X$ à l'origine $X^\mu\sim \sum_k \frac{1}{k!} \sum_{(\nu)} X^\mu_{\nu_1\dots \nu_k}$ $x^{\nu_1}\dots x^{\nu_k}$ pour $\nu_2\dots \nu_k$ fixés, la matrice $X^\mu \nu\nu_2\dots \nu_k$ appartient à $\fracg$. On est amené à la définition suivante.

\subsection{Prolongements de $\fracg \subset \End (\mathbb R^n)$}

Soit $\fracg$ une sous-algèbre de Lie de $\fracgl(n,\mathbb R)=\End(\mathbb R^n)$. Le {\sl prolongement d'ordre $k$} ($k=0,1,2,\dots $) de $\fracg$ est l'espace $\fracg^{(k)}$ des applications $(k+1)$-linéaires symétriques, $t$ de $\underbrace{\mathbb R^n \times\dots \times \mathbb R^n}_{(k+1)-\text{termes}}$ dans $\mathbb
R^n$ telles que $\forall v_1,\dots, v_k\in \mathbb R^n$ $v\mapsto t (v,v_1,\dots v_k)$ soit un élément de $\fracg (\subset \End(\mathbb R^n))$. Si $\fracg^{k_0}= \{0\}$, alors $\fracg^{(k)}=\{0\}, \>\>\> \forall k\geq k_0$ et on appelle {\sl ordre de } $\fracg$ le plus petit entier $k$ tel que $\fracg^{(k)}=\{0\}$. Si $\forall k$ $\fracg^{k}\not=\{0\}$, on dit que $\fracg$ {\sl est de type infini}. On dit enfin que $\fracg$ {\sl est elliptique} si $\fracg$ ne possède aucun élément de rang 1. {\sl On pose} $\fracg^{(-1)}=\mathbb R^n$, (on a $\fracg^{(0)}=\fracg$).

\subsection{Proposition}
{\sl Si $\fracg$ est d'ordre fini, $\fracg$ est elliptique.}\\

\noindent\underbar{Démonstration}. Soit $e\otimes \omega\not=0$ ($e\in \mathbb R^n$, $\omega\in (\mathbb R^n)^\ast)$ une matrice de rang 1. Si $e\otimes \omega\in \fracg$ alors $t$ définit par $(v_0,v_1,\dots, v_k)\mapsto e.\omega(v_0)\omega(v_1)\dots \omega(v_k)=t(v_0,\dots,v_k)$ et un élément de $\fracg^{(k)}$ différent de zéro. Par conséquent si $\fracg$ possède un élément de rang 1 i.e. si $\fracg$ n'est pas elliptique, $\fracg$ est de type infini. D'où la proposition.~$\square$

\subsection{L'algèbre de Lie graduée $\fracg^{(\bullet)}$} 

Les automorphismes infinitésimaux de $\mathbb R^n\times G$ forment une algèbre de Lie (pour le crochet des champs de vecteurs); on en déduit, par développement en $O\in \mathbb R^n$ une structure d'{\sl algèbre de Lie graduée} sur $\fracg^{(\bullet)}=\oplus_{k\geq -1}\fracg^{(k)}$ (que l'on peut identifier aux automorphismes infinitésimaux à composantes polynomiales en $(x^\mu)$) : $[\fracg^{(-1)},\fracg^{(-1)}]=0$ et $[t_p,t_q](v_0,v_1\dots,v_{p+q})=\\
=\frac{1}{p!(q+1)!}\sum_\pi t_p(t_q(v_{\pi(0)},\dots,v_{\pi(q)})$,$v_{\pi(q+1)}\dots v_{\pi(q+p)})$\\
$- \frac{1}{q!(p+1)}\sum_\pi t_q(t_p(v_{\pi(0)},\dots v_{\pi(p)})$, $v_{\pi(p+1)},\dots,v_{\pi(p+q)})$, $\forall t_p\in \fracg^{(p)}, t_q\in \fracg^{(q)}$

\subsection{Le complexe $\fracg^{(\bullet)}\otimes \wedge \mathbb R^n$}

Soit $\stackrel{r}{\vee}(\mathbb R^n)^\ast$ la r-ième puissance symétrique et $\stackrel{s}{\wedge}(\mathbb R^n)^\ast$ la s-ième puissance extérieure de $(\mathbb R^r)^\ast$ et\\
 $\delta:(\stackrel{r}{\vee}\mathbb R^n)^\ast)\otimes 
(\stackrel{s}{\wedge}\mathbb R^n)^\ast)\rightarrow (\stackrel{r-1}{\vee}(\mathbb R^n)^\ast)\otimes (\stackrel{s+1}{\wedge}(\mathbb R^n)^\ast)$\\
 l'application linéaire définie par :
\[
\delta ((\stackrel{r}{\otimes} v)\otimes (v_1\wedge\dots \wedge v_s))=(\stackrel{r-1}{\otimes}v)\otimes (v\wedge v_1\wedge \dots \wedge v_s)
\]
avec $v, v_i\in \mathbb R^{n\ast}$.
On a un complexe (i.e. $\delta^2=0$)
\[
0\rightarrow \stackrel{r}{\vee}(\mathbb R^n)^\ast\stackrel{\delta}{\rightarrow} \stackrel{r-1}{\vee}(\mathbb R^n)^\ast \otimes \stackrel{1}{\wedge}(\mathbb R^n)^\ast \stackrel{\delta}{\rightarrow}\dots \stackrel{\delta}{\rightarrow}\stackrel{r}{\wedge}(\mathbb R^n)^\ast\rightarrow 0
\]
qui n'est autre que le complexe des formes différentielles polynomiales sur $\mathbb R^n$.
Par tensorisation avec $\mathbb R^n$ en tenant compte de $\fracg\subset \mathbb R^n\otimes (\mathbb R^n)^\ast$ et des définitions, on obtient le sous-complexe :
\[
0\rightarrow \fracg^{(r-1)}\stackrel{\delta}{\rightarrow} \fracg^{(r-2)}\otimes (\mathbb R^n)\stackrel{\ast}{\rightarrow}\fracg^{(r-3)}\otimes \stackrel{2}{\wedge}(\mathbb R^n)\stackrel{x}{\rightarrow} \dots\stackrel{\delta}{\rightarrow}\fracg^{(-1)}\otimes\stackrel{r}{\wedge}(\mathbb R^n)^\ast\rightarrow 0
\]
\noindent\underbar{Posons} :\\
\[
H^{r,s}(\fracg)=\{K\in \fracg^{r-1}\otimes \stackrel{s}{\wedge}(\mathbb R^n)^\ast\vert \delta K=0\}/\delta (\fracg^{(r)}\otimes \stackrel{s-1}{\wedge}(\mathbb R^n)^\ast)
\]
i.e. ce sont les groupes d'homologie du complexe précédent. On a $H^{r,1}(\fracg)=0$ $\forall r\geq 1$. La partie importante pour l'intégrabilité formelle est $H^{r,2}(\fracg)$ pour $r\geq 0$ ($H^{\ell,2}(\fracg)=0$ si $\fracg$ est d'ordre $k$ et $\ell\geq k+1$).

\subsection{Les obstructions $c_k(P)$}
Soit $P$ une $G$-structure sur $V$. $G$ agit sur $\mathbb R^n$, sur $\fracg$ (action adjointe) et sur les $\fracg^{(k)}$ en général. On a donc {\sl des fibrés} (vectoriels) {\sl associés à $P$ sur $V$ correspondant} que nous noterons $\underline{\fracg}^{(k)}$, (on a $\underline{\fracg}^{(-1)}=T(V))$. De même $G$ agit sur les espaces $H^{r,s}(\fracg)$ et on a des fibrés associés à $P$ sur $V$ correspondant $\underline{H^{r,s}(\fracg)}$. On peut définir directement $\underline{H^{r,s}(\fracg)}$ comme l'homologie de 
\[
\rightarrow \underline{\fracg}^{(r)}\otimes \stackrel{s-1}{\wedge} T^\ast(V)\stackrel{\delta}{\rightarrow} \underline{\fracg}^{(r-1)}\otimes\stackrel{s}{\wedge}T^\ast(V) \stackrel{\delta}{\rightarrow}\underline{\fracg}^{(r-2)}\otimes \stackrel{s+1}{\wedge} T^\ast(V)\rightarrow \dots
\]
Les obstructions à l'intégrabilité de $P$ sont des sections $c_k(P)\in \Gamma\underline{(H^{k,2}(\fracg))}$ des fibrés $\underline{H^{k,2}(\fracg)}$ sur $V$ pour $k\geq 0$. Commençons par définir $c_0$.

\subsection{Définition de $c_0(P)$}
Soit $P$ une $G$-structure sur $V$ et soient $\omega$ et $\omega'$ deux formes de connexions sur $P$ de torsions $\tau$ et $\tau'$. On peut considérer que $\tau$ et $\tau'$ sont dans \\
$T(V)\otimes \stackrel{2}{\wedge}T^\ast(V)=\underline{\fracg}^{(-1)}\otimes \stackrel{2}{\wedge}T^\ast(V)$ et que $\omega-\omega'$ est dans  $\underline{\fracg}\otimes T^\ast(V)=\underline{\fracg}^{(0)}\otimes \stackrel{1}{\wedge}T^\ast(V)$. Par définition, on a
\[
\tau-\tau'=\delta(\omega-\omega'),\>\>\> \delta\tau=\delta\tau'=0
\]
autrement dit la classe de $\tau$ dans $\Gamma(\underline{H^{0,2}(\fracg))}$ ne dépend pas de la connexion $\omega$ choisie mais uniquement de $P$; on l'appelle {\sl première fonction de structure} de la $G$-structure $P$ sur $V$ et {\sl on la note} $c_0(P)$, (on peut aussi la définir comme fonction sur $P$ à valeurs dans $H^{0,2}(\fracg)$  au lieu de section de $\underline{H^{0,2}(\fracg)}$ comme ici).\\
L'annulation de $c_0(P)$ est, par construction équivalente à l'existence d'une connexion sans torsion sur $P$ et 1.8 montre que c'est une condition nécessaire à l'intégrabilité de $P$. C'est manifestement une condition du {\sl premier ordre} pour le système différentiel correspondant à l'intégrabilité.

\subsection{Prolongements de $G\subset \GL(n,\mathbb R)$}
Soit $G$ un sous-groupe de Lie de $\GL(n,\mathbb R)$ d'algèbre de Lie $\fracg$. \`A un élément $t$ de $\fracg^{(k)}$ ($k\geq 1$), on associe l'endomorphisme $T^{(k)}(t)$ de l'espace vectoriel $\oplus^{r=k}_{r=0}\fracg^{(r-1)}=\mathbb R^n\oplus \fracg\oplus \dots \oplus \fracg^{(k-1)}$ défini par : $T^{(k)}(t)v=v+t(v,\dots)\in \mathbb R^n\oplus \fracg^{(k-1)}$ si $v\in \mathbb R^n$ et
$T^k(t)x=x$ si $x\in \oplus^{r=k}_{r=1}\fracg^{r-1}$. $T^{(k)}$ ainsi défini est un homomorphisme de groupes, du groupe additif $\fracg^{(k)}$ dans le groupe linéaire réel $\GL(\oplus^{r=k}_{r=0}\fracg^{(r-1)})$ de l'espace vectoriel réel $\oplus^{r=k}_{r=0}\fracg^{(r-1)})$; son image $G^{(k)}=T^{(k)}(\fracg^{(k)})$ est un sous-groupe (abélien si $k\geq q$) du groupe linéaire $\GL(\otimes^{r=k}_{r=0} \fracg^{(r-1)})$
qui est appelé {\sl prolongement d'ordre $k$ de $G$}. On a $G^{(0)}=G$ et on notera que $\fracg^{(k)}$ est l'algèbre de Lie de $G^{(k)}$ et que $T^{(k)}$ est l'application exponentielle correspondante (ce qui justifie la terminologie).

Si $N^{(k)}=\dim\left(\oplus^{r=k}_{r=0}\fracg^{(r-1)}\right)$, $\GL\left(\oplus^{r=k}_{r=0}\fracg^{(r-1)}\right)=\GL(N^{(k)},\mathbb R)$ et $G^{(k+1)}$ {\sl est} canoniquement {\sl le prolongement d'ordre 1 de} $G^{(k)}$ i.e. $G^{(k+1)}\simeq (G^{(k)})^{(1)}$.

\subsection{Bases}
Soit $P$ une $G$-structure sur $V$ avec $\dim(V)=n$. Choisissons une fois pour toutes une base de l'algèbre de Lie $\fracg$ de $G$; à cette base correspond une base $\calb_u$ de l'espace des vecteurs tangents en $u\in P$ à la fibre de $P\stackrel{\pi}{\rightarrow}V$ passant par $u$ (i.e. $T_u(\pi^{-1}(\pi(u))$). Pour tout sous-espace $H$ de $T_u(P)$ supplémentaire à l'espace tangent à la fibre passant par $u$, il y a une base unique $\calb'_H$ de $H$ dont l'image par la projection $\pi$ est la base $u(\in P)$ de $T_{\pi(u)}(V)$; $(\calb_u,\calb'_H)$ est donc une base de $T_u(P)$ que nous noterons $u^{(1)}(u,H)$. L'ensemble des $u^{(1)}(u,H)$ est un sous-fibré $\tilde P^{(1)}$ du fibré des repères $L(P)$ de $P$ qui est invariant par $G^{(1)}\subset \GL(N^{(1)},\mathbb R)$, (on a évidemment $N^{(1)}=\dim(\mathbb R^n\oplus \fracg)=\dim(P)$).

\subsection{Prolongements de la $G$-structure $P$}
Soit $K$ un sous-espace de $\mathbb R^n\otimes \stackrel{2}{\wedge}(\mathbb R^n)^\ast=\fracg^{(-1)}\otimes \stackrel{2}{\wedge}(\mathbb R^n)^\ast$ supplémentaire au sous-espace $\delta(\fracg^{(0)}\otimes \stackrel{1}{\wedge}(\mathbb R^n)^\ast)=\delta(\fracg\otimes (\mathbb R^n)^\ast)$, (voir 1.13), et soit $P^{(1)}\subset \tilde P^{(1)}$ l'ensemble des $u^{(1)}(u,H)$ où $H$ satisfait $d\theta_u\circ i_H \in K$, $i_H$ désignant l'isomorphisme inverse de celui induit par $\theta_u\restriction H$ ($:H\stackrel{\simeq}{\rightarrow}\mathbb R^n$ et $\theta$ étant la 1-forme canonique à valeurs dans $\mathbb R^n$ ({\sl forme de soudure}) sur $P\subset L(V)$. $P^{(1)}$ est un sous-fibré principal de $L(P)$ dont le groupe de structure est $G^{(
1)}$ comme on le vérifie facilement ; c'est donc une $G^{(1)}$-structure sur $P$ que l'on appelle {\sl prolongement d'ordre 1 de la $G$-structure $P$ sur 
$V$}. On vérifiera que si on change le choix de $K$, on obtient une autre $G^{(1)}$-structure $P^{'(1)}$ sur $P$ qui est isomorphe à $P^{(1)}$, (comme $G^{(1)}$-fibré principal sur $P$), de sorte que le choix de $K$ n'est pas important.\\
Posons $P^{(0)}=P$ et définissons par récurrence les $P^{(k)}$ par $P^{(k+1)}=(P^ {(k)})^{(1)}$ de sorte $P^{(k)}$ est une $G^{(k)}$-structure sur $P^{(k-1)}$ appelée {\sl prolongement d'ordre $k$ de la $G$-structure $P$ sur $V$}. (On pourra poser $V=P^{(-1)}$ de la manière cohérente).

\subsection{Fonctions de structure de $P$ et définition des $c_k(P)$} 
La première fonction de structure $c_0(P^{(k)})$ de la $G^{(k)}$-structure $P^{(k)}$ sur $P^{(k-1)}$, $(k\geq 0)$, est parfois appelée $(k+1)$-{\sl ième fonction de structure de la $G$-structure $P$ sur $V$.} L'annulation de tous les $c_0(P^{(k)})$, $k\geq 0$, implique l'intégrabilité formelle des $P^{(k)}$ (c'est même équivalent), mais il y a des exemples où $P$ est intégrable sans que les $P^{(k)}$ le soient (même formellement), c'est pourquoi, nous éviterons cette terminologie. On peut cependant extraire de $c_0(P^{(k)})$ une partie invariante se projetant sur $V$ et correspondant à $c_k(P)\in \Gamma(\underline{H^{k,2}(\fracg)})$; c'est l'obstruction à l'intégrabilité formelle de la $G$-structure à l'ordre $k$.\\
L'annulation des $c_k(P)$ est équivalente à l'intégrabilité formelle de la $G$-structure. Lorsque $G$ est elliptique, un théorème de Malgrange implique que l'intégrabililté formelle est équivalente à l'intégrabilité. Plus récemment, ce résultat a été généralisé au cas de $G\subset \GL(n,\mathbb R)$ quelconque par Goldschmidt et Spencer et par Molino.

\subsection{Exercices et compléments}
\begin{enumerate}

\item
$\stackrel{r}{\vee}(\mathbb R^n)^\ast\otimes \stackrel{s}{\wedge}(\mathbb R^n)^\ast$ peut s'interpréter comme étant l'espace des formes différentielles de degré $s$ à coefficients polynômes homogènes de degré $r$ dans $\mathbb R^n$; l'opérateur $\delta:\stackrel{r}{\vee}(\mathbb R^n)^\ast\otimes \stackrel{s}{\wedge}(\mathbb R^n)^\ast\rightarrow \stackrel{r-1}{\vee}(\mathbb R^n)^\ast\otimes \stackrel{v+1}{\wedge}(\mathbb R^n)^\ast$ n'est autre que la différentielle extérieure sur les formes correspondantes. Le terme d'ordre $r$ du développement limité au voisinage de $O\in \mathbb R^n$ d'une forme différentielle de degré $s$ est donc un élément de\\
 $\stackrel{r}{\vee}(\mathbb R^n)^\ast\otimes \stackrel{s}{\wedge}(\mathbb R^n)^\ast$. La suite 
\[
\stackrel{r+1}{\vee}(\mathbb R^n)^\ast\otimes \stackrel{s-1}{\wedge}(\mathbb R^n)^\ast\stackrel{\delta}{\rightarrow}\stackrel{r}{\vee}(\mathbb R^n)\otimes \stackrel{s}{\wedge}(\mathbb R^n)^\ast\stackrel{\delta}{\rightarrow}\stackrel{r-1}{\vee}(\mathbb R^n)^\ast\otimes \stackrel{s+1}{\wedge}(\mathbb R^n)
\]
{\sl est exacte en} $\stackrel{r}{\vee}(\mathbb R^n)^\ast\otimes \stackrel{s}{\wedge}(\mathbb R^n)^\ast$, $\forall r,s\geq 0$, en posant
 $$\stackrel{r}{\vee}(\mathbb R^n)\otimes \stackrel{-1}{\wedge}(\mathbb R^n)^\ast=\left\{
\begin{array}{l}
0\ \text{si}\ r\geq 1\\
\mathbb R \ \text{si}\ r=0
\end{array}
\right.$$ 
c'est {\sl la version formelle} (i.e. au sens des séries formelles pour les dévelop\-pements limités en $O\in \mathbb R^n)$ {\sl du lemme de Poincaré}.
\item
Déduire de 1 (par tensorisation avec $\mathbb R^n)$ que $H^{r,s}(\fracgl(n,\mathbb R))=0$, $\forall r,s$ avec $r+s\geq 1$, où $\fracgl(n,\mathbb R)$ est l'algèbre de Lie ($\simeq \End\mathbb R^n)$ de $\GL(n,\mathbb R)$.
\item
Montrer que $H^{r,2}(\fracg)=0$ $\forall r\geq 1$ pour $\fracg$= algèbre de Lie de $Sp(\ell,\mathbb R)$ ou de $\GL(\ell,\mathbb C)$, (pour le dernier cas on peut s'inspirer de 1 et 2).
\item
L'algèbre de Lie $\pfraco(n)$ de $O(n)$ est d'ordre un $(\pfraco(n)^{(1)}=\{0\}$), donc $H^{r,s}(\pfraco(n))=0$ si $r\geq 2$. Montrer que $H^{0,2}(\pfraco(n))=0$ et que les éléments de $H^{1,2}(\pfraco(n))$ sont les tenseurs d'ordre quatre $R_{ab\ cd}$ possédant les symétries du tenseur de courbure de Riemann-Christoffel, (i.e. $R_{ab\ cd}=-R_{ba\ cd}=-R_{ab\ dc}=R_{cd\ ab}$ et $R_{ab\  cd}+R_{ac\ db}+R_{ad\ bc}=0$).
\item
Montrer que l'algèbre de Lie $\fracco(n)$ de $CO(n)$ (= le groupe linéaire conforme) est d'ordre deux; déterminer $\fracco(n)^{(1)}$. Montrer que les éléments de $H^{1,2}(\fracco (n))$ s'identifient canoniquement aux tenseurs d'ordre quatre $W_{ab\ cd}$ possédant les symétries du tenseur de courbure de Weyl. Montrer que $H^{r,2}(\fracco(n))=0$ si $r\not= 1$.
\item
Soient $\fracg_0 \subset \fracg$, alors $H^{0,k}(\fracg_0)=0\Rightarrow H^{0,k}(\fracg)=0$
\end{enumerate}

\newpage
\begin{center}
{\large Références pour le chapitre 1}
\end{center}
\vspace{1cm}
\begin{itemize}
\item
Kobayashi S. : Transformation groups in differential geometry. Springer Verlag, 1992.\\

\item
Dieudonné J. : \'Elements d'analyse, tome 4. Gauthier-Villars, 1991.\\

\item
Goldschmidt H. : The integrability problem for Lie equations. {\sl Bul. Am. Math. Soc.} {\bf 84} (1978), 531.\\

\item
Guillemin V. : The integrabiity problem for $G$-structures. {\sl Trans. Am. Math. Soc.} {\bf 116} (1965), 544.\\

\item
Malgrange B. :  \'Equations de Lie. {\sl J. Diff. Geom.} {\bf 6} (1972), 503, \\
\phantom{Malgrange, B. :  \'Equations de Lie. {\sl J. Diff. Geom.}} {\bf 7} (1972), 117.\\

\item
Molino P. :  Théorie des $G$-structures : le problème d'équivalence. Lecture Notes in Mathematics 988. Springer-Verlag, 1977.

\end{itemize}
\newpage

\section{Quelques théorèmes d'intégrabilité dans $\mathbb C^n$}
\subsection{Lemme}

{\sl
Soit $U$ un ouvert borné de $\mathbb C$ et $f$ une fonction de $C^k(U)$ bornée. La fonction $g$ sur $U$, définie par $g(z)=(2\pi)^{-1}\int_U(\xi-z)^{-1}f(\xi)d\xi\wedge d\bar \xi$ est dans $C^k(U)$ et, si $k\geq 1$, on a $\partial g/\partial \bar z=f$.}\\

\noindent\underbar{Démonstration}. La première partie $(g\in C^k(U))$ est immédiate. On a d'autre part, pour tout $D\subset U$ difféomorphe à un disque et toute fonction $g\in C^k(U)$ avec $k\geq 0$ et $z$ intérieur à $D$ ($z\in \stackrel{\circ}{D}$).
\[
\begin{array}{l}
\frac{1}{2\pi i}\left(\int_{\partial D}\frac{g(\xi)d\xi}{\xi-z}-\lim_{\varepsilon\rightarrow 0}\int_{\partial D_\varepsilon(z)} \frac{g(\xi)d\xi}{\xi-z}\right)=\\
\\
\frac{1}{2\pi i}\int_{\partial D} \frac{g(\xi)d\xi}{\xi-z}-g(z)=\frac{1}{2\pi i}\int_D \frac{\partial g(\xi)/\partial \bar\xi}{\xi-z}d\bar\xi \wedge d\xi,
\end{array}
\]
où $D_\varepsilon(z)$ est le disque de rayon $\varepsilon$ centré en $z$ et la dernière égalité provient de la formule de Stokes. La dernière partie du lemme résulte de la formule précédente, impliquant $\partial g/\partial \bar z=f$ dans $D$, par passage à la limite. $\square$

\subsection{Lemme}
{\sl
Soit $\calo$ un ouvert de $\mathbb C \times \mathbb R^m$ et soit $A:\calo\rightarrow \End(\mathbb C^N)$ une fonction de classe $C^k$ avec $k\geq 1$. Pour chaque $(z_0,\vec r_0)\in \calo$, on peut trouver un ouvert de $\calo$ de la forme $U_0\times W_0$, avec $z_0\in U_0\subset\mathbb C$ et $\vec r_0\in W_0 \subset \mathbb R^m$, et une fonction $G:U_0\times W_0\rightarrow \GL(N,\mathbb C)$ de classe $C^k$  tels que l'on ait :
$A=G^{-1}\frac{\partial G}{\partial \bar z}$ dans $U_0\times W_0$
}\\

\noindent \underbar{Démonstration}. Soit $K\in \mathbb R$ avec $0<K<1/2$. On peut trouver un ouvert borné de $\calo$ de la forme $U_1\times W_0$ avec $z_0\in U_1\subset \mathbb C$ et $\vec r_0\in W_0\subset \mathbb R^m$ tel que l'on ait :
\[
\sup_{(z,\vec r)\in U_1\times W_0}\left\{ \frac{1}{\pi}\int_{U_1}\bigparallel \frac{A(u+iv,\vec r)}{u+iv-z}\bigparallel du\  dv\right\}
 \leq K
 \]
  Soit $\chi$ une fonction à support compact dans $U_1$ à valeurs dans $[0,1]$ qui est égale à 1 dans un voisinage ouvert $U_0\subset U_1$ de $z_0$. Définissons les fonctions $F_n:\mathbb C\times W_0\rightarrow \End(\mathbb C^N)$ par
\[
F_0=\bbbone_{\mathbb C^N}\>\> \text{et}\>\> F_{n+1}(z,\vec r)=\frac{1}{2\pi i}\int \frac{F_n(\xi,\vec r)A(\xi,\vec r)\chi (\xi)}{\xi-z} d\xi\wedge d\bar \xi
\]  
On obtient, par récurrence sur $n$, que les $F_n$ sont de classe $C^k$ dans $U_1\times W_0$ et que, si $D$ est un opérateur différentiel d'ordre $p\leq k$ à coefficients constants, on a les estimations suivantes :
 \[
  \sup_{U_1\times W_0}\bigparallel F_n(z,\vec r)\bigparallel \leq K^n\>\>\> \text{et}\>\>\> \sup_{U_1\times W_0}\bigparallel DF_n(z,\vec r)\bigparallel \leq K^n \left(\sum^p_{q=0} C^q_n \frac{\beta_q(D)}{K^q}\right)
  \]
où les $\beta_q(D)$ sont des constantes. Ceci implique que $\sum_nF_n$ converge dans $U_1\times W_0$ vers une fonction de classe $C^k$ $G:U_1\times W_0\rightarrow \End(\mathbb C^N)$, et, comme on a 
\[
\bigparallel 1-G(z,\vec r)\bigparallel\leq \sum_{n\geq 1}\bigparallel F_n(z,\vec r)\bigparallel \leq \sum_{n\geq 1}K^n=\frac{K}{1-K}<1
\]
il s'ensuit que $G(z,\vec r)\in \GL(N,\mathbb C)$ sur $U_1\times W_0$. Il en résulte $\partial G/\partial \bar z=GA\Leftrightarrow A=G^{-1}\partial G/\partial \bar z$ dans $U_0\times W_0$.~$\square$

\subsection{Théorème}
{\sl
Soit $\calo$ un ouvert de $\mathbb C^n\times \mathbb R^m$ et soient $A_\alpha:\calo\rightarrow \End(\mathbb C^N)$ ($\alpha=1, 2, \dots, n$) $n$ fonctions de classe $C^k$ avec $k\geq 1$ satisfaisant les conditions 
\begin{equation}
\partial A_\beta/\partial\bar z_\alpha-\partial A_\alpha/\partial \bar z_\beta + [A_\alpha,A_\beta]=0\>\>\> \forall \alpha, \beta= 1,\dots,n.
\label{I}
\end{equation}
Alors chaque point $(\vec z_0,\vec r_0)\in \calo$ possède un voisinage ouvert $\calo_0\subset \calo$ sur lequel on a une fonction de classe $C^k$, $G:\calo_0\rightarrow \GL(N,\mathbb C)$, telle que $A_\alpha =G^{-1}\partial G/\partial \bar z_\alpha$ dans $\calo_0$}\\

\noindent \underbar{Démonstration}.
Le lemme 2.2 implique le résultat pour $n=1$ Supposons que ce résultat soit vrai pour $n=\nu-1$; $(\vec z_0,\vec r_0)\in \mathbb C^\nu\times \mathbb R^m=\mathbb C^{\nu-1}\times (\mathbb C\times \mathbb R^m)=\mathbb C^{\nu-1}\times \mathbb R^{m+2}$ a donc un voisinage ouvert $\tilde \calo_0\subset \calo$ dans lequel on a une fonction de classe $C^k,\tilde G:\tilde \calo_0\rightarrow \GL(N,\mathbb C)$ tels que $A_\alpha=\tilde G^{-1}\partial \tilde G/\partial \bar z_\alpha$ pour $\alpha=1,2,\dots,\nu-1$.\\
Définissons (``transformation de jauge inverse") $A'_\nu$ par $A_\nu=\tilde G^{-1}\partial \tilde G/\partial \bar z_\nu+ \tilde G^{-1}A'_\nu\tilde G$ dans $\tilde\calo_0$; (\ref{I}) est équivalent dans $\tilde \calo$ à $\partial A'_\nu/\partial z_\alpha=0$ pour $\alpha=1, 2, \dots,\nu-1$. On peut alors trouver un ouvert borné $U$ dans $\mathbb C$ et un ouvert $W$ de $\mathbb C^{\nu-1}\times \mathbb R^m$ tel que $\calo=\{(\vec z,\vec r)\in \mathbb C^\nu\times \mathbb R^m\vert z_\nu\in U, (z_1,\dots,z_{\nu-1},\vec r)\in W\}$ soit inclus dans $\calo_0$. Définissons la fonction $\phi$ de classe $C^k$ dans $\calo_0$ par : $\phi(\vec z,\vec r)=(2\pi i)^{-1}\int_U(\xi-z_\nu)^{-1}A'_\nu(z_1,\dots, z_{\nu-1},\xi,\vec r)d\xi\wedge d\bar \xi$. On a :
$A'_\nu=\partial\phi/\partial\bar z_\nu$, (d'après 2.1), et $\partial\phi/\partial \bar z_\alpha=0$ pour $a=1,2,\dots,\nu-1$ dans $\calo_0$. Posant $G=e^\phi\tilde G$ on a alors dans $\calo_0:A_\alpha=G^{-1}\partial G/\partial \bar z_\alpha$ pour $\alpha=1,\dots, \nu$. Le résultat est donc vrai pour $n=\nu$. Il est donc, par récurrence, vrai quel que soit $n$.~$\square$

\subsection{Formes différentielles pures}

L'algèbre extérieure des formes différentielles complexes de classe $C^k$ sur un ouvert $\calo$ de $\mathbb C^n$ est engendrée par les fonctions complexes de classe $C^k$ sur $\calo$ et les 1-formes $dz_\alpha,d\bar z_\alpha$ ($\alpha=1,2,\dots,n$). On a donc une notion de {\sl degrés partiels en $dz_\alpha$ et en $d\bar z_\alpha$} sur ces formes. Nous dirons qu'une {\sl forme est pure de type $(r,s)$} si elle est homogène de degré $r$ en $dz_\alpha$ et de degré $s$ en $d\bar z_\beta$ ($\alpha,\beta=1,\dots,n$); elle est alors évidemment de degré $p+q$.

La différentielle extérieure $d\lambda$ d'une forme pure $\lambda$ de type $(r,s)$ et de classe $C^k$ avec $k\geq 1$ sur $\calo$ est la somme d'une forme $\partial\lambda$ pure de type $(r+1,s)$ de d'une forme $\bar\partial \lambda$ pure de type $(r,s+1)$ de classe $C^{k-1}$ sur $\calo$. Par linéarité, {\sl on définit ainsi les opérateurs} $\partial$ et $\bar\partial$ sur les formes différentielles de classe $C^k$ avec $k\geq 1$ sur $\calo$ à valeurs dans les formes différentielles $C^{k-1}$ sur $\calo$; on a $d=\partial+\bar\partial$, et sur les formes de classe $C^k$ avec $k\geq 2$ on a $\partial^2=\bar \partial^2=\partial\bar\partial+\bar\partial\partial=0$.\\
Ces définitions se généralisent de manière évidente aux formes à valeurs dans des espaces vectoriels complexes.\\
{\sl Pour la suite}, pour simplifier les notations, nous nous placerons dans le cas $C^\infty$ et sauf mention spéciale, {\sl forme différentielle signifiera forme différentielle $C^\infty$}.

\subsection{Théorème}

{\sl
Soit $G$ un groupe de Lie complexe d'algèbre de Lie $\fracg$ et soit $\omega$ une forme différentielle à valeurs dans $\fracg$ pure de type (0,1) sur un ouvert $\calo$ de $\mathbb C^n$. Les conditions suivantes sont équivalentes :\\
(a) $\bar\partial\omega+\frac{1}{2}[\omega,\omega]=0$ sur $\calo$\\
(b) On a un recouvrement ouvert $(\calo_\alpha)$ de $\calo$ et des fonctions $C^\infty$ $g_\alpha:\calo_\alpha\rightarrow G$ tels que $\omega\restriction \calo_\alpha=g^{-1}_\alpha \bar\partial g_\alpha$ pour tout $\alpha$.
}

C'est un cas particulier de 2.3 car si dans 2.3 les $A_\alpha$ sont à valeurs dans une sous-algèbre de Lie complexe de $\End(\mathbb C^N)$, alors la fonction $G$ construite dans la démonstration de 2.3 est dans le sous-groupe de Lie de $\GL(N,\mathbb C)$ correspondant.

\subsection{Théorème}
{\sl
Soit $E$ un espace vectoriel complexe et soit $\lambda$ une forme différentielle à valeurs dans $E$ sur un ouvert $\calo$ de $\mathbb C^n$ qui est pure de type $(r,s)$ avec $s\geq 1$. Les conditions suivantes $\mathrm{(i)}$ et $\mathrm{(ii)}$ sont équivalentes :\\
$\mathrm{(i)}$ $\bar\partial \lambda=0$ sur $\calo$\\
$\mathrm{(ii)}$ On a un recouvrement ouvert $(\calo_\alpha)$ de $\calo$ et des formes différentielles à valeurs dans $E$ pures de type $(r,s-1)$ $\mu_\alpha $ sur les $\calo_\alpha$ tels que $\lambda\restriction \calo_\alpha=\bar\partial\mu_\alpha$ pour tout $\alpha$.
}\\

\noindent \underbar{Démonstration}. En remplaçant $E$ par $E\otimes (\stackrel{r}{\wedge}\mathbb C^n)$, il suffit de le démontrer pour $r=0$. Le cas $s=1$ résulte de 2.5 en identifiant $E$ à une algèbre de Lie abélienne, (on peut aussi le démontrer directement à partir de 2.1) ; dans tous les cas, on peut se ramener au cas de formes scalaires de type $(o,s)$ car la tensorisation avec $E$ est triviale. Le résultat est vrai, d'après ce qui précède, si les $d\bar z_2,\dots,d\bar z_n$ n'interviennent pas dans $\lambda$. Supposons-le vrai si 
les $d\bar{z}_m,d\bar{z}_{m+1},\dots,d\bar{z}_n$ n'interviennent pas et montrons qu'il est alors vrai pour un $\lambda$ où les $d \bar{z}_{m+1},\dots, d \bar z_n$ n'interviennent pas i.e. on a $\lambda=\varphi\wedge d \bar z_m+\psi$ où $\varphi$ et $\psi$ ne font pas intervenir $d\bar z_m,\dots, d\bar z_n$. D'après 2.1, on peut toujours résoudre localement $\varphi=\partial P/\partial \bar z_m\Rightarrow$ localement $\lambda =\bar\partial q+\psi'$ où $\psi'$ ne fait intervenir que $d \bar z_1,\dots, d\bar z_{m-1}$ et $\bar \partial\psi'=0$.\\
L'hypothèse de récurrence implique alors localement $\psi'=\bar\partial \chi$ et par consé\-quent localement $\Lambda=\bar\partial\mu$ ce qui est équivalent au résultat désiré.~$\square$

\subsection{Proposition}
{\sl
Soit $\lambda$ une forme différentielle fermée pure de type $(r,s)$ avec $r\geq 1$ et $s\geq 1$ sur un ouvert $\calo$ de $\mathbb C^n$. Alors on a un recouvrement ouvert $(\calo_\alpha)$ de $\calo$ et des formes $\varphi_\alpha$ pures de type $(r-1,s-1)$ sur les $\calo_\alpha$ tels que $\lambda\restriction \calo_\alpha=\partial\bar\partial\varphi_\alpha$ pour tout $\alpha$.}\\

\noindent\underbar{Démonstration}.
$d\lambda=0\Leftrightarrow \partial\lambda=0$ et $\bar\partial\lambda=0$ pour une forme pure $\lambda$. Localement on peut résoudre $d\lambda=0$ par $\lambda=d\mu$, $\mu=\sum_{k+\ell=r+s-1}\mu^{k,\ell}$ et comme $\lambda$ est de type $(r,s)$, $\mu=\mu^{r-1,s}+\mu^{r,s-1}$ avec $\bar\partial\mu^{r-1,s}=0$ et $\partial\mu^{r,s-1}=0\Rightarrow$ d'après 2.6, localement $\mu^{r-1,s}=\bar\partial\varphi_1$, $\mu^{r,s-1}=\partial \varphi_2\Rightarrow \lambda=\partial\bar\partial\varphi_1+\bar\partial\partial\varphi_2=\partial\bar\partial(\varphi_1-\varphi_2)$.~$\square$

\newpage

\begin{center}
{\large Références pour le chapitre 2}
\end{center}
\vspace{1cm}
\begin{itemize}
\item 
Hörmander L. : Complex analysis in several variables. Van Nostrand, 1966.\\

\item
Malgrange B. : Lectures on the theory of functions of several complex variables. Tata Institute, 1958.
\end{itemize}

\newpage

\section{Variétés complexes et presque complexes}

\subsection{Variétés complexes}

Rappelons, (voir Chapitre 1), qu'une {\sl variété complexe de dimension complexe} $n$ est un espace topologique séparé muni d'un atlas complet compatible avec le pseudo-groupe $\Gamma(\mathbb C^n)$ des transformations holomorphes de $\mathbb C^n$.\\
Soit $\omega$ une forme différentielle définie dans un ouvert $\calo$ d'une variété complexe $V$ de dimension $n$; si $x\in \calo$ et si $(U,\varphi)$ est une carte de l'atlas de $V$ (comme variété complexe) avec $x\in U$, $\varphi^{-1\ast}(\omega)$ est une forme définie dans un voisinage de $\varphi(x)$ dans $\mathbb C^n$. On a donc pour $\varphi^{-1\ast}(\omega)$ une notion de composante pure de type $(p,q)$ et des opérateurs $\partial$ et $\bar\partial$, comme dans le chapitre 2. Les images inverses par $\varphi$ de la composante pure de type $(p,q)$, $P^{p,q}\varphi^{-1\ast}(\omega)$, de $\varphi^{-1\ast}(\omega)$, de $\partial\varphi^{-1\ast}(\omega)$ et de $\partial\varphi^{-1\ast}(\omega)$ sont des formes différentielles définies au voisinage de $x$ dans $V$; leurs valeurs en $x$ ne dépendent pas de la carte $(U,\varphi)$ choisie avec $x\in U$ (en vertu de la compatibilité avec $\Gamma(\mathbb C^n)$). On a donc, pour une forme différentielle $\omega$ sur un ouvert $\calo$ de $V$, une notion {\sl de composante pure de type} $(r,s)$, $P^{r,s}\omega$, et des {\sl opérateurs} $\partial$ et $\bar\partial$ avec $\partial \omega=\sum_{r,s}P^{r+1,s}dP^{r,s}\omega$ et $d=\partial+\bar\partial$ qui correspondent aux notions définies sur $\mathbb C^n$ dans le chapitre 2.\\
Une variété complexe de dimension $n$ est en particulier une variété différentia\-ble de dimension $2n$ si on la munit de l'atlas complet compatible avec $\Gamma(\mathbb R^{2n})$).\\
Nous parlerons de {\sl cartes différentiables} pour désigner les cartes de cet atlas  et de{\sl cartes holomorphes} pour désigner es cartes de l'atlas compatibles avec $\Gamma(\mathbb C^n)=\Gamma_{\GL(n,\mathbb C)}(\mathbb R^{2n})$ définissant la structure de variété complexe de $V$.\\
Soit $x\in V$ et $(U,\varphi)$ une carte holomorphe avec $x\in U$; l'application $\varphi_{x^\ast}$ tangente en $x$ à $\varphi$ est un isomorphisme de $T_x(V)$ sur $\mathbb R^{2n}=\mathbb C^n$. \`A la multiplication par $i$ dans $\mathbb C^n$, i.e. l'endomorphisme $\left(\begin{array}{cc} 0 & -1_n\\ 1_n & 0\end{array}\right)$ de $\mathbb R^{2n}$ correspond par $\varphi_{x^\ast}$, l'endomorphisme $J_x=\varphi^{-1}_{x^\ast}(i \varphi_{x^\ast})$ de $T_x(V)$ satisfaisant $J_x^2=-1_x$. Il résulte de la compatibilité avec $\Gamma(\mathbb C^n)$ que $J_x$ ne dépend pas de la carte holomorphe $(U,\varphi)$ choisie avec $x\in U$. On a donc un champ (différentiable) $x\mapsto J_x$ d'endomorphismes des $T_x(V)$ satisfaisant $J^2=-1$ correspondant à la multiplication par $i$ dans $\mathbb C^n$.

\subsection{Structures presque complexes}

Soit $V$ une variété différentiable {\sl une structure presque complexe sur} $V$ est un champ différentiable, $J$, d'endomorphismes des espaces tangents à $V$ dont le carré est moins l'identité $J^2=-1$. Une variété munie d'une structure presque complexe est appelée {\sl variété presque complexe}. Il est évident qu'une variété presque complexe est de dimension paire (sur $\mathbb R$).\\
Si $V$ est une variété complexe de dimension complexe $n$, la structure presque complexe $J$ construite ci-dessus (correspondant à la multiplication par $i$ dans $\mathbb C^n$) sera appelée {\sl sa structure presque complexe canonique}. Si $(V,J)$ est la variété presque complexe sous-jacente à une variété complexe de dimension complexe $n$ ($\Rightarrow \dim_{\mathbb R}(V)=2n$), une carte différentiable $(U,\varphi)$ de $V$ est une carte holomorphe si et seulement si $\varphi_\ast \circ J\circ \varphi^{-1}_{\ast}=i=\left(\begin{array}{cc} 0 & -1_n\\ 1_n & 0\end{array}\right)$; autrement dit, la structure de variété complexe de $V$ est complètement déterminée par la structure presque complexe sous-jacente. On peut donc considérer que {\sl les variétés complexes sont des variétés presque complexes particulières}. En fait, ceci est un cas particulier de la situation décrite dans le chapitre 1. Se donner une structure presque complexe $J$ sur une variété $V$ de dimension $2n$ est équivalent à se donner une $\GL(n,\mathbb C)$-structure $P$ sur $V$ ($\GL(n,\mathbb C)\subset \GL(2n,\mathbb R)$); $P$ est l'ensemble des repères tangents à $V$ dans lesquels $J$ est représenté par la matrice $\left(\begin{array}{cc} 0  &-1_n\\ 1_n & 0\end{array}\right)$. $J$ est la structure presque complexe canonique d'une structure de variété complexe sur $V$ si et seulement si $P$ est intégrable.

\subsection{Connexions presque complexes}

Soit $(V,J)$ une variété presque complexe. {\sl Une connexion presque complexe sur $(V,J)$} est une connexion linéaire $\nabla$ sur $V$ telle que $\nabla J=0$; il revient au même de dire que c'est une connexion sur $\GL(n,\mathbb C)$-structure $P$ correspondante $(\dim(V)=2n)$.\\
La condition nécessaire et suffisante à l'intégrabilité d'une $\GL(n,\mathbb C)$-structure $P$ sur $V$ ($\dim V=2n$) étant $c_0(P)=0$, (voir chapitre 1), nous pouvons énoncer le théorème suivant.

\subsection{Théorème}

{\sl
Une variété presque complexe $(V,J)$ est (sous-jacente à) une variété complexe si et seulement si elle admet une connexion presque complexe sans torsion.
}

Dans ce cas, nous dirons, par abus de langage, que {\sl la structure presque complexe $J$ est intégrable}. Il est commode de caractériser l'intégrabilité de $J$ par l'annulation d'un tenseur construit directement à partir de $J$; pour cela, nous utiliserons la proposition suivante.

\subsection{Proposition}
{\sl Soit $(V,J)$  une variété presque complexe et soient $R$ et $T$ la courbure et la torsion d'une connexion presque complexe sur  $(V,J)$. On a :
\[
[J,R(X,Y)]=0
\]
et 
\[
T(JX,JY)-JT(JX,Y)-JT(X,JY)-T(X,Y)=-\frac{1}{2}N_J(X,Y)
\]
pour toute paire $X,Y$ de champs de vecteurs, où $N$ est le champ tensoriel (1-2) (ou plutôt la $2$-forme à valeurs vectorielles) défini par}
\[
N_J(X,Y)=2\{[JX,JY]-[X,Y]-J[X,JY]-J[JX,Y]\}
\]

\noindent\underbar{Démonstration}. C'est une conséquence immédiate de $\nabla J=0$ et des définitions de la courbure et de la torsion d'une connexion linéaire:
\[
R(X,Y)=[\nabla_X,\nabla_Y]-\nabla_{[X,Y]}
\]
et
\[
T(X,Y)=\nabla_X Y-\nabla_YX-[X,Y].~\square
\]
On remarquera que $N_J$ définit bien un champ tensoriel d'après la première partie de la proposition ($J$ est un champ tensoriel); ceci peut être vérifié directement. $N_J$ est appelé {\sl torsion de la structure presque complexe $J$ ou tenseur de Nijenhuis} de $J$. Cette terminologie est justifiée par la proposition suivante.

\subsection{Proposition}

{\sl Toute variété presque complexe $(V,J)$ admet une connexion presque complexe dont la torsion est $\frac{1}{8}N_J$.}\\

\noindent\underbar{Démonstration}. Soit $\tilde\nabla$ une connexion linéaire sans torsion sur $V$ et soit $Q$ le champ tensoriel défini par $4Q(X,Y)=(\tilde \nabla_{JY}J)X+J(\tilde\nabla_YJ)X+2J(\tilde\nabla_XJ)Y$ pour toute paire $X,Y$ de champs vectoriels; on définit une connexion $\nabla$ par :
$\nabla_XY=\tilde\nabla_XY-Q(X,Y)$. $\nabla$ est presque complexe et sa torsion est $\frac{1}{8}N_J$. $\square$\\
En combinant 3.5, 3.6 et 3.4 on a le théorème suivant.

\subsection{Théorème}

{\sl Une variété presque complexe $(V,J)$ est (sous-jacente à) une variété complexe si et seulement si $N_J=0$.}\\
Autrement dit $J$ est intégrable si et seulement si $N_J=0$.\\
Soit $(V,J)$ une variété presque complexe de dimension $2n$, on peut munir l'espace tangent $T_x(V)$ à $V$ en $X$ d'une structure d'espace vectoriel complexe de dimension $n$ en posant
\[
zV=xV+yJ_xV,\ \ \ \forall V\in T_x(V)\ \ \ \text{et}\ \ \forall z=x+iy\in \mathbb C
\]
Le fibré tangent $T(V)$ d'une variété presque complexe $(V,J)$ est ainsi canoniquement un fibré vectoriel complexe de rang $n$ si $\dim(V)=2n$.\\
On a donc, par dualité, pour les formes linéaires à valeurs complexes une notion de forme complexe-linéaire; plus généralement, si $\omega$ est une forme différentielle (à valeurs complexes) sur $V$, on a une notion de composantes pures $P^{r,s}\omega$ de type $(r,s)$ de $\omega$ ($r$-complexe-linéaire et $s$-complexe-antilinéaire) qui généralise la notion déjà introduite pour les variétés complexes. On peut définir sur les formes différentielles les opérateurs $\partial$ et $\bar\partial$ par
$\partial=\sum_{r,s}P^{r+1,s}\circ d\circ P^{r,s}$ et $\bar\partial=\sum_{r,s}P^{r,s+1}\circ d\circ P^{r,s}$; on a bien $\bar\partial \omega=\overline{(\partial\bar\omega)}$ mais généralement $d\not=\partial+\bar\partial$ et $\partial^2\not=0$. Une autre caractérisation des variétés complexes est, en fait la suivante.

\subsection{Théorème}. {\sl Une variété presque complexe $(V,J)$ est une variété complexe si et seulement si l'une des propriétés suivantes (qui sont équivalentes) est satisfaite.\\
a) $\Im(d\circ P^{1,0})\subset\Im(P^{2,0}) + \Im(P^{1,1})$\\
b) $\Im(d\circ P^{r,s})\subset \Im(P^{r+1,s})+ \Im(P^{r,s+1})$\\
c) $d=\partial +\bar\partial$\\
e) $\partial^2=0$\\
(ou l'une des propriétés obtenue par conjugaison complexe)}

\subsection{Espaces projectifs complexes}
Soit $P_n(\mathbb C)$ l'ensemble des sous-espaces vectoriels complexes de dimension un de $\mathbb C^{n+1}$. $P_n(\mathbb C)$ est le quotient de $\mathbb C^{n+1}\backslash \{0\}$ par la relation d'équivalence $\sim$ : $v\sim v'$ dans $\mathbb C^{n+1}\backslash \{0\}$ si $\exists \lambda \in \mathbb C\backslash \{0\}$ tel que $v=\lambda v'$. Nous munirons $P_n(\mathbb C)$ de la topologie quotient; c'est un espace topologique séparé. Soit $P\in P_n(\mathbb C)$ et $(z_0,z_1,\dots,z_n)\in \mathbb C^{n+1}\backslash \{0\}$ un élément de $P$ non nul; ses composantes $z_k$, $k\in \{0,1,\dots, n\}$ sont appelées {\sl coordonnées homogènes de} $P\in P_n(\mathbb C)$, elles sont déterminées par $P$ à un facteur complexe différent de zéro mutiplicatif près (indépendant de $k$). Soit $\calo_k$ l'ensemble ouvert dans $P_n(\mathbb C)$ des sous-espaces vectoriels complexes de dimension un qui ne sont pas contenus dans l'hyperplan $z_k=0$ et, pour $P\in \calo_k$, soit $\varphi_k(P)=(z_0/z_k, \dots,z_{k-1}/z_k,z_{k+1}/z_k,z_n/z_k)\in \mathbb C^n$,où $(z_0,\dots, z_n)\in \mathbb C^{n+1}\backslash\{0\}$ est un système de coordonnées homogènes de $P$ (i.e. un élément non nul de $\mathbb C^{n+1}$ appartenant à $P\subset\mathbb C^{n+1})$. $\varphi_k$ est un homéomorphisme de $\calo_k$ sur $\mathbb C^n$ (indépen\-dant du choix de $(z_0,\dots,z_n)\in P\backslash \{0\}\subset \mathbb C^{n+1})$ et on a $\cup^{k=n}_{k=0} \calo_k=P_n(\mathbb C)$. $\varphi_r\circ \varphi^{-1}_s$ est holomorphe de $\varphi_s(\calo_r\cap \calo_s)\subset \mathbb C^n$ sur $\varphi_r(\calo_r\cap \calo_s)$, (multiplication par $z_s/z_r$ à une réindexation près), ainsi que son inverse. Autrement dit, $(\calo_k,\varphi_k)_{k\in\{0,1,\dots,n\}}$ est un atlas de $P_n(\mathbb C)$ compatible avec $\Gamma(\mathbb C^n)$. En complétant cet atlas compatible avec $\Gamma(\mathbb C^n)$, on obtient une structure de variété complexe de dimension complexe $n$ sur $P_n(\mathbb C)$; $P_n(\mathbb C)$ muni de cette structure est {\sl l'espace projectif complexe de dimension $n$.}\\
On peut aussi identifier $P_n(\mathbb C)$ à l'ensemble $\{P\in \End(\mathbb C^{n+1})\vert P^2=P=P^+$ et $\Tr(P)=1\}$, des projecteurs hermitiens de rang un dans $\mathbb C^{n+1}$. La structure (sous-jacente) de variété différentiable de $P_n(\mathbb C)$ est induite par la structure d'espace vectoriel réel des endomorphismes hermitiens $\End^h(\mathbb C^{n+1})$ de $\mathbb C^{n+1}$; $P_n(\mathbb C)$ est ainsi une sous-variété différentiable de $\End^h(\mathbb C^{n+1})=\mathbb R^{(n+1)^2}$. Par différentiation, l'espace tangent $T_P(P_n(\mathbb C))$ en $P$ à $P_n(\mathbb C)$ s'identifie à l'espace (vectoriel réel) des endomorphismes hermitiens de trace nul $A$ de $\mathbb C^{n+1}$ satisfaisant à $PA+AP=A$, ou ce qui revient au même à $[P,[P,A]]=A$~:
\[
T_p(P_n(\mathbb C))=\{A\in \End(\mathbb C^{n+1})\vert A=A^+=[P,[P,A]], \Tr(A)=0\}
\]
L'application $A\mapsto J_P(A)=-i[P,A]$ définit alors un endomorphisme de $T_P(P_n(\mathbb C))$ dont le carré est moins l'identité; {\bf la structure presque complexe, $P\mapsto J_P$, ainsi définie sur $P_n(\mathbb C)$ n'est autre que la structure presque complexe canonique de la variété complexe} $P_n(\mathbb C)$. En effet, si $Z$ est une matrice colonne $(n+1)\times 1$ dont les éléments forment un système de coordonnées homogènes de $P$, on peut écrire $P=\frac{1}{\vert\vert Z\vert\vert^2}ZZ^+$; par différentiation, un vecteur tangent en $P$ à $P_n(\mathbb C)$ et une matrice $A$ de la forme $A=\left(\frac{1}{\vert\vert Z \vert\vert^2})^2(\vert\vert Z\vert\vert^2\{ZZ^{\prime +}+Z'Z^+\}-\{\langle Z'\vert Z\rangle+\langle Z\vert Z'\rangle\}ZZ^+\right)$ où $Z'$ est un vecteur de $\mathbb C^{n+1}$. La multiplication par $i$ de $Z'$ (et par $-i$ de $Z^{\prime +}$ correspond à la structure presque complexe canonique de $P_n(\mathbb C)$ en $P$ et est donné comme on le vérifie immédiatement par $A\mapsto -i[P,A]$.

\subsection{Grassmanniennes complexes}
Soit $G_{N,p}(\mathbb C)$ l'ensemble des sous-espaces vectoriels complexes de dimension $p$ de $\mathbb C^N$. L'action de $U(N)$ sur $G_{N,p}(\mathbb C)$ est transitive et le stabilisateur de $\mathbb C^p(+\{O_{N-p}\})\subset \mathbb C^N$ est $U(p)\times U(N-p)$ de sorte que $G_{N,p}(\mathbb C)\simeq U(N)/ U(p)\times U(N-p)$. On peut en généralisant 3.9 construire un atlas de $G_{N,p}(\mathbb C)$ compatible avec $\Gamma(\mathbb C^{p(N-p)})$, (c'est un peu plus laborieux (voir 3.18 b). $G_{N,p}(\mathbb C)$ est ainsi une variété complexe de dimension complexe $p(N-p)$ {\sl appelée grassmannienne des $p$-plans de $\mathbb C^N$}.\\
Soit $P\in G_{N,p}(\mathbb C)$ et $(e_1,\dots,e_p)$ une base de $P$ $(e_k\in \mathbb C^N)$, $e_1\wedge e_2\wedge \dots \wedge e_p\in \wedge^p\mathbb C^N=\mathbb C^{C^p_N}$ représente complètement $P$ et $P$ détermine cet élément de $\stackrel{p}{\wedge}\mathbb C^N\backslash \{0\}$ à un nombre complexe multiplicatif non nul près. On a donc une injection $j:G_{N,p}(\mathbb C)\rightarrow P_{C^p_N-1}(\mathbb C)$ de $G_{N,p}(\mathbb C)$ dans l'espace projectif complexe de dimension $C^p_N-1$; cet injection permet en fait d'identifier $G_{N,p}(\mathbb C)$ à une sous-variété complexe de $P_{C^p_N-1}(\mathbb C)$. Si $P\in G_{N,p}(\mathbb C)$, un système de  coordonnées homogènes de $j(P)$ est appelé système de {\sl coordonnées Plücker} de $P$ (i.e. $(e_1\wedge \dots e_p)_{i_1\dots i_p}$; $i_1> i_2>\dots >i_p$). La condition pour qu'un élément de $\stackrel{p}{\wedge}\mathbb C^N$ représente un élément de $G_{N,p}(\mathbb C)$ est qu'il soit décomposable, ce qui correspond à un certain nombre de conditions homogènes quadratiques dans $\stackrel{p}{\wedge}\mathbb C^N=\mathbb C^{C^p_N}$.\\
On peut identifier $G_{N_p}(\mathbb C)$ à l'ensemble, $\{P\in \End(\mathbb C^N)\vert P=P^+=P^2$, $\Tr(P)=p\}$, des projecteurs hermitiens de rang $p$ dans $\mathbb C^N$. L'espace\linebreak[4] $T_p(G_{N,p}\mathbb C))$ s'identifie alors, comme dans 3.9, à 
\[
T_p(G_{N,p}(\mathbb C))=\{A\in \End(\mathbb C^N)\vert A=A^+=[P,[P,A]],\Tr(A)=0\}
\]
et $A\mapsto J_p(A)=-i[P,A]$ est la structure presque complexe canonique de la variété complexe $G_{N,p}(\mathbb C)$.

\subsection{Remarque} $G_{n+1,1}(\mathbb C)=P_n(\mathbb C)$ et $\vec u\mapsto \frac{1-\sigma(\vec u)}{2}$ est un difféomorphisme de 
$$S^2=\{\vec u\in \mathbb R^3\vert \vec u^2=1\}$$
 sur 
\[
P_1(\mathbb C)=\{P\in \End (\mathbb C^2)\vert P=P^+=P^2,\Tr(P)=1\} ;
\]
$\sigma(\vec u)=\sum^3_1\sigma_k u_k$ où $\sigma_1,\sigma_2,\sigma_3$ sont les matrices de Pauli\\
 $P_1(\mathbb C) (\simeq S^2)$ est la {\sl sphère de Riemann}.

\subsection{Variétés de drapeaux complexes}
On généralise 3.9 et 3.10 en considérant, pour $n_0\geq n_1\geq \dots \geq n_k$ entiers, l'ensemble $G_{n_0,n_1,\dots,n_k}(\mathbb C)$ des suites de sous-espaces vectoriels complexes de $\mathbb C^{n_0}$ emboitée $\mathbb C^{n_0} \supset E_1\supset E_2 \supset \dots \supset E_k$ avec $\dim(E_\ell)=n_\ell$. On peut munir $G_{n_0,\dots,n_k}(\mathbb C)$ d'une structure de variétés complexes. Les $G_{n_0,\dots,n_k}(\mathbb C)$ s'appellent {\sl variétés de drapeaux}. En ``oubliant" certains des entiers $n_\ell$ ($\ell\geq 1$), on a des projections $G_{n_0,n_1,\dots,n_k}(\mathbb C)\rightarrow G_{n_0,n_{i_1},\dots,n_{i_p}}(\mathbb C)$ ($n_{i_1}\geq \dots \geq n_{i_p},n_{i_\ell}\in \{n_1,n_2,\dots,n_k\}$); ces projections sont holomorphes. L'application $(E_1,\dots,E_k)\in G_{{n_0},n_1,\dots,n_k}(\mathbb C)\mapsto E_1\oplus E_2\oplus\dots\oplus E_k\subset \mathbb C^{kn_0}$ permet d'identifier $G_{n_0,n_1,\dots,n_k}(\mathbb C)$ à un sous-ensemble  de l'espace projectif complexe de dimension $C^{n_1+\dots+n_k}_{kn_0}-1$. Ces variétés de drapeaux sont à nouveau des sous-variétés complexes des espaces projectifs complexes caractérisées, en termes de coordonnées homogènes par des conditions algébriques homogènes, (annulation de polynômes homogènes).

\subsection{Structure presque complexe de $T(V)$}
Soit $V$ une variété différentiable de dimension $n$, $T(V)$ son fibré tangent et $\nabla$ une connexion linéaire sur $V$. Soit $\pi:T(V)\rightarrow V$ la projection; si $v\in T(V)$, l'ensemble des vecteurs tangents en $v$ à $T(V)$ se projetant par $\pi_\ast$ sur $0\in T_{\pi(\sigma)}(v)$ (vecteurs verticaux), peut être identifié à $T_{\pi(v)}(V)$ et $\pi_\ast$ définit d'autre part un isomorphisme de l'espace des vecteurs horizontaux pour $\nabla$ en $v$ sur $T_{\pi(v)}(V)$. On a donc un isomorphisme $\varphi_v$ de $T_v(T(V))$ sur $T_{\pi(v)}(V)\oplus T_{\pi(v)}(V)$. L'application 
$(X,Y)\mapsto (-Y,X)$ de $T_{\pi(v)}(V)\oplus T_{\pi(v)}(V)$ sur lui-même permet de définir, via l'isomorphisme $\varphi_v$ un endomorphisme $J_v$ de $T_v(T(V))$ satisfaisant $J^2_v=-\bbbone$. $v\mapsto J_v$ est une structure presque complexe qui n'est intégrable que si la courbure et la torsion de $\nabla$ sont nulles (comme on le voit en calculant explicitement $N_J$; nous ferons ultérieurement un calcul analogue.

\subsection{Variétés (presque) hermitiennes}
Une {\sl variété presque hermitienne} est une variété riemannienne, $V$, munie d'une structure presque complexe isométrique $J$. Dans le cas où $J$ est intégrable ($V$ est complexe) $V$ est appelée variété {\sl hermitienne}.\\
Une variété presque hermitienne de dimension $2n$ n'est autre qu'une variété de dimension $2n$ munie d'une $U(n)$-structure $(U(n)\subset \GL(2n,\mathbb R))$ mais on remarquera que la {\sl $U(n)$-structure d'une variété hermitienne n'est intégrable que si sa structure riemannienne est plate.}\\
Soit $V$ une variété presque hermitienne de métrique $g$ et soit $J$ sa structure presque complexe ; alors $(X,Y)\mapsto g(X,JY)=-g(Y,JX)$ est une 2-forme différentielle qui est pure de type (1,1) on l'appelle {\sl la 2-forme fondamentale de la variété presque hermitienne}; nous la noterons $\psi$. L'espace tangent en $x\in V$ à $V$ est canoniquement, comme nous l'avons vu dans 3.7, un espace vectoriel complexe de dimension $n$ (si $\dim_\mathbb R(V)=2n)$; on le munit d'une structure d'espace hilbertien en posant :
\[
\langle X\vert Y\rangle=g(X,Y)-ig(X,JY)=g(X,Y)-i\psi(X,Y)
\]
$\forall X,Y\in T_x(V)$, (on notera que $\vert\vert X\vert\vert^2=g(X,X)=\langle X \vert X\rangle$). Le fibré tangent d'une variété presque hermitienne est donc canoniquement un fibré en espaces hilbertiens (={\sl fibré hermitien}) de rang $n$ si $\dim_\mathbb RV=2n$. Il en est de même du fibré cotangent $T^\ast(V)$ et du fibré en algèbres extérieures $\wedge T^\ast(V)=\oplus^{2n}_{p=0}\stackrel{p}{\wedge}T^\ast (V)$ et la décomposition $\wedge T^\ast(V)=\oplus^n_{r,s=0} \wedge^{r,s}T^\ast(V)$ en composantes pures devient une décomposition orthogonale.\\
Soit $V$ une variété presque hermitienne, en général sa connexion riemannienne $\nabla$ n'est pas presque complexe, (i.e. $\nabla J\not= 0$), et l'on sait que si elle l'est alors $V$ est complexe en vertu de 3.4 puisque $\nabla$ est sans torsion. On a le théorème général suivant. 

\subsection{Théorème}

Soit $V$ {\sl une variété presque hermitienne. Les conditions $a$ et $b$ suivantes sont équivalentes}.\\
a) {\sl La connexion riemannienne de $V$ est une connexion presque complexe, $(\nabla J=0)$.}\\
b) {\sl La structure presque complexe de $V$ est intégrable et sa 2-forme fondamentale est fermée, $(d\psi=0)$}\\

\noindent \underbar{Démonstration}. Supposons $\nabla J=0$ alors, comme $\nabla$ est sans torsion, 3.4 implique $J$ est intégrable; on a $d\psi(X,Y,Z)=\nabla\psi(X,Y,Z)$ et comme $\psi(X,Y)=g(X,JY)$ et que $\nabla g=0$ et $\nabla J=0$ on a $\nabla\psi=d\psi=0$.\\
On a d'autre part en utilisant la définition de $\psi$, de $d\psi$ et de $N_J$ :
\[
4g((\nabla_XJ)Y,Z)=6d\psi(X,JY,JZ)-6d\psi(X,Y,Z)+g(N_J(Y,Z),JX)
\]
$\Rightarrow \nabla J=0$ si $d\psi=0$ et $N_J=0$.~$\square$

\subsection{Variétés kähleriennes}

Une variété satisfaisant aux conditions de 3.15 est appelée {\sl variété Kählerienne}. Une variété presque hermitienne dont la 2-forme fondamentale est fermée sera appelée {\sl variété presque kählerienne}.\\
Une variété kählerienne $V$ est donc en particulier une variété complexe. Sa 2-forme fondamentale $\psi$ est fermée pure de type (1,1); on peut donc localement lui appliquer la proposition 2.7 et obtenir la représentation $\psi=i\partial\bar \partial F_\alpha$ sur $\calo_\alpha$ où $(\calo_\alpha)_\alpha$ est un recouvrement ouvert de $V$ (par des domaines de cartes holomorphes) et $F_\alpha:\calo_\alpha\rightarrow \mathbb R$ sont des fonctions différentiables à valeurs réelles (pour que $\psi$ le soit) sur les $\calo_\alpha$ appelées {\sl potentiels de Kähler}. La métrique de $V$ est évidemment complètement déterminée par $\psi$ donc par les potentiels de Kähler. Si $F$ est définie sur le domaine $\calo$ de la carte holomorphe $(\calo,\varphi)$ avec $\psi=i\partial\bar\partial F$ et si les $z^k=\varphi^k$ désignent les coordonnées complexes correspondantes dans $\mathbb C^n$, on a : $\psi=i(\partial_k\bar\partial_\ell F)dz^k\wedge \overline{dz^\ell}$ et $G=(\bar\partial_k\partial_\ell F)\overline{dz^k}\otimes dz^\ell$ est la forme sesquilinéaire représentant la métrique hermitienne de $V$; la matrice hermitienne $(\bar\partial_k\partial_\ell F)$ doit donc être strictement positive en chaque point si $F$ est un potentiel de Kähler.\\

\subsection{Exemple : la grassmannienne $G_{N,p}(\mathbb C)$}

Identifions encore la grassmannienne $G_{N,p}(\mathbb C)$ avec l'ensemble des projecteurs $P$ hermitiens de rang $p$ dans $\mathbb C^N$, (voir 3.10); $G_{N,p}(\mathbb C)$ est ainsi identifiée comme variété différentiable à une sous-variété de l'espace vectoriel réel $\End^h(\mathbb C^N)\simeq \mathbb R^{N^2}$ des endomorphismes hermitiens de $\mathbb C^N$. $\End^h(\mathbb C^N)$ est un espace euclidien pour le produit scalaire $(A,B)=\Tr(AB)$ et la métrique riemannienne induite sur $G_{N,p}(\mathbb C)$ s'écrit symboliquement $ds^2=\Tr((dP)^2)$. On a $dP\circ J=-i[P,dP]$, d'après 3.10, il en résulte $\Tr((dP\circ J)^2)=-\Tr([P,dP]^2)=\Tr(dP[P[PdP]])=\Tr((dP)^2)$ (car $dP=[P[P,dP]])$; munie de cette métrique $G_{N,p}(\mathbb C)$est donc une variété hermitienne. Sa 2-forme fondamentale est $\psi=\Tr(dP[P,dP])=2(\Tr(PdP\wedge dP)$. On a par différentiation extérieure $d\psi=2i\Tr(dP\wedge d P\wedge dP)=0$ car $dPP+PdP=dP\Leftrightarrow (dP)P=(1-P)dP$ et, par conséquent, la trace du produit d'un nombre impair de dérivées de $P$ est nulle. $G_{N,p}(\mathbb C)$ est donc une variété Kählérienne.

\subsection{Remarques}

a) $A\mapsto [P[P,A]]$ est, pour $P\in G_{N,p}(\mathbb C)$ un opérateur symétrique dans l'espace euclidien $\End^h(\mathbb C^N)$; on vérifie facilement que c'est la projection orthogonale de $\End^h(\mathbb C^N)$ sur l'espace tangent en $P$ à $G_{N,p}(\mathbb C)$. Cette remarque est très utile pour les calculs de dérivées covariantes sur $G_{N,p}(\mathbb C)$.\\

\noindent b) Soit $\calo_{1,2,\dots,p}$ l'ouvert de $G_{N,p}(\mathbb C)$ défini de la manière suivante : $\calo_{1,2,\dots,p}$ est l'ensemble des sous-espaces vectoriels de dimension $p$ de $\mathbb C^N$ qui ne contiennent aucun vecteur non nul orthogonal à $\mathbb C^p=\{z^1,\dots, z^p\}\subset \mathbb C^N=\mathbb C^p\oplus \mathbb C^{N-p}$. En termes de projecteurs hermitiens de rang $p$,\\
 $\calo_{1,2,\dots,p}=\{P\in G_{N,p}(\mathbb C)\vert PP_{1,2,\dots,p}\}$ est de rang $p\}$,\\
où $P_{1,2,\dots,p}=\left( \begin{array}{cc} 1_p & 0\\ 0 & 0\end{array}\right)$ est le projecteur sur $\mathbb C^p(\subset \mathbb C^N)$. Un sous-espace de $\calo_{1,2,\dots,p}$ admet une base unique de la forme\\
 $((1,0,\dots,0,z^p_1,\dots,z^N_1),\dots,(0,0,\dots,1, z^p_p,\dots,z^N_p))$ et les $z^\ell_k$, ($k\in \{1,\dots,p\},$ $\ell\in \{p,\dots,N\})$ forment un système de coordonnées holomorphes de $G_{N,p}(\mathbb C)$ sur $\calo_{1,\dots,p}$. Soit $Z=\left(\begin{array}{ccc}z^p_1 &\dots & z^p_p\\ \vdots & & \vdots\\
s^N_1 & \dots & z^N_p\end{array}\right)$ la matrice $(N-p)\times p$ correspondante, on a pour le projecteur $P$ correspondant :
\[
P=\left(\begin{array}{c} 1_p\\ Z\end{array}\right)(1_p+Z^+Z)^{-1}(1_p,Z^+),
\]
$1_p+Z^+Z$ est une matrice $p\times p$ hermitienne inversible. Un potentiel de Kähler sur $\calo_{1,\dots,p}$ est donné par $F_{1,\dots,p}=\log (\det (1_p+Z^+Z)]$.

\subsection{Structure presque kählerienne de $T(V)$}

Reprenons l'exemple 3.13 dans le cas où $V$ est une variété riemannienne et où $\nabla$ est la connexion riemannienne sur $V$. On a une structure riemannienne sur $V$. On a une structure riemannienne canonique sur $T(V)$ obtenue en relevant horizontalement la métrique de $V$, en munissant les vecteurs verticaux de leur métrique naturelle comme vecteurs
tangents à des espaces euclidiens (les $T_x(V)$) et en prenant la somme directe. La structure presque complexe construite en 3.13 sur $T(V)$ est, par définitions isométrique de sorte que $T(V)$ est une variété presque hermitienne. On vérifiera que la 2-forme fondamentale correspondante est fermée de sorte que le fibré tangent à une variété riemannienne $V$ est canoniquement une variété presque kählerienne qui n'est kählerienne que si $V$ est plate. {\bf La 2-forme fondamentale de $T(V)$ n'est autre que l'image inverse de la 2-forme symplectique canonique $\omega=d\theta$ des fibrés cotangents $T^\ast(V)$ par l'isomorphisme $f:T(V)\rightarrow T^\ast(V)$ induit par la métrique de $V,(f(X)=g(X,\bullet)\in T^\ast_x(V)$,  $\forall X\in T_x(V))$}.

\subsection{Le cas des variétés de dimension 2}

En dimension 2, une structure presque complexe est toujours intégrable et {\bf 
est une structure conforme orientée}. Si $(V,J)$ est presque complexe avec $\dim V=2$, pour tout champ de vecteur $X,(X,JX)$ est une base des espaces tangents aux points où $X\not= 0$ et $N_J(X,JX)=0$ comme on le vérifie dans le cas général $(\forall X)$. On a, d'autre part, $CO^+(2)=\GL(1,\mathbb C)\subset \GL(2,\mathbb R)$, donc une structure complexe est une structure conforme orientée qui est toujours intégrable en dimension 2, c'est-à-dire plate.\\
Une manière de se donner une structure conforme est de se donner une métrique riemannienne. En dimension 2, la structure conforme orientée, i.e. la structure complexe, sous-jacente est évidemment isométrique; une variété riemannienne orientée est donc, canoniquement, une variété hermitienne. La 2-forme fondamentale correspondante est toujours fermée puisqu'elle est de degré maximum. {\bf C'est donc une variété kählerienne}. L'atlas compatible avec $\Gamma(\mathbb C)$ d'une variété riemannienne de dimension 2, $V$ de métrique $g$, est constitué par les cartes $(U,z=x+iy)$ où $(U,(x,y))$ est une carte différentiable telle que $g=A(dx\otimes dx+dy\otimes dy)$; dans une telle carte, un potentiel Kähler est une solution de $\Delta F=4\frac{\partial^2 F}{\partial z\partial \bar z}=4A$.\\
En résumé, la notion de variété complexe de dimension complexe 1 est identique à celle de variété orientée munie d'une structure conforme et la notion de variété Kählerienne de dimension complexe 1 est identique à celle de variété riemannienne.

\subsection{Applications (presque) (anti-) holomorphes}

Soient $(V_1,J_1)$ et $(V_2,J_2)$ deux variétés presque complexes de dimension $2n_1$ et $2n_2$ respectivement. Une application différentiable $f:V_1\rightarrow V_2$ est appelée {\sl application presque complexe} si $J_2\circ f_\ast=f_\ast\circ J_1$ ($f_\ast:T(V_1)\rightarrow T(V_2))$); nous parlerons aussi d'application {\sl presque holomorphe} pour une telle application et d'application {\sl presque anti-holomorphe} dans le cas où $f$ est presque complexe de $(V_1,J_1)$ dans $(V_2,-J_2)$, (i.e. $f_\ast\circ J_1=-J_2\circ f_\ast$). Dans le cas où $V_1$ et $V_2$ sont des variétés complexes nous parlerons simplement {\sl d'applications holomorphes} (et {\sl anti-holomorphe}).\\

Supposons $V_1$ et $V_2$ complexes (de dimensions complexes $n_1$ et $n_2$) et soit $f:V_1\rightarrow V_2$ une application holomorphe. Soient $(U_1,\varphi_1)$ et $(U_2, \varphi_2)$ des cartes holomorphes de $V_1$ et $V_2$ telles que $f(U_1)\cap U_2\not= \emptyset$; alors la restriction\linebreak[4] $\varphi_2 \circ f \circ \varphi^{-1}_1\restriction \varphi_1\circ f^{-1}(U_2)$ est une fonction holomorphe sur l'ouvert $\varphi_1(f^{-1}(U_2))$ de $\mathbb C^{n_1}$ à valeurs dans $\mathbb C^{n_2}$.

\subsection{Lemme}

{\sl Soit $\calo$ un ouvert connexe du plan complexe $\mathbb C$ et soit $P:\calo\rightarrow P_n(\mathbb C)$ une application différentiable. Alors $P$ est holomorphe ou anti-holomorphe si et seulement si on a $\left(\frac{\partial P}{\partial z}\right)^2=0$ dans $\calo$}.\\
(Pour la matrice $(n+1)\times (n+1)$ complexe $\partial P/\partial z$; $P(z)$ est un projecteur hermitien de rang 1 dans $\mathbb C^{n+1}$).\\

\noindent\underbar{Démonstration}. D'après 3.9 et la définition 3.21, $P$ est holomorphe ou anti-holomorphe si et seulement si on a : $[P,\frac{\partial P}{\partial z}]=\mp \frac{\partial P}{\partial z}$ ($-$ correspond à l'holomorphie). Si $P$ est holomorphe ou anti-holomorphe, on a donc
\[
\left(\frac{\partial P}{\partial z}\right)^2=\left([P,\frac{\partial P}{\partial z}]\right)^2=P\frac{\partial P}{\partial z} P\frac{\partial P}{\partial z}+\frac{\partial P}{\partial z} P\frac{\partial P}{\partial z}P-P\left(\frac{\partial P}{\partial z}\right)^2 -\frac{\partial P}{\partial z}P\frac{\partial P}{\partial z}
\]
$= -\left(\frac{\partial P}{\partial z}\right)^2=0$ car $P\frac{\partial P}{\partial z}=\frac{\partial P}{\partial z}(1-P)$ et $(1-P) \frac{\partial P}{\partial z}=\frac{\partial P}{\partial z}P.$\\
Inversement supposons que $\left(\frac{\partial P}{\partial z}\right)^2=0$. Soit $U_z$ un unitaire dans $\mathbb C^{n+1}$ tel que $U_z P(z)U_z^+=\left(\begin{array}{cc} 1 & 0\\ 0 & 0_p\end{array}\right)$ alors $U(z) \left(\frac{\partial P}{\partial z}(z)\right) U(z)^+=\left(\begin{array}{cc} 0 & \mu^t\\ \lambda & 0_p\end{array}\right)$ avec $\lambda, \mu\in \mathbb C^n$ car $\frac{\partial P}{\partial z}=P\frac{\partial P}{\partial z}(1-P)+(1-P)\frac{\partial P}{\partial z}P$ en différentiant $P^2=P$. On a donc :
\[
0=U(z)\left( \frac{\partial P}{\partial z}(z)\right)^2 U(z)^+=\left(\begin{array}{cc} \mu^t\lambda & 0\\ 0 & \lambda\mu^t\end{array}\right)\Rightarrow
\]
On a soit $\lambda=0$ soit $\mu=0$, mais
\[
U_z\left\{[P,\frac{\partial P}{\partial z}](z)+\frac{\partial P}{\partial z}(z)\right\}U^+_z=2\left(\begin{array}{cc} 0 & \mu^t\\ 0 & 0_p\end{array}\right)\]
 et
\[
U_z\left\{[P,\frac{\partial P}{\partial z}](z)-\frac{\partial P}{\partial z}(z)\right\}U^+_z=-2\left(\begin{array}{cc} 0 & 0\\ \lambda & 0_p\end{array}\right)\]
$\Rightarrow$ On a soit $[P,\frac{\partial P}{\partial z}]=-\frac{\partial P}{\partial z}$, (holomorphie), soit $[P,\frac{\partial P}{\partial z}]=+\frac{\partial P}{\partial z}$ (anti-holomorphie) en un point arbitraire de $\calo$ et, comme $\calo$ est connexe on a l'une de ces identités dans $\calo$.~$\square$\\

L'holomorphie où l'anti-holomorphie de $P$, $[P,\frac{\partial P}{\partial z}]=\mp\frac{\partial P}{\partial Z}$, peut aussi s'écrire $-iPdP=\pm \ast PdP$ où $\ast$ est la dualité de Hodge. Ce sont les solutions ``(i) self-duales ou anti-self-duales du modèle $\mathbb CP^n$" (=$P_n(\mathbb C))$.

\subsection{Applications harmoniques}

Soient $V_1$ et $V_2$ deux variétés riemanniennes orientées. Une application différentiable $f:V_1\rightarrow V_2$ est appelée {\sl application harmonique} si elle est solution des équations de Euler-Lagrange correspondant à $\int \vert\vert f_\ast\vert\vert^2 \vol_{V_1}$ où, en coordonnées locales, $\vert\vert f_\ast \vert\vert^2= g^{(2)}_{AB} (\partial_\mu f^A)(\partial _\nu f^B)g^{(1)\mu\nu}$ et\\
$\vol_{V_1}=\sqrt{\det (g^{(1)}_{\rho\sigma})} dx^1\wedge\dots \wedge dx^{n_1}$. 
\[
\Delta _1 f^A + \Gamma^{(2)A}_{BC} g^{(1)_{\mu\nu}} (\partial_\mu f^B)(\partial_\nu f^C)=0
\]

\subsection{Exemples}

a) $V_2=\mathbb R$; les applications harmoniques sont les fonctions harmoniques sur $V_1$.

b) $V_1=\mathbb R$; les applications harmoniques sont les géodésiques de $V_2$, (paramètres affine).

c) Soient $V_1$ et $V_2$ des variétés presque kähleriennes; alors les applications presque complexes de $V_1$ dans $V_2$ sont harmoniques. Cela est évidemment aussi vrai pour les applications presque anti-holomorphes, en changeant de signe une des structures presque complexes.

\subsection{Remarque : modèles $\sigma$}

Si $V_1$ est de dimension 2, il est clair, d'après la définition 3.23, que la notion d'application harmonique {\bf ne dépend que de la structure conforme de } $V_1$, i.e. comme on a supposé $V_1$ orientée {\bf que de sa structure complexe} canonique (voir 3.20); nous pouvons donc parler des {\sl applicatins harmoniques d'une variété complexe de dimension complexe 1 dans une variété riemannienne.} Ce sont les solutions du ``modèle $\sigma$ à valeurs dans $V_2$" sur $V_1$; dans le cas où $V_2$ est kählerienne, les applications holomorphes et anti-holomorphes sont harmoniques : ce sont les ``instantons du modèle".

\subsection{Lemme}

{\sl Soit $\calo$ un ouvert de $\mathbb R^m$, soit $X$ un champ de vecteurs à valeurs complexes sur $\calo$ et soit $P:\calo\rightarrow P_n(\mathbb C)$ une application différentiable. On a $(D_XP)^3=\frac{1}{2}\{\Tr(D_XP)^2\}D_XP$ où $D_X=\sum^n_{k=1}X^k\frac{\partial}{\partial x^k}$}.\\

\noindent\underbar{Démonstration}.  
$\forall x\in \calo$ on a :
\[
P(x)=U(x)^+\left(\begin{array}{cc}1 & 0\\ 0 & 0_n\end{array}\right) U(x),
\]
où $U(x)^+U(x)=\bbbone_{n+1}$.
Comme 
\[
P(x)(D_XP(x))P(x)=(1-P(x))(D_XP(x))(1-P(x))=0
\]
on a :
\[
D_X P(x)=U(x)^+\left(\begin{array}{cc} 0 & \lambda^t\\ \mu & 0_n\end{array}\right)U(x)
\]
 où $\mu,\lambda \in \mathbb C^n$. Il en résulte
\[
(D_XP(x))^3=\lambda^t\mu D_XP(x)=\frac{1}{2}\{\Tr(D_X P(x))^2\}D_XP(x).~\square
\]

\subsection{Proposition (modèle $P_n(\mathbb C)$)}

{\sl Soit $\calo$ un ouvert de $\mathbb C$ et $P:\calo\rightarrow P_n(\mathbb C)$ une application harmonique, (i.e. une solution du ``modèle $P_n(\mathbb C)$"). On a : $\frac{\partial}{\partial \bar z} \Tr \left(\frac{\partial P}{\partial z}\right)^2=0$ dans $\calo$}.\\

\noindent\underbar{Démonstration}. Dans ce cas, les équations 3.23 se réduisent à $[P,\frac{\partial^2P}{\partial \bar z \partial z}]=0\Leftrightarrow$
$\frac{\partial^2P}{\partial \bar z \partial z} = P \frac{\partial^2P}{\partial \bar z\partial z} P + (1-P) \frac{\partial^2P}{\partial z\partial z} (1-P)$. D'autre part 
\[
\frac{\partial}{\partial \bar z} 
\Tr
\left(\frac{\partial P}{\partial z}\right)^2
=2\Tr\left(
\frac{\partial P}{\partial z}\frac{\partial^2P}{\partial \bar z\partial z}\right)
\]
\[
=2\Tr \left\{
\left(P\frac{\partial P}{\partial z}P\frac{\partial^2P}{\partial \bar z\partial z}\right)
+\left((1-P)
\frac{\partial P}{\partial z}(1-P)\frac{\partial^2 P}{\partial \bar z\partial z} 
\right)\right\}
\]
où, pour la dernière égalité, on a utilisé l'équation précédente et l'invariance de la trace par permutation circulaire. Comme (par différentiation de $P^2=P$) on a $P(dP)P=(1-P)dP(1-P)=0$, on a $\frac{\partial}{\partial \bar z}\Tr \left(\frac{\partial P}{\partial z}\right)^2=0$~$.\square$

\subsection{Corollaire (solutions du modèle $P_1(\mathbb C)$ sur $S^2$)}

{\sl Soit $P:\mathbb C\rightarrow P_n(\mathbb C)$ une application harmonique telle que\linebreak[4] $\lim_{\vert z\vert\rightarrow \infty}(dP(z))=0$. On a $\left(\frac{\partial P}{\partial z}\right)^3=0$ et, si $n=1$ on a $\left(\frac{\partial P}{\partial z}\right)^2=0$ (i.e. $P$ est $\pm$ holomorphe en vertu de 3.22) En particulier, toute application harmonique de $P_1(\mathbb C)$ dans $P_1(\mathbb C)$ est $\pm$ holomorphe, (i.e. toutes les solutions du ``modèle $P_1(\mathbb C)$ sur $S^2$"$\simeq P_1(\mathbb C)$ sont les ``instantons" et les ``anti-instantons")}.\\

\noindent\underbar{Démonstration}. Le théorème de Liouville implique que toute fonction holomorphe sur $\mathbb C$ tendant vers zéro à l'infini est nulle. On a donc $\left(\frac{\partial P}{\partial z}\right)^3=0$ en appliquant 3.27 et 3.26. Dans le cas $n=1$,  $\frac{\partial P}{\partial z}$ est une matrice 2$\times$2 et $\left(\frac{\partial P}{\partial z}\right)^3=0\Rightarrow \left(\frac{\partial P}{\partial z}\right)^2=0$.
La projection stéréographique de $S^2$ privée d'un point sur $\mathbb C$ est biholomorphe, donc conforme, et, si $P(z)$ correspond sur $\mathbb C$ à une application de $S^2$ dans $P_1(\mathbb C)$, on a $\lim_{\vert z\vert\rightarrow \infty}(dP(z))=0\Rightarrow$ le résultat pour $S^2$ puisque les équations sont invariantes conformes (3.24).~$\square$

\subsection{Exercices et compléments}

1. Dans $\mathbb R^{2\ell}$, on a $U(\ell)=Sp(\ell, \mathbb R)\cup O(2\ell)=\GL(\ell,\mathbb C)\cap O(2\ell)$. Une variété presque hermitienne $V$ de dimension $2\ell$ est une variété de dimension 2$\ell$ munie d'une $U(\ell)$-structure : Dire que $V$ est presque kählerienne revient à dire que la $Sp(\ell, \mathbb R)$-structure contenant sa $U(\ell)$-structure est intégrable; dire que $V$ est hermitienne revient à dire que la $\GL(\ell,\mathbb C)$-structure contenant sa $U(\ell)$-structure est intégrable; dire que $V$ est Kählerienne revient à dire que la $Sp(\ell,\mathbb R)$-structure et la $\GL(\ell,\mathbb C)$-structure contenant sa $U(\ell)$-structure sont toutes les deux intégrables.\\

\noindent 2. Soit $V$ une variété kählerienne, $(\calo,\varphi)$ une carte holomorphe et soient $(z^k=\varphi^k)$. On utilise la base $dz^1,\dots, dz^n,d\bar z^1,\dots,d\bar z^n$ pour les différents objets. Montrer que, pour les symboles de Christoffel, on a $\Gamma^\alpha_{\beta\bar\gamma}=\Gamma^{\bar a}_{\beta\gamma}=0$. En déduire, en utilisant la formule de 3.23, que toute application holomorphe d'une variété kählerienne dans une autre est harmonique.\\

\noindent 3. On reprend les hypothèses et les notations ci-dessus de 2. Montrer que le tenseur de Ricci est donné par : $R_{\alpha\beta}=R_{\bar\alpha\bar\beta}=0$ et $R_{\alpha\bar\beta}=\frac{\partial^2}{\partial z^\alpha\partial\bar z^\beta}\log\det(G)$ où $G=(G_{\tau\bar\lambda})$ est la matrice hermitienne positive telle que $g(X,Y)=G_{\alpha\bar\beta}(X^\alpha Y^{\bar\beta}+Y^\alpha X^{\bar\beta})$ pour le produit scalaire de deux vecteurs $X,Y$ tangents. En déduire qu'une variété kählerienne est Ricci-plate si et seulement si au voisinage de chaque point on a des coordonnées holomorphes $z^1,\dots,z^n$ et un potentiel de Kähler (3.16) $F$ tels que l'on ait 
$$\det\left(\frac{\partial^2 F}{\partial z^\alpha\partial \bar z^\beta}\right)=1$$
Plus généralement, une variété kählerienne est une variété d'Einstein, i.e. Ricci=$\Lambda g$, si et seulement si au voisinage de chaque point on a des coordonnées holomorphes ($z^\alpha$) et un potentiel de Kähler $F$ satisfaisant
\[
\det \left(\frac{\partial^2F(z)}{\partial z^\alpha\partial\bar z^\beta}\right)=e^{\Lambda F(z)}
\]

\noindent 4. Toute sous-variété complexe d'une variété kählerienne est une sous-variété minimale (c'est relié étroitement au fait que l'injection correspondante est harmonique).\\

\noindent 5. $P:\mathbb C\rightarrow G_{N,p}(\mathbb C)$ une application harmonique ($\Rightarrow [P,\frac{\partial^2P}{\partial z\partial\bar z}]=0$). Montrer que l'on a $\frac{\partial}{\partial\bar z}\Tr\left\{(\frac{\partial P}{\partial z})^k\right\}=0$ pour tout entier $k$. En déduire que si, en plus, $\lim_{\vert z\vert\rightarrow \infty}(dP(z))=0$, alors on a $\left(\frac{\partial P}{\partial z}\right)^{2p+1}=0$, (et bien sûr $\left(\frac{\partial P}{\partial z}\right)^N=0$, ce qui n'est intéressant que lorsque $N<2p+1$). Une application holomorphe ou anti-holomorphe $P:\mathbb C \rightarrow G_{N,p}(\mathbb C)$ vérifie $\left(\frac{\partial P}{\partial z}\right)^2=0$, 
mais la réciproque n'est vraie que pour $p=1$, (i.e. le cas des espaces projectifs, lemme 3.22).
\newpage

\begin{center}
{\large Références pour le chapitre 3}
\end{center}
\vspace{1cm}
\begin{itemize}
\item
Kobayashi S., Nomizu K. : Foundations of differential geometry, Vol. II. Interscience Publishers, 1969.\\

\item
Wells R.O. : Differential analysis on complex manifolds. Springer Verlag, 1980.\\

\item
Chern S.S. : Complex manifolds without potential theory. Van Nostand, 1967.\\

\item
Eells J.,  Lemaire L. : A report on harmonic maps. {\sl Bull. London Math. Soc.\/ \bf 10} (1978), 1-68.\\

\item
Newlander A., Nirenberg L. : Complex analytic coordinates in almost complex manifolds. {\sl Ann. of Math. \/ \bf 65} (1957), 391-404.
\end{itemize}

\newpage
\section{Fibrés holomorphes}

\subsection{Fibrés principaux holomorphes}
Soit $V$ une variété complexe et $G$ un groupe de Lie complexe. $P(V,G)$ étant un $G$-fibré principal sur $V$, un système de trivialisations locales $(\calo_\alpha, s_\alpha)_{\alpha\in I}$ de $P(V,G)$ (i.e. $\calo_\alpha$ est un recouvrement ouvert de $V$ et $s_\alpha:\calo_\alpha\rightarrow P(V,G)$ est une section de $P(V,G)\restriction \calo_\alpha$, $\forall \alpha \in I$) sera appelé {\sl système de trivialisations locales compatible avec les structures complexes de $V$ et de $G$} si pour $\alpha, \beta\in I$ avec $\calo_\alpha\cap \calo_\beta\not=\emptyset$ les {\sl fonctions (de transition)} $g_{\alpha\beta}:\calo_\alpha\cap \calo_\beta\rightarrow G$ définies par $s_\alpha g_{\alpha\beta}=s_\beta$ sur $\calo_\alpha\cap \calo_\beta$ sont holomorphes. En adjoignant au système précédent les sections $s:\calo\rightarrow P(V,G)$ où $\calo$ est un ouvert de $V$ et où $\forall\alpha\in I$ avec $\calo\cap\calo_\alpha\not=\emptyset$, $g_\alpha:\calo\cap\calo_\alpha\rightarrow G$, telle que $s=s_\alpha g_\alpha$ sur $\calo\cap\calo_\alpha$, est holomorphe, on obtient un {\sl système complet de trivialisations locales compatibles avec les structures complexes de $V$ et $G$} contenant le système initial: {\sl Un $G$-fibré principal holomorphe sur} $V$ est $G$-fibré principal sur $V$ muni d'un système complet de trivialisations locales compatibles avec les structures complexes de $V$ et $G$; les sections locales (i.e. trivialisations locales) du système sont appelées {\sl sections holomorphes}.

\subsection{Exemples}

a) Soit $V$ une variété complexe; son fibré tangent $T(V)$ est canoniquement un fibré vectoriel complexe et le fibré des repères complexes correspondants est canoniquement un $\GL(\ell,\mathbb C)$-fibré principal holomorphe sur $V$ ($\ell=\dim_{\mathbb C}V$), (les repères naturels correspondant aux cartes holomorphes de $V$ forment un système de trivialisations locales compatibles avec la structure complexe de $V$ et celle de $\GL(\ell,\mathbb C)$.\\

b) Le fibré des repères complexes du fibré tautologique sur la grassmannienne $G_{N,p}(\mathbb C)$ peut être identifié à l'ensemble des matrices $N\times p$ complexes de rang $p$, $H\in M_{Np}(\mathbb C)$ une section sur $\calo\subset G_{N,p}(\mathbb C)$ est (avec les notations de 3), une application $P\mapsto H(P)$ telle que $PH=H$ ($P=H(H^+H)^{-1}H^+$). Les sections $P\mapsto H(P)$ qui sont holomorphes comme fonctions à valeurs dans $M_{Np}(\mathbb C)=\mathbb C^{Np}$ définissent une structure de $GL(p,\mathbb C)$-fibré principal holomorphe sur $G_{Np}(\mathbb C)$ pour ce fibré de repères.

\subsection{Fibrés vectoriels holomorphes}

Soit $E\rightarrow V$ un fibré vectoriel complexe de rang $p$ sur une variété complexe $V$; nous dirons que $E$ est {\sl un fibré vectoriel holomorphe sur} $V$ si son fibré de repères est muni d'une structure de $\GL(p,\mathbb C)$-fibré principal holomorphe sur $V$. On a alors une notion évidente de sections holomorphes de $E$ (les composantes par rapport aux repères holomorphes sont holomorphes).
Si $F\stackrel{\pi}{\rightarrow} V$ est un fibré principal ou vectoriel holomorphe sur $V$, on a une structure unique de variété complexe sur $F$ (que nous appellerons sa {\sl structure complexe canonique}) pour laquelle la projection et les sections holomorphes sont des applications holomorphes (voir 3.21).

\subsection{Exemples}

a) Soit $P(V,G)$ un $G$-fibré principal holomorphe sur $V$ et $\rho:G\rightarrow \End(\mathbb C^p)$ une représentation linéaire complexe de $G$ dans $\mathbb C^p$. Alors le fibré associé $E=\{(r,v)\in P(V,G)\times \mathbb C^p\}/\{(rg^{-1},\rho(g)v)\vert g\in G\}$ est un fibré vectoriel holomorphe de rang $p$ sur $V$.\\

b) Si $E$ et $F$ sont des fibrés vectoriels holomorphes sur $V$, il en est canoniquement de même pour $E\oplus F$, $E\otimes F$, $E^\ast$, $\vee^pE$, ${\wedge}^p E,\dots$\\

c) Soient $V$ et $W$ deux variétés complexes, $f:V\rightarrow W$ une application holomorphe et $F\stackrel{\pi}{\rightarrow} W$ un fibré vectoriel (resp. un $G$-fibré principal) holomorphe sur $W$. Alors l'image inverse $f^\ast(F)=\{(v,\varphi)\in V\times F\vert f(v)=\pi(\varphi)\}$ est un fibré vectoriel (resp. un $G$-fibré principal) holomorphe sur $V$.

\subsection{Homomorphismes et formes horizontales}

On définira les homomorphismes de fibrés holomorphes (de types précé\-dents) comme étant les homomorphismes de fibrés qui sont des applications holomorphes pour leurs structures complexes canoniques.\\
Soit $F\stackrel{\pi}{\rightarrow} V$ un fibré; une $p$-forme en $\xi\in F$, $\varphi\in \wedge^p T^\ast_xi(F)$, est appelée $p$-{\sl forme horizontale en } $\xi\in F$ si elle s'annule dès que l'un des vecteurs de $T_\xi(F)$ auxquels on l'applique est {\sl vertical} i.e. tangent à la fibre $\pi^{-1}(\pi(\xi))$ en $\xi$. Si $\varphi$ est une $p$-forme horizontale en $\xi\in F$ elle 
est l'image inverse $\pi^\ast(\varphi_0)$ par $\pi$ d'une $p$-forme $\varphi_0$ de $V$ en $\pi(\xi)$. On a, par tensorisation, la notion de $p$-forme horizontale à valeurs dans un espace vectoriel et la notion de $p$-forme différentielle (=champ de $p$-formes) horizontale à valeurs dans un espace vectoriel.\\
Dans le cas où $V$ est une variété presque complexe, on peut définir la notion de composante pure de type $(r,s)$ pour une $(r+s)$-forme horizontale puisqu'elle est l'image inverse d'une forme de $V$. De même si $\varphi$ est une $p$-forme horizontale à valeurs dans un espace vectoriel complexe on peut écrire $\varphi=\sum_{r+s=p} P^{r,s}\varphi$.\\
Rappelons que si $P(V,G)$ est un $G$-fibré principal sur $V$, la courbure $\Omega=d\omega+\frac{1}{2}[\omega,\omega]$ d'une connexion $\omega$ sur $P(V,G)$ est une 2-forme différentiable horizontale sur $P(V,G)$ à valeurs dans l'algèbre de Lie $\fracg$ de $G$. Lorsque $V$ et $G$ sont des variétés complexes ($\Rightarrow \fracg$ est un e.v. complexe) on peut écrire $\Omega=P^{2,0}\Omega+P^{1,1}\Omega+P^{0,2}\Omega$. Avec ces conventions, on a le théorème de Koszul et Malgrange suivant.

\subsection{Théorème}

{\sl Soit $P(V,G)$ un $G$-fibré principal sur $V$, où $V$ est une variété complexe et $G$ est un groupe de Lie complexe, et soit $\omega$ une forme de connexion sur $P(V,G)$ dont la forme de courbure $\Omega$ satisfait $P^{0,2}\Omega=0$. Alors $P(V,G)$ a une structure unique de $G$-fibré principal holomorphe sur $V$ telle qu'une section locale $s$ est holomorphe si et seulement si $P^{0,1}s^\ast(\omega)=0$.}\\

\noindent\underbar{Démonstration}. Il suffit de montrer que tout point de $V$ est dans le domaine d'une section locale $h$ satisfaisant $P^{0,1}h^\ast(\omega)=0$ et que si $h$ et $h'=hg$ sont deux telles sections définies sur le même ouvert $\calo$ de $V$ alors $g:\calo\rightarrow G$ est holomorphe.\\
Si $h$ et $hg$ sont deux sections sur $\calo$, on a $(hg)^\ast(\omega)=ad(g^{-1})(h^\ast(\omega))+g^{-1}dg$; si $P^{0,1}h^\ast(\omega)=0$ et $P^{0,1}(hg)^\ast(\omega)=0$ ceci implique $g^{-1}\bar\partial g=0\Rightarrow \bar\partial g=0$, d'où la dernière assertion.\\
D'autre part, si $s$ est une section quelconque sur un ouvert de carte holomorphe de $V$, on a $0=P^{0,2}s^\ast(\Omega)=\bar\partial P^{0,1}s^\ast(\omega)+\frac{1}{2}[P^{0,1}s^\ast(\omega),P^{0,1}s^\ast (\omega)]$, et on peut appliquer le théorème 2.5. On en déduit que tout point du domaine de $s$ est dans ouvert $\calo$ sur lequel on a une fonction à valeur dans $G$, $g:\calo\rightarrow G$, tel que $P^{0,1}s^\ast(\omega)=0$ sur $\calo$. Il en résulte que tout point de $V$ est dans le domaine d'une section locale $h$ telle que $P^{0,1}h^\ast (\omega)=0$.~$\square$

\subsection{Remarques}

a) En appliquant directement 2.3, on voit que l'on a le résultat précédent en supposant simplement que $P(V,G)$ est de classe $C^2$ et $\omega$ de classe $C^1$ (pour pouvoir définir $\Omega$).\\

\noindent b) Une autre manière de formuler le théorème est de dire que sous les hypothèses de 4.6, {\sl on a une structure unique de fibré principal holomorphe sur $P(V,G)$ pour laquelle $\omega$ est pure de type (1,0) (ou, ce  qui revient au même, pour laquelle les sous-espaces horizontaux sont des sous-espaces vectoriels complexes des espaces tangents à $P(V,G)$)}.\\

\noindent c) $P(V,G)$ et $\omega$ étant comme dans 4.6, si $\alpha$ est une 1-forme tensorielle de type $ad$ sur $P(V,G)$ la connexion $\omega+\alpha$ satisfait $P^{0,1}(\omega+\alpha)=0$ pour la même structure holomorphe de $P(V,G)$ si et seulement si $P^{0,1}\alpha=0$; {\bf cette dernière condition ne dépend en fait que de la structure complexe de} $V$ puisque $\alpha$ est horizontale. C'est le seul arbitraire sur $\omega$ pour le fibré holomorphe donné; pour le fixer il suffit d'imposer à $\omega$ une condition de réalité. Il reste alors à savoir si une telle connexion existe sur un fibré holomorphe donné. Le théorème suivant de Singer répond essentiellement à ces deux points.

\subsection{Théorème}  

{\sl Soit $P(V,G)$ un $G$-fibré principal holomorphe sur $V$ et soit $Q(V,U)$ un sous-fibré réel de $P(V,G)$ dont le groupe de structure $U$ est tel que son algèbre de Lie $\fracu$ est une forme réelle de l'algèbre de Lie $\fracg$ de $G$, (i.e. $\fracg=\fracu \oplus_{\mathbb R} i\fracu$ ``comme algèbres réelles"). Alors, il existe une connexion unique sur $Q(V,U)$ dont l'extension  à $P(V,G)$ a une forme de connexion $\omega$ satisfaisant $P^{0,1}\omega=0$}.\\

\noindent \underbar{Démonstration}. Sur $\fracg$ on a une involution antilinéaire de  \[
u_1+iu_2\mapsto (u_1+iu_2)^c=u_1-iu_2
\]
 correspondant à la décomposition réelle $\fracg=\fracu\oplus i\fracu$; $\fracu$ est l'ensemble des points fixes de $\fracg$ pour cette involution. Soit $(\calo_\alpha)$ un recouvrement ouvert de $V$ sur lequel on a des sections $s_\alpha:\calo_\alpha\rightarrow Q(V,U)$ et des sections holomorphes $h_\alpha:\calo_\alpha\rightarrow P(V,G)$; on a alors des fonctions $g_\alpha:\calo_\alpha\rightarrow G$ telles que $s_\alpha=h_\alpha g_\alpha$ sur $\calo_\alpha$. Supposons qu'il existe $\omega$ comme dans l'énoncé de 4.8, on doit avoir $P^{0,1}s^\ast_\alpha(\omega)=P^{0,1}(ad g^{-1}_\alpha)h^\ast_\alpha(\omega)+g^{-1}_\alpha g_\alpha)=g^{-1}_\alpha\bar \partial g_\alpha$ et, par conséquent, $s^\ast_\alpha(\omega)=g^{-1}_\alpha\bar{\partial g}_\alpha+(g_\alpha^{-1}\bar\partial g_\alpha)^c$ puisque $s^\ast_\alpha(\omega)$ est à valeurs dans $\fracu$ et $P^{1,0}s^\ast_\alpha(\omega)=(P^{0,1}s^\ast_\alpha(\omega))^c$. Il en résulte que $\omega$ est unique et donnée dans les sections $s_\alpha$ par la formule précédente. Il reste à vérifier que cette formule définit bien une connexion (i.e. que l'on a $s^\ast_\beta(\omega)=(s_\alpha g_{\alpha\beta})^\ast(\omega)=ad g^{-1}_{\alpha\beta}s^\ast_\alpha(\omega)+ g^{-1}_{\alpha\beta}dg_{\alpha\beta}$ dans $\calo_\alpha\cap \calo_\beta$), ce qui est facile.~$\square$\\
Terminons par une définition courante.

\subsection{Fibré canonique}

Soit $V$ une variété complexe de dimension complexe $\ell$. Le fibré $K=\wedge^{\ell,0}T^\ast(V)$ est un fibré holomorphe de rang 1 sur $V$ appelé {\sl fibré canonique} de $V$. On remarquera que si $V$ est presque complexe de dimension $2\ell$, $\wedge^{\ell,0}T^\ast(V)$ est encore un fibré vectoriel complexe de rang 1.

\subsection{Exercices et compléments}

1. Soit $V$ une variété complexe et soit $E\rightarrow V$ un fibré vectoriel holomorphe hermitien sur $V$ (i.e. les fibres sont des espaces de Hilbert). Montrer en utilisant 4.8 qu'il existe une connexion unique sur $E$ qui préserve le produit scalaire dans les fibres et est pure de type (1,0).\\

\noindent 2. Soit $U\rightarrow P_n(\mathbb C)$ le fibré tautologique sur $P_n(\mathbb C)$ et soit $K\rightarrow P_n(\mathbb C)$ le fibré canonique. Ce sont des fibrés holomorphes de rang 1 et on a : $K\simeq U^{\otimes (n+1)}$. Tout fibré de rang 1 holomorphe sur $P_n(\mathbb C)$ est isomorphe à une puissance tensorielle de $U$ ou de son dual; on pose $U^{\ast\otimes k}=\calo(k)$ et $U^{\otimes k}=\calo(-k)$. La raison de cette notation est que les sections locales de $\calo(k)$ ($k\in \mathbb Z$) sont les fonctions homogènes de degré $k$ dans les ouverts correspondants de $\mathbb C^{n+1}\backslash \{0\}$. Soit $E\rightarrow P_1(\mathbb C)$ un fibré vectoriel holomorphe de rang $n$. Alors $E\simeq \calo(k_1)\oplus\dots\oplus \calo(k_n)$ pour $k_1\leq k_2\dots \leq k_n$ dans $\mathbb Z$.

\newpage
\begin{center}
{\large Références pour le chapitre 4}
\end{center}
\vspace{1cm}
\begin{itemize}
\item
Koszul J.-L., Malgrange B. : Sur certaines structures fibrées complexes. {\sl Arch. Math.\/ \bf 9} (1958), 102-109.\\

\item
Singer I.M. : Geometric interpretation of a special connection. {\sl Pacific. J. Math.\/ \bf 9 (1959)}, 585-590.
\end{itemize}

\newpage
\section{Structures complexes sur les espaces vectoriels réels}

Pour fixer les notations, nous rappelons dans ce paragraphe un certain nombre de définitions et de résultats concernant les structures complexes sur les espaces vectories réels de dimensions finies.

\subsection{Définition}

Soit $E$ un espace vectoriel réel; {\sl une structure complexe sur $E$} est un endomorphisme $J$ satisfaisant $J^2=-\bbbone$ ($\bbbone$ est l'application identique de $E$). Nous désignerons par $\cali(E)$ l'ensemble des structures complexes sur $E$.\\
Si $J$ est une structure complexe sur $E$ on munit $E$ d'une structure d'espace vectoriel complexe en posant $zV=xV+yJV$ pour $z=x+iy\in \mathbb C$ et $V\in E$; nous désignerons l'espace vectoriel complexe ainsi obtenu par $E_J$.\\
Inversement, si $F$ est un espace vectoriel complexe, la multiplication par $i\in \mathbb C$ définit une structure complexe $J$ sur l'espace vectoriel réel sous-jacent $E$ et on a $F=E_J$.\\
Par la suite, nous ne nous intéresserons qu'au cas des espaces de dimensions finies; comme $\dim E=2\dim E_J$, on supposera que $\dim E=2\ell$ ($\ell\in \mathbb N$) afin que $\cali(E)$ soit non vide. $E$ est isomorphe à $\mathbb R^{2\ell}$ et $E_J$ à $\mathbb C^\ell$.

\subsection{Lemme}

{\sl $E, J$ et $E_J$ sont comme précédemment.\\

\noindent a) Pour toute base $(e_1,\dots,e_\ell)$ de $E_J$, ($e_1,\dots,e_\ell,Je_1,\dots,Je_\ell$) est une base de $E$.}\\

\noindent b) {\sl Soit $L\in \End(E)$ un endomorphisme de $E$; alors $L$ définit un endomorphisme $L_J$ de $E_J$ ($L_J\in \End(E_J))$ si et seulement si $LJ=JL$. Dans ce cas on a : $\det(L)=\vert \det(L_J)\vert^2$ où $\det(\bullet)$ désigne le déterminant.}\\

Il résulte de la partie a) de ce lemme que le groupe linéaire $\GL(E)$ opère transitivement sur $\cali(E)$ (par $J\mapsto GJG^{-1}$); il résulte de b) que le stabilisateur d'un point $J$ de $\cali(E)$ peut être identifié au groupe linéaire complexe $\GL(E_J)$ et que, considéré comme plongé dans $\GL(E)$, $\GL(E_J)$ est un sous-groupe de $\GL_+(E)=\{G\in \GL(E)\vert \det(G)>0\}$. Ceci implique que si $J\in \cali(E)$, l'orientation de la base (réelle) de $E$ ($e_1,\dots, e_\ell,Je_1,\dots,Je_\ell$) associée à la base ($e_1,\dots,e_\ell$) de $E_J$ est indépendante de la base ($e_1,\dots,e_\ell$) de $E_J$ choisie mais ne dépend que de $J$; on est conduit à la définition suivante~:

\subsection{Orientation de $J$}

L'orientation $o(J)$ de $E$ définie ci-dessus pour $J\in \cali(E)$ sera appelée {\sl l'orientation de} $J$; si $o$ est une orientation de $E$, nous dirons que $J$ {\sl est compatible avec } $o$ si $o=o(J)$.

\subsection{Remarque et notations}

Il y a d'autres conventions possibles pour définir l'orientation de $J\in \cali(E)$; une autre convention couramment utilisée consiste à choisir l'orientation de ($e_1, Je_1,e_2,Je_2,\dots,e_\ell,Je_\ell$) au lieu de celle de ($e_1,\dots,e_2,Je_1,\dots,Je_\ell$). (Ces deux conventions sont opposées si $\dim(E)=4$ ($\ell=2$) par exemple ! Il faut faire un peu attention).\\

Si $o$ est une orientation de $E$, $\cali(E,o)$ désignera l'ensemble des structures complexes sur $E$ compatibles avec $o$. Si $E$ est orienté nous utiliserons la notation $\cali_+(E)$ (resp. $\cali_-(E)$) pour désigner l'ensemble des structures complexes compatibles avec l'orientation de $E$ (resp. l'orientation opposée à celle de $E$); par exemple, nous utiliserons les notations $\cali_+(\mathbb R^{2\ell})$ et $\cali_-(\mathbb R^{2\ell})$.\linebreak[4] $J\in \cali(E)$ étant donnée, l'application $G\mapsto GJG^{-1}$ de $\GL(E)$ sur $\cali(E)$ permet les identifications 
\[
\cali(E)\simeq \GL(E)/\GL(E_J)\ \text{et}\  \cali(E,o(J))\simeq \GL_+(E)/\GL(E_J).
\]
 Ces identifications permettent de munir $\cali(E)$ de structures topologiques, différentiables, etc. On vérifiera que ces structures sont compatibles avec celles qui sont induites par l'inclusion $\cali(E)\subset \End(E)$, $(\End(E)\simeq \mathbb R^{(2\ell)^2})$.\\

Soit $E^\ast_c=E^\ast\otimes_\mathbb R \mathbb C$ le complexifié du dual de $E$ muni de sa conjugaison complexe canonique comme complexifié d'un espace vectoriel réel. Le dual de $E_J$ pour $J\in \cali(E)$ (comme espace vectoriel complexe) s'identifie au sous-espace $\wedge^{1,0}E^\ast_J$ de $E^\ast_c$ défini par :
\[
\wedge^{1,0}E^\ast_J=\{\omega\in E^\ast_c \vert \omega\circ J=i\omega\}.
\]
Nous désignerons par $\wedge^{0,1}E^\ast_J$ le complexe conjugué dans $E^\ast_c$ de $\wedge^{1,0}E^\ast_J$. On a manifestement $E^\ast_c=\wedge^{1,0}E^\ast_J\oplus \wedge^{0,1}E^\ast_J$, $\forall J\in \cali(E)$, et plus précisément :

\subsection{Proposition}

{\sl L'application $J\mapsto \wedge^{1,0}E^\ast_J$ est une bijection de $\cali(E)$ sur l'ensemble des sous-espaces (complexes) $F$ de $E^\ast_c$ qui sont supplémentaires à leurs complexes conjugués dans $E^\ast_c$ (i.e. $E^\ast_c=F\oplus \bar F$).}\\

L'intérêt de cette proposition est qu'elle permet d'identifier $\cali(E)$ à un ouvert dense de la grassmannienne $G_\ell(E^\ast_c)$ des sous-espaces de dimension = $\ell$ de $E^\ast_c$. $G_\ell(E^\ast_c)\simeq G_{2\ell,\ell}(\mathbb C)$ est une variété algébrique complexe et nous munirons désormais $\cali(E)$ de la struture de variété analytique complexe induite. Si $o$ est une orientation, la frontière commune dans $G_\ell(E^\ast_c)$ des deux ouverts disjoints $\cali(E,o)$ et $\cali(E,-o)$ contient l'ensemble des sous-espaces de $E^\ast_c$ invariant par conjugaison complexe de dimension =$\ell$ que l'on peut identifier à la grassmannienne réelle $G_\ell(E^\ast)$ des sous-espaces
de dimension =$\ell$ dans l'espace réel $E^\ast$, $(G_\ell(E^\ast)\subset G_\ell(E^\ast_c)$ est d'ailleurs une partie très particulière de la frontière de $\cali(E,o)$ dans $G_\ell(E^\ast_c)$).

\subsection{Exemple : le cas de $\mathbb R^2$} 

Dans le cas $E=\mathbb R^2$ on peut facilement visualiser la situation. Soient $e_1=(1,0)$ et $e_2=(0,1)$ les éléments de la base canonique $(e_1,e_2)$ et considérons $\mathbb R^2$ comme orienté par cette base. Soient $I_+$ et $I_-$ les structures complexes satisfaisant à $I_+e_1=e_2$ et $I_-e_2=e_1$ $(I_-=-I_+$ et $I_+$ et $I_-$ sont les seules structures complexes isométriques pour la structure euclidienne pour laquelle $(e_1,e_2)$ est orthonormale). On a $G_\ell(E^\ast_c)\simeq P_1(\mathbb C)\simeq S^2$, $G_\ell(E^\ast)\simeq P_1(\mathbb R)\simeq S^1$; $I_+$ et $I_-$ sont deux points diamétralement opposés de $S^2$ tandis que la grassmannienne réelle $G_\ell(E^\ast)\sim S^1$ est le grand cercle de $S^2$ perpendiculaire à l'axe $I_+I_-$, $\cali_+(\mathbb R^2)$ (resp. $\cali_-(\mathbb R^2)$) est la calotte de $S^2$ contenant $I_+$ (resp. $I_-$) délimitée par $G_\ell(E^\ast)\simeq S^1\subset G_\ell (E^\ast_c)\simeq S^2$.

\subsection{Cas général}

Revenons au cas général. Pour tout $J\in \cali(E)$, l'espace $\wedge^{1,0}E^\ast_J$ est la fibre en $J$ {\sl d'un fibré vectoriel complexe holomorphe de rang (complexe) $\ell$ sur $\cali(E)$} qui est la restriction à $\cali(E)\subset G_\ell(E^\ast_c)$ du fibré tautologique sur la grassmannienne complexe $G_\ell(E^\ast_c)$ ($\simeq G_{2\ell,\ell}(\mathbb C)$); {\sl ce fibré sera noté} $\wedge^{1,0}E^\ast$. On a de même {\sl un fibré vectoriel anti-holomorphe $\wedge^{0,1}E^\ast$ sur $\cali(E)$} dont la fibre en $J$ est $\wedge^{0,1}E^\ast_J$. La somme directe de $\wedge^{1,0}E^\ast$ et $\wedge^{0,1}E^\ast$ est le fibré trivial $E^\ast_c\times \cali(E)$ sur $\cali(E)$ de fibre $E^\ast_c$:
\[
\wedge^{1,0}E^\ast\oplus \wedge^{0,1}E^\ast=E^\ast_c\times \cali(E)
\]
Soit $\wedge E^\ast_c=\oplus_{r\geq 0}    \wedge^r E^\ast_c$ l'algèbre extérieure (complexe) construite sur $E^\ast_c$. Pour $J\in \cali(E)$ on définit les sous-espaces $\wedge^{p,q}E^\ast_J$ de $\wedge E^\ast_c$ par récurrence sur les entiers $p$ et $q$ à partir de $\wedge^{1,0}E^\ast_J$ et $\wedge^{0,1}E^\ast_J$ par :
\[
\wedge^{p+p',q+q'}E^\ast_J = (\wedge^{p,q}E^\ast_J)\wedge (\wedge^{p',q'}E^\ast_J).
\]
$\wedge^{q,p}E^\ast_J$ est le complexe conjugué dans $\wedge E^\ast_c$ de $\wedge^{p,q}E^\ast_J$.\\
 De même que précédemment, on a les {\sl fibrés vectoriels complexes $\wedge^{p,q}E^\ast$ sur} $\cali(E)$ de fibres en $J\in \cali(E)$ $\wedge^{p,q}E^\ast_J$; $\wedge^{p,q}E^\ast\simeq \stackrel{p}{\wedge} (\wedge^{1,0}E^\ast)\otimes \stackrel{q}{\wedge}(\wedge^{0,1}E^\ast)$ et
\[
\oplus_{p+q=r} \wedge^{p,q}E^\ast=\wedge^r E^\ast_c\times \cali(E).
\]
Les fibrés $\wedge^{r,0}E^\ast$ sont des fibrés vectoriels complexes holomorphes sur $\cali(E)$ pour $r=0,1,\dots,\ell$.

\newpage

\section{Fibrés de structures complexes}

Nous avons vu que, si $E$ est un espace vectoriel réel de dimension $2\ell$, le groupe linéaire de $E$ agit (transitivement) sur $\cali(E)\subset \End(E)$. Il en résulte, que si $E(M)$ est un {\sl fibré} vectoriel réel de rang $2\ell$ sur $M$, on a un {\sl fibré associé $\cali(E(M))$ sur $M$} dont la fibre en $x\in M$ est $\cali(E_x)$; le groupe structural est le groupe linéaire $\GL(2\ell, \mathbb R)$ 
. $\cali(E(M))$ est un sous-fibré du fibré $\End(E(M))\simeq E(M)\otimes E^\ast(M)$ par construction.

\subsection{Définition}

Soit $M$ une variété différentiable de dimension $2\ell$ et soit $T(M)$ son fibré tangent; le fibré $\cali(T(M))$ défini comme au-dessus (et muni de sa structure différentiable canonique) sera appelé {\sl le fibré des structures complexes au-dessus de $M$}. Si $M$ est orientable et si $o$ est une orientation de $M$, on a un sous-fibré $\cali(T(M),o)$ qui est le {\sl fibré des structures complexes au-dessus de $M$ compatibles avec l'orientation $o$}. Lorsque $M$ est orientée et que cela n'entraine pas de confusion, nous dénoterons par $\cali_+(T(M))$ (resp. $\cali_-(T(M))$) le fibré des structures complexes au-dessus de $M$ compatibles avec l'orientation de $M$ (resp. l'orientation opposée); $\cali_+(T(M))$ et $\cali_-(T(M))$ sont les deux composantes connexes isomorphes de $\cali(T(M))$.

\subsection{Exemple : le cas de la variété $\mathbb R^{2\ell}$}

Examinons le cas où $M=\mathbb R^{2\ell}$ muni de son orientation canonique. La fibre type de $\cali(T(M))$ est $\cali(\mathbb R^{2\ell})$. Chaque $J\in \cali(\mathbb R^{2\ell})$ est une structure complexe sur $\mathbb R^{2\ell}$ et, par conséquent, définit canoniquement une section $s_J$ du fibré $\cali(T(\mathbb R^{2\ell}))$ sur $\mathbb R^{2\ell}$. Ces sections sont les sections horizontales de $\cali(T(\mathbb R^{2\ell}))$ pour la connexion linéaire plate canonique de $\mathbb R^{2\ell}$. Chaque $s_J$ est isomorphe à $\mathbb R^{2\ell}$ et peut être considéré comme un espace vectoriel complexe de dimension $\ell$ en munissant $\mathbb R^{2\ell}$ de $J\in \cali(\mathbb R^{2\ell})$. On obtient ainsi une fibration de l'espace $\cali(T(\mathbb R^{2\ell}))$ sur $\cali(\mathbb R^{2\ell})$ par des espaces vectoriels complexes de dimension $\ell$. On a le diagramme :

\[
\begin{diagram}
\node{\cali(T(\mathbb R^{2\ell}))}\arrow{e,t}{\mathbb C^\ell}
\arrow{s,r}{\cali(\mathbb R^{2\ell})}\node{\cali(\mathbb R^{2\ell})}\\
\node{\mathbb R^{2\ell}}
\end{diagram}
\]

En fait, cette fibration de $\cali(T(\mathbb R^{2\ell}))$ sur $\cali(\mathbb R^{2\ell})$ est un fibré vectoriel complexe holomorphe sur $\cali(\mathbb R^{2\ell})$; c'est la restriction à $\cali(\mathbb R^{2\ell})$ du dual du fibré tautologique sur la grassmannienne $G_\ell(\mathbb C^{2\ell})=G_{2\ell,\ell}(\mathbb C)$. Pour chaque $x\in \mathbb R^{2\ell}$, $\cali(T_x(\mathbb R^{2\ell}))$ (i.e. la fibre de $\cali(T(\mathbb R^{2\ell}))$ est donc de manière naturelle une variété complexe holomorphe.\\
Nous allons, plus généralement, construire une structure presque complexe sur $\cali(T(M))$ pour toute variété $M$ de dimension $2\ell$ munie d'une connexion linéaire. Cette structure sera compatible avec la structure complexe naturelle des fibres $\cali(T_x(M))$ définie en 5.5; la structure complexe de $\cali(E)$ est invariante par $\GL(E)\Rightarrow$ structure complexe sur les fibres de $\cali(T(M))$.

\subsection{Structure presque complexes et connexions}
Soit $M$ une variété de dimension $2\ell$ munie d'une connexion linéaire, soit $x$ un point de $M$ et $J_x\in \cali(T(M))$ une structure complexe sur $T_x(M)$, $(J_x\in \cali(T_x(M)))$. $\cali(T_x(M))$ étant une variété complexe, on a une structure complexe naturelle sur l'espace tangent en $J_x$ à la fibre $\cali(T_x(M))$ (i.e. les vecteurs verticaux). $J_x$ étant une structure complexe sur $T_x(M)$, on peut, par relèvement horizontal, munir le sous-espace horizontal en $J_x$ de la structure complexe correspondante. Par somme directe, on a ainsi défini, pour chaque $J_x\in \cali(T(M))$ une structure complexe dans l'espace tangent en $J_x$ à $\cali(T(M))$ qui dépend différentiablement de $J_x$. On a donc une structure presque complexe $j$ sur $\cali(T(M))$ que nous appellerons {\sl la structure presque complexe naturelle sur $\cali(T(M))$ associée à la connexion linéaire de $M$.}\\
Dans le cas où $M=\mathbb R^{2\ell}$ muni de sa connexion plate canonique, cette structure presque complexe est intégrable et l'on retrouve ainsi la structure de variété complexe de $\cali(T(\mathbb R^{2\ell})$ décrite dans (6.2). \\
Rappelons que les connexions linéaires sur $M$ forment un espace affine et que la différence de deux connexions linéaires est canoniquement un champ tensoriel $K$ de type (1-2) sur $M$, (i.e. une section de $T(M)\otimes T^\ast(M)\otimes T^\ast(M)$), défini en termes des différentielles covariantes $\nabla$ et $\nabla'$ correspondantes par :
\[
K(X,Y)=\nabla'_XY-\nabla_X Y;
\]
pour toute paire de champs vectoriels ($X,Y$) sur $M$. Le lemme de comparaison suivant est la clef pour analyser la manière dont la structure presque complexe sur $\cali(T(M))$ dépend de la connexion linéaire de $M$.

\subsection{Lemme}
{\sl Soit $x$ un point de $M$ et soit $J\in \cali(T(M))$ une structure complexe sur $T_x(M)$. Les structures presque complexes naturelles sur $\cali(T(M))$ associées à deux connexions linéaires $\nabla$ et $\nabla'$ sur $M$ coincident en $J$ si et seulement si leur différence $K=\nabla'-\nabla$ satisfait l'identité suivante :
\[
K(JX,JY)-JK(JX,Y)-JK(X,JY)-K(X,Y)=0
\]
pour tout couple $(X,Y)$ de vecteurs tangents à $M$ en $X$.}

\noindent \underbar{Démonstration}. L'espace $V_J=T_J(\cali(T_x(M)))$ des vecteurs verticaux en $J$ peut être identifié à l'espace des endomorphismes $A$ de $T_x(M)$ satisfaisant $AJ+JA=0$; l'application $A\mapsto JA$ est une structure complexe naturelle sur $V_J$ qui coincide avec la structure complexe naturelle de $V_J$ comme espace tangent en $J$ à la variété complexe $\cali(T_x(M))$.\\
L'endomorphisme de $T_x(M)$ défini par $Y\mapsto K(X,JY)-JK(X,Y)$ est dans $V_J$ ; c'est, par définition, la projection verticale $\calV_J(X)$ pour $\nabla$ du relèvement horizontal pour $\nabla'$ de $X$ en $J$, ($X\in T_x(M))$.\\
Pour que les structures presque complexes naturelles sur $\cali(T(M))$ associée à $\nabla$ et $\nabla'$ coincident en $J$, il faut et il suffit qu'elles coincident sur l'espace des vecteurs horizontaux pour $\nabla'$ en $J$; ceci équivaut à $J\calV_J(X)=\calV_J(JX)$, $\forall X\in T_x(M)$ et par conséquent à $JK(X,J,Y)+K(X,Y)=K(JX,JY)-JK(JX,Y)\>\>$ $\forall X,Y\in T_x(M)$.~$\square$\\

Rappelons que (voir Chapitre 3) si $j$ est une structure pesque complexe sur une variété $V$, on définit un champ tensoriel de type (1-2) sur $V$ (ou plutôt une 2-forme à valeurs vectorielles), appelée torsion de $j$, par 
\[
N(X,Y)=2([jX,jY]-j[X,jY]-j[jX,Y]-[X,Y]),
\]
pourtout couple $X,Y$ de champs vectoriels sur $V$. L'annulation de $N$ est la condition nécessaire et suffisante pour que $j$ soit induite par une structure, alors unique, de variété complexe sur $V$.\\
Pour être plus complet, nous donnons sans démonstration dans le lemme suivant l'expression de la torsion $N$ de la structure presque complexe naturelle sur $\cali(T(M))$ associée à une connexion linéaire sur $M$ en termes de la courbure et de la torsion de cette connexion.

\subsection{Lemme}
{\sl
Soit $M$ une variété de dimension $2\ell$ munie d'une connexion linéaire de courbure $R$ et de torsion $T$ et soit $N$ la torsion de la structure presque complexe naturelle sur $\cali(T(M))$ associée à la connexion linéaire de $M$. Soit $x$ un point de $M$ et soient $\tilde X$ et $\tilde Y$ deux vecteurs tangents en $J\in \cali(T_x(M))$ à $\cali(T(M))$ de projections $X$ et $Y$ sir $T_x(M)$. Alors, la projection sur $T_x(M)$ de $N(\tilde X,\tilde Y)$ est donnée par}
\[
2(T(X,Y)+KT(X,JY)+JT(JX,Y)-T(JX,JY))
\]
{\sl tandis que sa composante verticale est}
\[
\begin{array}{lll}
4(R(X,Y)J&-&JR(X,Y)+JR(JX,JY)-R(JX,JY)J\\
&+&JR(X,JY)J+R(X,JY)+JR(JX,Y)J+R(JX,Y))
\end{array}
\]

\subsection{Exemples}

Soient $\nabla$ et $\nabla'$ deux connexions linéaires sur $M$ telles que l'on ait pour tout champ de vecteurs $X$ et $Y$ sur $M$ :\\
$\nabla'_XY-\nabla_XY=K(X,Y)=\alpha(X)Y+\beta(Y) X$, où $\alpha$ et $\beta$ sont des 1-formes sur $M$.\\
Alors, le lemme 6.4 implique les structures presque complexes naturelles sur $\cali(T(M))$ associées à $\nabla$ et à $\nabla'$ sont identiques. C'est le cas si $\nabla$ et $\nabla'$ sont ``projectivement reliées", i.e. $\alpha=\beta$.\\
Ceci implique en particulier que $\cali(T(M))$ est une variété complexe pour la structure presque complexe associée à une connexion linéaire projectivement plate sur $M$. On peut d'ailleurs le vérifier directement en appliquant le lemme 6.5.

\subsection{Intérêt des lemmes 6.4 et 6.5}

Un intérêt des lemmes 6.4 et 6.5 est que ce sont des résultats ponctuels dans $\cali(T(M))$ ce qui permet de les utiliser dans le cas où on a une réduction du groupe structural à un sous-groupe $G$ de $\GL(2\ell,\mathbb R)$ et que l'on considère un sous-fibré $\cali_G(T(M))$ de $\cali(T(M))$ associé à cette réduction tel que $\ \cali_G(T_x(M))$ soit une sous-variété complexe de $\cali(T_x(M))$.\\
Explicitement soit $G$ un sous-groupe du groupe linéaire $\GL(2\ell, \mathbb R)$, soit $\cali_G(\mathbb R^{2\ell})$ une sous-variété complexe de $\cali(\mathbb R^{2\ell})$ invariante par $G$ et soit $P(M,G)$ une $G$-structure sur $M$ (i.e. une réduction à $G$ du fibré des repères de $M$). On a alors un fibré associé à $P(M,G)$, $\cali_G(T(M))$, dont la fibre type est $\cali_G(\mathbb R^{2\ell})$, qui est un sous-fibré de $\cali(T(M))$. $\cali_G(T(M))$ est une sous-variété presque complexe de $\cali(T(M))$ pour la structure presque complexe naturelle associée à une connexion sur $P(M,G)$. Deux telles connexions peuvent donner la même structure presque complexe sur $\cali_G(T(M))$ sans que cela soit le cas pour $\cali(T(M))$ tout entier, on a alors l'identité du lemme 6.5 pour les $J\in \cali_G(T(M))$ et pas pour tous les $J\in \cali(T(M))$.\\
Bien qu'il y ait d'autres cas intéressants, dans la suite de cet exposé nous ne nous intéresserons qu'au cas
où $G=SO(2\ell), O(2\ell)$ ou le groupe conforme et où $\cali_G(\mathbb R^{2\ell})$ est l'ensemble des structures complexes isométriques sur $\mathbb R^{2\ell}$. Pour cela nous rappelons, dans le paragraphe suivant, un certain nombre de propriétés relatives aux structures complexes isométriques sur les espaces euclidiens.\\
Avant cela, terminons ce paragraphe par quelques définitions qui seront utiles par la suite.

\subsection{Définitions}

Soit $M$ une variété de dimension $2\ell$; les éléments du fibré cotangent à $\cali(T(M))$ s'annulant sur les vecteurs tangents aux fibres de $\cali(T(M))\rightarrow M$ (i.e. sur les vecteurs verticaux) seront appelés {\sl formes horizontales}. L'ensemble des formes horizontales est canoniquement un sous-fibré $\wedge^1_H \cali(T(M))$ du fibré cotangent à $\cali(T(M))$; $\wedge^1_H
\cali(T(M))$ n'est autre que l'image inverse du fibré cotangent à $M$ par la projection $\cali(T(M))\rightarrow M$. On a de même une notion de formes horizontales complexes ($\in \wedge^1_H\cali(T(M))\otimes \mathbb C)$. Une forme horizontale complexe en $J_x\in \cali(T_x(M))$ sera dite {\sl pure de type} (1,0), (resp. (0,1)), si elle est l'image inverse par $\cali(T(M))\rightarrow M$ d'un élément de $\wedge^{1,0}T_x(M)^\ast_{J_x}$ (resp. $\wedge^{0,1}T_x(M)^\ast_{J_x}$). L'ensemble $\wedge^{1,0}_H\cali(T(M))$  des formes horizontales pures de type (1,0) et l'ensemble $\wedge^{0,1}_H\cali(T(M))$ des formes horizontales pures de type (0,1) sont deux sous-fibrés supplémentaires de $\wedge^1_H\cali(T(M))\otimes \mathbb C$. On définira $\wedge^{p,q}_H\cali(T(M))$ par récurrence et par produits extérieurs (comme dans 5.7). On a donc défini une notion de {\sl forme horizontale pure de type $(p,q)$ sur $\cali(T(M))$} qui ne se réfère pas à une structure presque complexe particulière de $\cali(T(M))$; cette définition est justifiée par le lemme suivant dont la démonstration est immédiate.

\subsection{Lemme}

{\sl Si $j$ est la structure presque complexe naturelle sur $\cali(T(M))$ associée à une connexion linéaire sur $M$, alors les formes horizontales pures de type $(p,q)$ au sens de la définition précédente sont les formes horizontales (de degré $p+q$) qui sont pures de type $(p,q)$
pour la structure presque complexe $j$ sur $\cali(T(M))$.}

\subsection{Remarque}

$\wedge^{\ell,0}_H\cali(T(M))$ est un fibré vectoriel complexe de rang 1 sur $\cali(T(M))$ et on remarquera que sa restriction à une section $J$ de $\cali(T(M))\rightarrow M$ (si une telle section existe) n'est autre que le fibré canonique $\wedge^{\ell,0} T^\ast(M,J)$ de la variété presque complexe ($M,J$). Le fibré $\wedge^{\ell,0}_H\cali(T(M))$ jouera un certain rôle par la suite; l'existence d'une ``racine carrée" de ce fibré, i.e. d'un fibré vectoriel complexe $\eta$ de rang 1 tel que $\eta\otimes \eta\simeq \wedge^{\ell,0}_H\cali(T(M))$, caractérise l'annulation de la 2ème classe de Stiefel-Whitney de $M$, les ``racines carrées" étant alors classifiées par $H^1(M,\mathbb Z/2)$.

\newpage

\section{Structures complexes isométriques sur les espaces euclidiens}

Dans tout ce paragraphe, $E$ désigne un espace euclidien (réel) de dimension $2\ell$ dont le produit scalaire sera noté $(\bullet, \bullet)$.

\subsection{Définition de $\calh(E)$}

Soit $\calh(E)$ l'ensemble des structures complexes isométriques sur $E$. On vérifie facilement qu'une structure complexe sur $E$ est une transformation conforme si et seulement si elle est isométrique; c'est pourquoi nous dirons qu'une structure complexe de $\calh(E)$ est {\sl compatible avec la structure conforme de $E$}.\\
Si $J$ est dans $\calh(E)$ nous munirons $E_J$ d'une structure d'espace hilbertien en définissant son produit scalaire sesquilinéaire $\langle\bullet\vert\bullet\rangle$ par :
\[
\langle X\vert Y\rangle=(X,Y)-i(X,JY),\ \ \forall X,Y\in E\ \ (=E_J\ \text{comme ensemble})
\]
On a $(X,X)=\vert\vert X\vert\vert^2=\langle X\vert X\rangle$, car $\calh(E)$ est  aussi l'ensemble des rotations anti-symétriques de $E$.\\
Comme sous-groupe de $\GL(E)$ le groupe $O(E)$ des rotations de $E$ agit sur $\cali(E)$; $\calh(E)$ est une orbite pour cette action et le stabilisateur de $J\in \calh(E)$ s'identifie canoniquement au groupe unitaire $U(E_J)$ de $E_J$. On peut identifier $\calh(E)$ à $O(E)/U(E_J)$. De même, si $o$ est une orientation de $E$, $\calh(E,o)=\calh(E)\cap \cali(E,o)$ s'identifie à $SO(E)/U(E_J)$.\\
Si $E$ est orienté, $\calh_+(E)$ (resp. $\calh_-(E)$) désignera l'ensemble des structures complexes isométriques compatibles avec l'orientation de $E$ (resp. l'orientation opposée à celle de $E$).\\
Reprenons les notations de 5.5. Le dual $E^\ast$ de $E$ est un espace euclidien comme dual d'espace euclidien; son produit scalaire peut s'étendre par bilinéarité complexe au complexifié $E^\ast_c$ et nous utiliserons encore la notation $(\bullet,\bullet)$ pour désigner la forme bilinéaire symétrique correspondante sur $E^\ast_c$. Il est facile de voir que si $J$ appartient à $\calh(E)$, $\wedge^{1,0}E^\ast_J$ est maximal isotrope dans $E^\ast_c$ et que, plus précisément, on a :

\subsection{Proposition}

{\sl L'application $J\mapsto \wedge^{1,0}E^\ast$ est une bijection de $\calh(E)$ sur l'ensemble des sous-espaces isotropes de dimension $\ell$ de $E^\ast_c$ pour $(\bullet,\bullet)$.}\\
Cette condition (isotropie) étant algébrique, il en résulte que $\calh(E)$ est canoniquement une variété algébrique complexe qui possède deux composantes connexes isomorphes $\calh(E,\pm o)$, ($o$ étant une orientation de $E$.)

\subsection{Exemples}

\noindent $\calh_+(\mathbb R^2)\simeq \{I_+\}$ (voir 5.6)\\
$\calh_+(\mathbb R^4)\simeq P_1(\mathbb C) (\simeq SO(4)/U(2)$; localement $SO(4)\simeq SU(2)\times SU(2))$\\
$\calh_+(\mathbb R^6)\simeq P_3(\mathbb C)(\simeq SO(6)/U(3)$; localement $SO(6)\simeq SU(4))$\\
Pour $\ell\geq 4$, $\calh_+(\mathbb R^{2\ell})$ est plus compliqué; on comprendra pourquoi et on retrouvera les isomorphismes précédents lorsque nous ferons le lien avec la théorie des spineurs.

\subsection{Structure kählerienne de $\calh(E)$}

Le produit scalaire de $E^\ast$ peut se prolonger aussi, par sesquilinéarité, en un produit scalaire hermitien sur $E^\ast$, que nous noterons $\langle\bullet\vert\bullet\rangle$. On a :
\[
\langle \omega_1\vert \omega_2\rangle=(\bar \omega_1,\omega_2),\ \ \ \forall \omega_1,\omega_2\in E^\ast_c,
\]
et $E^\ast_c$ muni de ce produit scalaire est un espace hilbertien.\\
Pour $J$ dans $\calh(E)$, $E_J$ est un espace hilbertien; son dual $\wedge^{1,0}E^\ast_J$ est donc aussi un espace hilbertien et on vérifie facilement que son produit scalaire coincide avec la restriction de celui de $E^\ast_c$, (i.e. le dual de l'espace hilbertien $E_J$ est un sous-espace hilbertien de $E^\ast_c$).\\
La condition d'isotropie de $\wedge^{0,1}E^\ast_J$ pour $(\bullet,\bullet)$ $(J\in \calh(E))$ est équivalente à l'orthogonalité de $\wedge^{1,0}E^\ast_J$ avec son complexe conjugué $\wedge^{0,1}E^\ast_J$ dans l'espace hilbertien $E^\ast_c$. On a donc au sens des sommes directes hilbertiennes $E^\ast_c=\wedge^{1,0}E^\ast_J\oplus\wedge^{0,1}E^\ast_J$, et plus généralement $\wedge^rE^\ast_c=\oplus_{p+q=r}\wedge^{p,q}E^\ast_J$.\\
Par la suite, nous nous intéresserons presque exclusivement à $\calh(E)$. C'est pourquoi nous utiliserons la notation $\wedge^{p,q}E^\ast$ pour désigner la restriction à $\calh(E)\subset \cali(E)$ du fibré correspondant (voir 5.7).\\
Soit $P^{p,q}_J$ le projecteur hermitien sur $\wedge^{p,q}E^\ast_J$ dans $\wedge E^\ast_c=\oplus_r\wedge^rE^\ast_c$, ($J\in \calh(E))$. La métrique Kählerienne de $\calh(E)$, (induite par celle de la grassmannienne $G_\ell(E^\ast_c)$), est donnée par :
\[
ds^2=1/4 \Tr(dP^{1,0})^2; (J\mapsto 1/4 \Tr (dP^{1,0}_J)^2=(ds^2)_J,\ J\in \calh(E))
\]
le facteur 1/4 est choisi par commodité ici.\\
Le fibré $\wedge^{1,0}T^\ast(\calh(E))$ des formes $\mathbb C$-linéaires sur les espaces tangents à $\calh(E)$ est un fibré vectoriel hermitien holomorphe sur $\calh(E)$ de rang $\frac{\ell(\ell-1)}{2}$, il en est de même de la 2ème puissance extérieure de $\wedge^{1,0}E^\ast$ (i.e. $\wedge^{2,0}E^\ast$). En fait, ils sont isomorphes.

\subsection{Théorème}

{\sl On a un isomorphisme canonique $\chi:\wedge^{2,0}E^\ast\rightarrow \wedge^{1,0}T^\ast(\calh(E))$ de fibrés vectoriels hermitiens holomorphes sur $\calh(E)$}.\\

\noindent\underbar{Démonstration}
Soient $\alpha$ et $\beta$ deux sections locales de $\wedge^{1,0}E^\ast$ définies sur un même ouvert de $\calh(E)$; $\alpha$ et $\beta$ sont en particulier deux fonctions à valeurs dans $E^\ast_c$ et on peut considérer les différentielles $d\alpha$ et $d\beta$ de ces fonctions. Comme $(\alpha,\beta)=0$ (isotropies des $\wedge^{1,0}E^\ast_J$), la forme différentielle à valeurs complexes $(\alpha,d\beta)$ est antisymétrique en $\alpha$ et $\beta$ et est bilinéaire relativement à la multiplication par des fonctions $C^\infty$ sur $\calh(E)$ à valeurs complexes.En d'autres termes, on a un homomorphisme $\chi$ de fibrés vectoriels complexes de $\wedge^{2,0}E^\ast$ dans $T^\ast_c(\calh(E))$ tel que $\chi(\alpha\wedge \beta)=(\alpha,d\beta)$ pour toute paire $\alpha,\beta$ de sections locales de $\wedge^{1,0}E^\ast$.\\
L'holomorphe d'une section locale $\omega$ de $\wedge^{1,0}E^\ast$ est caractérisée par $\bar\partial\omega=0$ (i.e. $d\omega=\partial\omega$). En utilisant une base locale de sections holomorphes de $\wedge^{1,0}E^\ast$ on voit que l'on a $\chi(\alpha\wedge\beta)=(\alpha,\partial \beta)$ ce qui implique que $\chi(\wedge^{2,0}E^\ast)$ est contenu dans $\wedge^{1,0}T^\ast(\calh(E))$ ainsi que l'holomorphie de $\chi:\wedge^{2,0}E^\ast\rightarrow \wedge^{1,0}T^\ast(\calh(E))$.\\
Il nous reste à montrer que $\chi$ induit des isométries sur les fibres. Soit $\omega^1,\dots,\omega^\ell$ une base orthonormale de sections locales de $\wedge^{1,0}E^\ast$. On a :
\[
\begin{array}{ccc}
&&1/4\Tr(dP^{1,0})^2=1/2\Tr P^{1,0}(dP^{1,0})^2=1/2\sum_n\langle dP^{1,0}\omega^n\vert dP^{1,0}\omega^n\rangle\\
\\
&& = 1/2\sum_n\langle d\omega^n\vert (1-P^{1,0})d\omega^n\rangle=1/2\sum_{n,m}\vert (\omega_m,d\omega_n)\vert^2
\end{array}
\]
Ceci implique que la transposée de $\chi$ et par conséquent $\chi$ est une isométrie sur les fibres.~$\square$\\

\subsection{Remarque} 

La 2-forme 
$$\phi=\sum_{m<n} 1/i (\bar\omega^m,d\bar \omega^n)\wedge (\omega^m,d\omega^n)=1/2i \Tr (P^{1,0}dP^{1,0}\wedge dP^{1,0})$$
est une forme fermée de type (1,1) sur $\calh(E)$; c'est la 2-forme fondamentale de la variété kählerienne $\calh(E)$, (voir Chapitre 3).\\
Supposons maintenant que $E$ est orienté. On définit une bijection linéaire $\ast$ de $\wedge E^\ast=\oplus_r\wedge^r E^\ast$ sur lui-même appelée {\sl dualité de Hodge} en posant $\ast(\omega^1\wedge \dots \omega^p)=\omega^{p+1}\wedge\dots\wedge\omega^{2\ell}$ pour toute base orthonormale d'orientation positive $(\omega^1,\dots,\omega^{2\ell})$ de $E^\ast$ et $1\leq p\leq 2\ell$, (on vérifie que cela est cohérent). On étend cette opération, par linéarité complexe, en un isomorphisme $\ast:\wedge E^\ast_c\rightarrow \wedge E^\ast_c$ d'espaces vectoriels complexes. On notera que $\ast$ est isométrique, que $\ast(\wedge^rE^\ast_c)=\wedge^{2\ell-r}E^\ast_c$ et que l'on a :
\[
\alpha\wedge\ast\beta=(\alpha,\beta)\  \vol,\>\> \forall \alpha,\beta \in \wedge^rE^\ast_c
\]
où $\vol\in \wedge^{2\ell}E^\ast$ est l'élément de volume de $E^\ast$, ($\vol=\omega^1\wedge\dots \wedge\omega^{2\ell}, (\omega^1,\dots,\omega^{2\ell})$ étant comme au-dessus). En particulier, $\ast$ est un isomorphisme de $\wedge^\ell E^\ast_c$ sur lui-même et on a le lemme algébrique suivant.

\subsection{Lemme} 

{\sl Soit $\Omega$ un élément de $\wedge^\ell E^\ast_c$. Alors on a $\Omega+i^\ell\ast\Omega=0$, (resp. $\Omega-i^\ell\ast \Omega=0$), si et seulement si $P^{0,\ell}_J\Omega=0$ $\forall J\in \calh_+(E)$, (resp. $\forall J\in \calh_-(E)$)}.\\

\noindent\underbar{Démonstration}. Soit $J\in \calh_+(E)$ et soit ($\alpha^1,\dots,\alpha^\ell$) une base orthonormale de l'espace hilbertien $\wedge^{1,0}E^\ast_J$; alors 
\[
\left(\frac{\alpha^1+\bar\alpha^1}{\sqrt{2}},\dots,\frac{\alpha^\ell+\bar\alpha^\ell}{\sqrt{2}},i\frac{\alpha^1+\bar\alpha^1}{\sqrt{2}},\dots,i\frac{\alpha^\ell+\bar\alpha^\ell}{\sqrt{2}}\right)
\]
est une base orthonormale d'orientation positive de $E^\ast$. Il en résulte que l'on a : $\bar\alpha^1\wedge\dots\wedge\bar\alpha^\ell=i^\ell\ast (\bar \alpha^1\wedge\dots\wedge\bar\alpha^\ell$).\\
Si $\Omega=-i^\ell\ast\Omega$, on a :
\[
\bar\Omega\wedge\ast(\bar\alpha^1\wedge\dots\wedge\bar\alpha^\ell)=\langle\Omega\vert\bar\alpha^1\wedge \dots \wedge\bar\alpha^\ell\rangle\vol,
\]
mais 
\[
\bar\Omega\wedge\bar\alpha^1\wedge\dots\wedge\bar\alpha^\ell=i^\ell\bar\Omega\wedge\ast (\bar\alpha_1\wedge\dots \wedge\bar\alpha^\ell)=\bar\alpha_1\wedge\dots\wedge \bar\alpha^\ell \wedge i^\ell\ast \bar\Omega
\]
\[
=-\bar\Omega\wedge\bar\alpha^1\wedge\dots \wedge \bar\alpha^\ell=0.
\]
Par conséquent, $\Omega=-i^\ell\ast\Omega\Rightarrow \langle \Omega\vert \bar \alpha^1\wedge\dots \wedge \bar\alpha^\ell\rangle=0$, ce qui signifie $P^{0,\ell}_J\Omega=0$. Comme $J\in \calh_+(E)$ est arbitraire on a montré que $\Omega+i^\ell\ast\Omega=0\Rightarrow P^{0,\ell}_J\Omega=0$, $\forall J\in \calh_+(E)$. La réciproque résulte du fait que $\calh_+(E)$ est une orbite dans l'action de $SO(E)$ sur les structures complexes et que les $\Omega$ satisfaisant $\Omega+i^\ell\ast\Omega=0$ forment un sous-espace irréductible pour l'action de $SO(E)$ sur $\wedge^\ell E^\ast_c$. On raisonne de même pour le cas $\Omega-i^\ell\ast\Omega=0$.~$\square$

\subsection{Remarques}

\noindent a) Par produit tensoriel, on a le même énoncé pour les $\ell$-formes à valeurs dans un espace vectoriel complexe $F$, (i.e. les $\Omega\in \wedge^\ell E^\ast\otimes F$).\\

\noindent b) {\sl Définissons $\wedge^\ell_+E^\ast_c$ et $\wedge^\ell_- E^\ast$ par} :
\[
\wedge^\ell_\pm E^\ast_c=\{\Omega\in \wedge^\ell E^\ast_c\vert \Omega\pm i^\ell \ast \Omega=0\}
\]
On a $\wedge^\ell E^\ast_c=\wedge^\ell_+E^\ast_c \oplus \wedge^\ell_- E^\ast_c$ au sens des sommes directes hilbertiennes car $\Omega_\pm \in \wedge^\ell_\pm E^\ast_c\Rightarrow \langle \Omega_+\vert \Omega_-\rangle=0$.\\

\noindent c) {\sl L'énoncé 7.7 est invariant conforme} : les structures complexes isométriques, i.e. $\calh(E)$, et la restriction de la dualité de Hodge, $\ast$, à $\wedge^\ell E^\ast_c$ ne changent pas si on remplace le produit scalaire $(\bullet,\bullet)$ de $E$ par $k(\bullet,\bullet)$ où $k$ est un nombre réel strictement positif.

\subsection{Rétraction par déformation de $\cali(E)$ sur $\calh(E)$}

La décomposition polaire de $J\in \cali(E)$, $J=J_0 e^S$ ($e^{2S}=J^tJ$), donne une bijection de $\cali(E)$ sur l'ensemble des paires $(J_0,S)$ satisfaisant $J_0\in \calh(E)$, $J_0S+SJ_0=0$ et $S^t=S$ qui peut être identifié au fibré normal à $\calh(E)$ dans $\cali(E)$; dans cette bijection, $\calh(E)$ s'identifie à la section nulle de ce fibré normal. L'homotétie dans les fibres donne {\sl une rétraction par déformation de $\cali(E)$ sur $\calh(E)\subset \cali(E)$ $(J_\tau=J_0 e^{\tau S})$.}

\newpage

\section{Fibrés de structures complexes et structures conformes}

Soit $M$ une variété riemannienne; on a un {\sl fibré $\calh(T(M))$ des structures complexes isométriques au-dessus de $M$} dont la fibre en $x\in M$ est l'espace $\calh(T_x(M))$ des structures complexes isométriques sur l'espace euclidien $T_x(M)$, (voir Chapitre 7). $\calh(T(M))$ ne dépend en fait que de la structure conforme de $M$ induite par la métrique.

\subsection{Définition}

Soit $M$ une variété de dimension $2\ell$ munie d'une structure conforme. Le fibré $\calh(T(M))$ des structures complexes au-dessus de $M$ qui sont isométriques pour une métrique compatible avec (i.e. induisant) la structure conforme de $M$ sera appelé {\sl fibré des structures complexes au-dessus de $M$ compatibles avec la structure conforme de $M$}.\\
$\calh(T(M))$ est un sous-fibré de $\cali(T(M))$, et sa fibre $\calh(T_x(M))$ en $x\in M$  est une sous-variété complexe de $\cali(T_x(M))$. Si $g$ est une métrique compatible avec la structure conforme de $M$ alors, $\calh(T(M))$ est une sous-variété presque complexe de $\cali(T(M))$ pour la structure presque complexe naturelle associée à la connexion riemannienne correspondant à $g$. Le théorème suivant montre que la structure presque complexe ainsi définie sur $\calh(T(M))$ ne dépend en fait que de la structure conforme de $M$.

\subsection{Théorème}

{\sl Soit $M$ une variété de dimension $2\ell$ munie d'une structure conforme et soient $g$ et $g'$ deux métriques riemanniennes compatibles avec la structure conforme de $M$. Les structures presque complexes naturelles associées aux connexions riemanniennes de $g$ et de $g'$ coïncident sur $\calh(T(M))$.}\\

\noindent\underbar{Démonstration}. 
Il suffit de vérifier, en vertu du lemme 6.4, que la différence $K=\nabla'-\nabla$ des connexions riemanniennes $\nabla'$ de $g'=\exp(2\varphi).g$ ($\varphi\in C^\infty(M))$ et $\nabla$ de $g$ vérifie :
\[
g(Z,K(JX,JY))-g(Z,JK(JX,Y))-g(Z,JK(X,JY))-g(Z,K(X,Y)=0
\]
\[
\forall X, Y, Z\in T_x(M),\>\> \forall J\in \calh(T_x(M)),\>\> \forall x\in M
\]
Or on a :
\[
g(Z,K(X,Y))=g(Z,X)\langle Y,d\varphi\rangle+g(Z,Y)\langle X,d\varphi\rangle-g(X,Y)\langle Z,d\varphi\rangle
\]
et l'identité cherchée découle de $J^2=-1$ et $g(JX,JY)=g(X,Y)$.~$\square$

\subsection{Remarque}

Si $J\in \cali(T(M))$ n'est pas dans $\calh(T(M))$ (i.e. pas isométrique pour $g$), l'identité précédente n'est pas satisfaite; autrement dit, les structures presque complexes naturelles associées à $\nabla'$ et à $\nabla$ {\bf ne coïncident pas sur $\cali(T(M))$ tout entier}.\\
Par la suite nous considérerons $\calh(T(M))$ comme une variété presque complexe pour la {\sl structure presque complexe naturelle} définie ci-dessus.\\
Si $M$ est conforme plate, $\calh(T(M))$ est une variété complexe.

\subsection{Proposition}

{\sl La variété complexe $\calh(T(S^{2\ell}))$ des structures complexes au-dessus de la sphère $S^{2\ell}$ compatibles avec la structure conforme plate usuelle de $S^{2\ell}$ est isomorphe à la variété algébrique complexe $\calh(\mathbb R^{2\ell+2})$ des structures complexes isométriques sur l'espace euclidien de dimension $2\ell+2:\calh(T(S^{2\ell}))=\calh(\mathbb R^{2\ell+2})$.}\\

\noindent \underbar{Démonstration}. Un élément $J$ de $\calh(\mathbb R^{2\ell+2})$ est canoniquement une matrice réelle $(2\ell+2)\times (2\ell+2)$ orthogonale antisymétrique: on peut donc l'écrire $J=\left(\begin{array}{cc} \tilde J & u\\ - u^t & 0\end{array}\right)$
où $u\in \mathbb R^{2\ell+1}$ est unitaire, i.e. $u\in S^{2\ell}$ et où $\tilde J\in \calh(T_u(S^{2\ell}))$ i.e. $\tilde J u=0=\tilde J^t u$ et $\tilde J\restriction (\{u\}^\perp_{\mathbb R^{2\ell+1}}=T_u(S^{2\ell}))$ est orthogonale antisymétrique. Cette correspondance établit une bijection entre $\calh(\mathbb R^{2\ell+2})$ et $\calh(T(S^{2\ell}))$ dont on vérifie l'holomorphie en utilisant le fait que la connexion riemannienne de $S^{2\ell}$ est induite par le plongement isométrique $S^{2\ell}\subset \mathbb R^{2\ell+1}\subset \mathbb R^{2\ell+2}$.~$\square$\\
On a donc une fibration
\[
\calh(\mathbb R^{2\ell+2})\underset{\calh(\mathbb R^{2\ell})}\longrightarrow S^{2\ell}
\]
de $\calh(\mathbb R^{2\ell+2})$ par $\calh(\mathbb R^{2\ell})$ au-dessus de $S^{2\ell}$ dont les fibres ($\simeq \calh(\mathbb R^{2\ell})$) sont des sous-variétés complexes de $\calh(\mathbb R^{2\ell+2})$ ($=\calh(T(S^{2\ell})$). La suite exacte d'homotopie correspondante implique le corollaire suivant.

\subsection{Corollaire}

{\sl Les groupes d'homotopie $\pi_n(\calh(\mathbb R^{2\ell}))$ se stabilisent de la manière suivante} :
\[
\pi_n(\calh(\mathbb R^{2\ell}))= \pi_n(\calh(\mathbb R^{\infty})),\>\> \forall \ell>\frac{n+1}{2} \geq 1.
\]
Ceci résulte de $\pi_{n+1}(S^{2\ell})=\pi_n(S^{2\ell})=\{0\}$ pour $\ell>\frac{n+1}{2}\geq 1$. D'autre part $\pi_0(\calh(\mathbb R^{2\ell}))$ est un ensemble à deux éléments puisque $\calh(\mathbb R^{2\ell})$ a deux composantes connexes correspondant aux deux orientations de $\mathbb R^{2\ell}$ pour $\ell\geq 1$.\\
On a évidemment $\calh_+(T(S^{2\ell}))=\calh_+(\mathbb R^{2\ell+2})$.\\
En ce qui concerne l'homologie, la proposition 8.4 a le corollaire suivant.

\subsection{Corollaire}

{\sl L'homologie $H(\calh(\mathbb R^{2\ell}))$ de $\calh(\mathbb R^{2\ell})$ est identique à l'homologie\linebreak[4] $H(\prod^{k=\ell-1}_{k=0}S^{2k})$ du produit topologique des sphères $S^{2k}$ pour $0\leq k\leq \ell-1$; de même $H(\calh_+(\mathbb R^{2\ell}))=H(\prod^{\ell-1}_{k=1} S^{2k})$, (on identifie $S^0$ à l'ensemble à deux élements).}\\

\noindent\underbar{Démonstration}. La sphère $S^{2\ell}$ admet une décomposition cellulaire $S^{2\ell}\sim \mathbb R^{2\ell}\cup \mathbb R^0$; ceci implique que l'on a $\calh(\mathbb R^{2\ell+2})\sim \calh(\mathbb R^{2\ell})\times \mathbb R^{2\ell}\cup\calh(\mathbb R^{2\ell})$. Il en résulte, par récurrence sur $\ell$, que $\calh(\mathbb R^{2\ell})$ admet une décomposition cellulaire dans laquelle on a seulement des cellules de dimensions paires et que, si $p_\ell(t)$ désigne le polynôme de Poincaré de $\calh(\mathbb R^{2\ell})$ on a : $p_{\ell+1}(t)=p_\ell(t). (t^{2\ell}+1)$. Il en résulte $p_\ell(t)=\prod^{\ell-1}_{k=0}(1+t^{2k})$ ce qui est le polynôme de Poincaré du produit $\prod^{\ell-1}_{k=0}S^{2k}$. L'homologie de $\calh(\mathbb R^{2\ell})$ est donc libre et égale à celle de $\prod^{\ell-1}_{k=0}S^{2k}$.~$\square$

\subsection{Exemples}

On a :
\[
\calh_+(T(S^2))=\calh_+(\mathbb R^4)=P_1(\mathbb C)
\]
\[
\calh_+(T(S^4))=\calh_+(\mathbb R^6)=P_3(\mathbb C)
\]
Pour l'homologie, on a : $H_n(\calh_+(\mathbb R^4))=H_n(P_1(\mathbb C))=H_n(S^2)$ et $H_n(\calh_+(\mathbb R^6)=H_n(P_3(\mathbb C))=H_n(S^2\times S^4)$.\\
Pour l'homotopie le corollaire 8.5 done :
\[
\begin{array}{llllllll}
\pi_1(\calh_+(\mathbb R^4)) & = & \pi_1(P_1(\mathbb C)) & = & \{0\} & = & \pi_1(\calh_+(\mathbb R^{2\ell})), & \forall \ell\geq 2\\
\pi_2(\calh_+(\mathbb R^4)) & = & \pi_2(P_1(\mathbb C)) & = & \mathbb Z & = & \pi_2(\calh_+(\mathbb R^{2\ell})), & \forall \ell\geq 2\\
\pi_3(\calh_+(\mathbb R^6)) & = & \pi_3(P_3(\mathbb C)) & = & \{0\} & = & \pi_3(\calh_+(\mathbb R^{2\ell})), & \forall \ell\geq 3\\
\pi_4(\calh_+(\mathbb R^6)) & = & \pi_4(P_3(\mathbb C)) & = & \{0\} & = & \pi_4(\calh_+(\mathbb R^{2\ell})), & \forall \ell\geq 3
\end{array}
\]
On remarquera que $\pi_3(\calh_+(\mathbb R^4))=\pi_3(P_1(\mathbb C))=\mathbb Z$ alors que $\pi_3(\calh_+(\mathbb R^{2\ell}))$ est nul pour $\ell\geq 3$, de sorte que, dans ce cas la borne donnée par le corollaire 8.5, (pour $\pi_3$), est saturée.

\subsection{Formes horizontales pures et auto-dualité}

Soit $M$ une variété orientée de dimension 2$\ell$ munie d'une structure conforme. $\calh(T(M))$, $\calh_+(T(M))$ et $\calh_-(T(M))$ sont des variétés presque complexes et on a, par conséquent, une décomposition des formes extérieures sur ces espaces en composantes pures de type $(p,q)$; nous désignerons par $P^{p,q}$ les projecteurs correspondants. $\calh(T(M)), \calh_+(T(M))$ et $\calh_-(T(M))$ sont, d'autre part, des fibrés sur $M$; {\sl nous désignerons les projections de ces fibrés par} $\pi, \pi_+$ et $\pi_-$ respectivement. Il résulte de la remarque 7.8.c)  que pour les\linebreak[4] $\ell$-formes sur $M$, la dualité de Hodge correspondant à une métrique compatible avec la structure conforme de $M$ ne dépend en fait que la structure conforme de $M$; {\sl on a donc une opération $\ast$ naturelle sur les $\ell$-formes}. Une conséquence immédiate du lemme 7.7 est que, pour une $\ell$-forme $\Omega$ sur $M$, on a équivalence entre $\Omega+i^\ell \ast \Omega=0$, (resp. $\Omega-i^\ell\ast \Omega$), et l'annulation de la composante pure de type $(0,\ell)$ de son image inverse $\pi^\ast_+(\Omega)$, (resp. $\pi^\ast_-(\Omega))$, par $\pi_+$, (resp. $\pi_-$).\\
Nous voulons généraliser l'énoncé précédent. Rappelons qu'une $p$-forme à valeurs vectorielles définie sur un fibré est dite {\sl horizontale} si elle s'annule dès que l'un des $(p)$ vecteurs auxquels on l'applique est vertical i.e. tangent à la fibre; pour tout point du fibré la projection de ce fibré induit un isomorphisme des formes horizontales sur les formes au point correspondant de la base. Les opérations algébriques définies pour les formes sur la base d'un fibré se relèvent ainsi aux formes horizontales sur le fibré.\\
Soit $P\rightarrow M$ un fibré sur $M$ et $E$ un espace vectoriel complexe. L'opération $\ast$ définie pour les $\ell$-formes sur $M$ à valeurs dans $E$ se relève en une opération, encore notée $\ast$, pour les $\ell$-formes horizontales sur $P$ à valeur dans $E$. Si $\Omega$ est une telle $\ell$-forme horizontale sur $P$ il lui correspond, par image inverse, une $\ell$-forme $\pi^\ast_+(\Omega)$, (resp. $\pi^\ast_-(\Omega))$ sur le fibré $\pi^\ast_+(P)$, (resp. $\pi_-^\ast (P))$, image inverse de $P$ par $\pi_+$ resp. $\pi_-$; ces formes sont horizontales et l'on peut parler de leurs composantes de type $(r,\ell-r)$, etc. Avec ces conventions, on a :

\subsection{Lemme}

{\sl Soit $P$ un fibré sur $M$, $E$ un espace vectoriel complexe, et $\Omega$ une $\ell$-forme horizontale sur $P$ à valeurs dans $E$. On a équivalence entre $\Omega+i^\ell\ast \Omega=0$ (resp. $\Omega - i^\ell\ast \Omega=0$), et l'annulation de la composante pure de type $(0,\ell)$ de $\pi^\ast_+(\Omega)$ (resp. $\pi^\ast_-(\Omega)$).}\\
Ce lemme est encore une conséquence immédiate du lemme 7.7.\\

Le lemme précédent, le théorème 4.6 et le théorème 4.8 conduisent, pour la dimension 4 ($\ell=2$) à la forme générale du théorème de Atiyah et Ward relatif à la transformation de Penrose pour les instantons.

\subsection{Théorème (transformation de Penrose pour les instantons)}

{\sl Soit $M$ une variété orientée de dimension 4 munie d'une structure conforme telle que $\calh_+(T(M))$, (resp. $\calh_-(T(M))$), soit une variété complexe; soient $G$ un groupe de Lie complexe et $U$ un sous-groupe de Lie réel de $G$ dont l'algèbre de Lie est une forme réelle de celle de $G$. La projection $\pi_+:\calh_+(T(M))\rightarrow M$, (resp. $\pi_-:\calh_-(T(M))\rightarrow M$), induit, par image inverse, une bijection entre les $U$-fibrés principaux sur $M$ munis de connexions dont les formes de courbure satisfont $\Omega=\ast\Omega$, (resp. $\Omega=-\ast\Omega$), et les $U$-fibrés principaux sur $\calh_+(T(M))$, (resp. $\calh_-(T(M)))$ tels que les $G$-fibrés principaux correspondant par extension du groupe structural sont des fibré holomorphes sur $\calh_+(T(M))$, (resp. $\calh_-(T(M))$), holomorphiquement triviaux sur les fibres de la fibration $\pi_+=\calh_+(T(M))\rightarrow M$ (resp. $\pi_-:\calh_-(T(M))\rightarrow M)$}.\\

\noindent \underbar{Démonstration}. Soit $R$ un fibré principal sur $M$ de groupe structural $U$ muni d'une connexion $\omega$ dont la courbure $\Omega$ satisfait $\Omega=\ast\Omega$, soit $Q=\pi^\ast_+(R)$ l'image inverse par $\pi_+$ de $R$ et $P$ le $G$-fibré principal sur $\calh_+(T(M))$ contenant $Q$. L'image inverse de $\omega$ par $\pi_+$ définie, par extension à $P$, une connexion sur $P$ telle que sa courbure a, en vertu du lemme 8.9, une composante pure de type (0,2) qui est nulle; il résulte de 4.6 que $P$ a une structure unique de fibré holomorphe pour laquelle cette connexion est pure de type (1,0). Par construction, la courbure de cette connexion sur $P$ s'annule sur les fibres de la fibration $\pi_+:\calh_+(T(M))\rightarrow M$; il en résulte (en appliquant 4.6 aux fibres) que $P$ est holomorphiquement trivial sur les fibres de cette fibration.\\
Inversement si $P$ est un fibré principal holomorphe sur $\calh_+(T(M))$ de groupe structural $G$ et si $Q$ est un sous-fibré réel de $P$ de groupe structural $U$, on a, d'après 4.8, une connexion unique sur $Q$ dont l'extension à $P$ est pure de type (1,0). Si $P$ est holomorphiquement trivial sur les fibres de $\pi_+:\calh_+(T(M))\rightarrow M$, cette connexion est plate sur les fibres de $\pi_+:\calh_+(T(M))\rightarrow M$. Il en résulte que l'on peut ``pousser la situation sur $M$"; i.e. on a un fibré principal $R$ sur $M$ de groupe structural $U$ et une connexion sur $R$ tels que $Q$ est l'image inverse par $\pi_+$ de $R$ et que la connexion de $Q$ est l'image inverse par $\pi_+$ de celle de $R$. Comme, par construction, la composante pure de type (0,2) de la courbure de la connexion de $Q$ étendue à $P$ est nulle, le lemme 8.9 implique que la courbure $\Omega$ de la connexion de $R$ vérifie $\Omega=\ast\Omega$.\\
Aux isomorphismes de $U$-fibrés avec connexions sur $M$ correspond, par la construction ci-dessus les isomorphismes de $G$-fibrés holomorphes sur $\calh_+(T(M))$ avec la structure supplémentaire décrite dans l'énoncé du théorème. On raisonne de même avec les signes opposés.~$\square$\\
On retrouve le cas le plus usuel en prenat $U=SU(2)$, $G=SL(2,\mathbb C)$ et $M=S^4$, ce qui implique $\calh_+(T(M))\simeq P_3(\mathbb C)\simeq \calh_-(T(M))$.

\subsection{Intégrabilité : cas de la dimension 4}

Soit $M$ une variété orientée de dimension 2$\ell$ munie d'une structure conforme. D'après le lemme 6.5, la torsion de la structure presque complexe naturelle de $\calh(T(M))$ s'annule en $J\in \calh(T_x(M))$ ($x\in M$), si et seulement si on a, $\forall X,Y \in T_x(M)$
\[
[W(X,Y)+JW(X,JY)+JW(JX,Y)-W(JX,JY), J]=0,
\]
où $[\bullet, \bullet]$ est le commutateur des endomorphismes de $T_x(M)$ correspondant et où $W$ est la courbure conforme de Weyl. Dans le cas où $\ell=2$ $W$ a une décomposition $W=W_++W_-$ avec $W_\pm=\pm\ast W_\pm$. Si $W_-=0$, (resp. $W_+=0$) la {\sl structure conforme} de $M$ est dite {\sl auto-duale} (resp. {\sl anti-auto-duale}). La symétrie de $W$ implique, avec nos conventions d'orientation $[W,J]=[W_-,J]$ si $J\in \calh_+(T(M))$ et $[W,J]=[W_+,J]$ si $J\in \calh_-(T(M))$; $(\calh_+(\mathbb R)$ est l'ensemble des matrices orthogonales anti-symétriques qui sont anti-auto-duales avec nos conventions d'orientation). Il en résulte que si $W_-=0$ (resp. $W_+=0$), la structure presque complexe naturelle de $\calh_+(T(M))$, (resp. $\calh_-(T(M))$, est intégrable; la réciproque est également vraie et résulte d'un argument d'irréductibilité. (Dans les cas $\ell\not=2$, $\calh_+(T(M))$ est complexe si et seulement si $M$ est conforme plate (ce qui implique que $\calh_-(T(M))$ est complexe)).

\subsection{Définition}

Soit $M$ une variété de dimension $2\ell$ munie d'une structure conforme. Nous désignerons par $\wedge^{p,q}_H\calh(T(M))$ la restriction à $\calh(T(M))\subset \cali(T(M))$ de $\wedge^{p,q}_H\cali(T(M))$; il résulte de 6.9 que ses éléments qui seront appelés {\sl formes horizontales pures de type $(p,q)$ sur $\calh(T(M))$} sont des formes pures de type $(p,q)$ pour la structure presque complexe naturelle de $\calh(T(M))$.

\subsection{Remarque}

En utilisant  7.9, on construit facilement une rétraction par déformation de $\cali(T(M))$ sur $\calh(T(M))$. Ceci implique en particulier que la classification des fibré sur $\cali(T(M))$ est la même que sur $\calh(T(M))$.

\newpage

\section{Relation avec la théorie des spineurs}

Dans ce paragraphe $E$ désigne un espace euclidien de dimension $2\ell$ et nous utiliserons les notations et conventions du paragraphe 7. Par la suite, il sera commode d'adopter la définition suivante concernant les spineurs associés à $E$.

\subsection{Spineurs}
On appellera {\sl espace de spineurs associé à $E$} un espace hilbertien $S$ dans lequel on a une représentation hermitienne irréductible $\gamma$ de l'algèbre de Clifford $C(E^\ast)$ du dual $E^\ast$ de $E$. Les éléments de $S$ seront appelés {\sl spineurs}.\\
Autrement dit, on a une application linéaire (réelle) $\gamma$ de $E^\ast$ dans l'espace des opérateurs hermitiens sur $S$ telle que $\gamma (\omega_1)\gamma (\omega_2)+\gamma(\omega_2)\gamma(\omega_1)=2(\omega_1,\omega_2)\bbbone$ $\forall \omega_1,\omega_2\in E^\ast$ et telle que le commutant $\gamma(E^\ast)'$ de $\gamma(E^\ast)$ dans l'algèbre $\End(S)$ se compose des multiples de l'identité $\bbbone$ de $S$. Il revient au même de dire que l'on a une application linéaire complexe $\gamma:E^\ast_c\rightarrow \End(S)$ satisfaisant à ce qui précède et à :
$\gamma(\omega)^\ast=\gamma(\bar\omega)$, $\forall \omega\in E^\ast_c$.\\
Il est classique, et facile à démontrer, que ces conditions spécifient $(\gamma, S)$ à une équivalence unitaire près. On a : $\dim(S)=2^\ell$.\\
Si $\psi$ est un spineur non nul, l'ensemble des $\omega\in E^\ast_c$ satisfaisant à $\gamma(\omega)\psi=0$ est un sous-espace totalement isotrope $I(\psi)$ de $E^\ast_c$ pour $(\bullet,\bullet)$. On est amené à la définition suivante due à E. Cartan.

\subsection{Spineurs simples}

Un spineur non nul $\psi$ sera appelé {\sl spineur simple} ({\sl spineur pur} dans la terminologie de C. Chevalley) si $I(\psi)$ est maximal isotrope, i.e. si $\dim(I(\psi))=\ell$. L'ensemble des spineurs simples sera désigné par $\fracS$.\\
L'ensembe $P(\fracS)$ des directions de spineurs simples est une sous-variété algé\-brique complexe de l'espace projectif complexe $P(S)\simeq P_{2^\ell-1}(\mathbb C)$.\\
Il résulte de la proposition 7.2 que à tout $\psi\in \fracS$ correspond un $J\in \calh(E)$ unique pour lequel on a $I(\psi)=\wedge^{1,0}E^\ast_J$; $J$ ne dépend en fait que de la direction $[\psi]$ de $\psi$. On définit ainsi, ($[\psi]\mapsto J=i([\psi])$), une application algébrique $i$ de $P(\fracS)$ dans $\calh(E)$; $i:P(\fracS)\rightarrow \calh(E)$. Par ailleurs, et c'est le point de départ de l'approche géométrique à la théorie des spineurs de E. Cartan, les sous-espaces isotropes de dimension $\ell$ de $E^\ast_c$ sont paramétrés par $P(\fracS)$. On peut donc résumer la situation par le théorème suivant.

\subsection{Théorème}

{\sl L'application $i:P(\fracS)\rightarrow \calh(E)$ est un isomorphisme de variétés complexes}.\\

Par la suite, lorsque cela n'entraîne pas de confusion, nous identifierons $P(\fracS)$ à $\calh(E)$ par cet isomorphisme.

\subsection{Définition : le fibré naturel $L$}

La restriction à $P(\fracS)$ du fibré tautologique sur $P(S)$ est un fibré vectoriel complexe holomorphe hermitien de rang 1 sur $\calh(E)$ que nous noterons $L\rightarrow \calh (E$) et qui sera appelé {\sl fibré naturel de rang 1 sur $\calh(E)$}, ($L=\fracS \cup \{\text{section nulle}\}$). Nous reviendrons ultérieurement sur l'étude des propriétés de $L$ et sur ses relations avec le fibré canonique de $\calh(E)$ et le fibré $\wedge^{\ell,0}E^\ast$ qui sont tous deux également des fibrés holomorphes hermitiens de rang 1 sur $\calh(E)$.

\subsection{Proposition}

{\sl Le fibré hermitien $(\oplus_{k\geq 0}\wedge^{0,k}E^\ast)\otimes L$ est isomorphe au fibré trivial de fibre $S$ sur $\calh(E)$ par un isomorphisme tel que, si $\Phi$ désigne l'application à valeurs dans $S$ correspondante, on a} :
\[
\Phi(\omega\wedge\varphi)=1/{\sqrt 2} \gamma(\omega)\Phi(\varphi),\>\> \forall \omega \in \wedge^{0,1}E^\ast_J,\>\> \forall \varphi\in(\oplus_{k\geq0}\wedge^{0,k}E^\ast_J)\otimes L_J,\>\> \forall J\in \calh(E)
\]

\noindent\underbar{Démonstration}. Soit $J\in \calh(E)$, on a
\[
\gamma\left(\frac{\omega_1}{\sqrt{2}}\right)\gamma\left(\frac{\omega_2}{\sqrt{2}}\right)+\gamma\left(\frac{\omega_2}{\sqrt{2}}\right)\gamma\left(\frac{\omega_1}{\sqrt{2}}\right)= 0,
\]
et
\[
\gamma\left(\frac{\omega_1}{\sqrt{2}}\right)^\ast \gamma\left(\frac{\omega_2}{\sqrt{2}}\right)+\gamma\left(\frac{\omega_2}{\sqrt{2}}\right)\gamma\left(\frac{\omega_1}{\sqrt{2}}\right)^\ast=\langle\omega_1\vert\omega_2\rangle\bbbone,\>\>\>  \forall \omega_1,\omega_2\in \wedge^{1,0}E^\ast_J.
\]
Il en résulte que si $\psi\in L_J$, $\lVert \gamma\left(\frac{\omega_1}{\sqrt{2}}\right)\dots \gamma \left(\frac{\omega_k}{\sqrt{2}}\right)\psi\rVert=
\lVert \omega_1\wedge\dots \wedge \omega_k\lVert \bullet\rVert\psi\rVert$ et on définit $\Phi$ par $\Phi(\omega_1,\wedge\dots\wedge\omega_k)\otimes \psi)=(1/\sqrt{2})^k\gamma(\omega_1)\dots \gamma(\omega_k)\psi$, $\forall \omega_1,\dots, \omega_k\in \wedge^{0,1}E^\ast_J$.~$\square$

\subsection{Remarque}

Les $\gamma\left(\frac{\omega}{\sqrt{2}}\right)$, $\omega\in \wedge^{1,0}E^\ast_J$ et leurs adjoints satisfont aux relations d'anti-commutation canonique et sont irréductibles; la direction correspondante $L_J\in P(\fracS)$ est la direction unique du vide associé. Ceci est une manière de montrer la bijectivité de l'application $i$ du théorème 9.3. L'action des rotations sur $\calh(E)$ correspond aux transformations de Bogoliubov, les spineurs simples sont tous les vides possibles et le contenu de 9.5 est la représentation de Fock familière.\\
Aucune des terminologies ``spineur simple" (E. Cartan) et ``spineur pur" (C. Chevalley) n'est satisfaisante et cette remarque montre que la terminologie {\sl spineur de Fock} serait bien mieux adaptée.

\subsection{Semi-spineurs}

Supposons $E$ orienté. L'action de $SO(E)$ sur $E^\ast$ correspond à l'action adjointe de son revêtement d'ordre 2, $\text{Spin}(E)$, sur $\gamma(E^\ast)$; $\text{Spin}(E)$ est ainsi un sous-groupe de $SU(S)$ et $S$ est somme directe de deux sous-espaces isomorphes $S_+$ et $S_-$ irréductibles pour l'action de spin$(E)$. Les éléments de ces deux sous-espaces sont appelés {\sl semi-spineurs}. L'action projective de Spin$(E)$ qui est une action de $SO(E)$ laisse invariant $P(\fracS)$ et correspond sur $P(\fracS)$ à l'action de $SO(E)$ sur $\calh(E)$. Comme, pour cette action, $\calh(E)$ est composé des deux orbites $\calh_+(E)$ et $\calh_-(E)$, {\sl les spineurs simples sont nécessairement des semi-spineurs}. On a $\fracS=\fracS_+ \cup \fracS_-$ avec $P(\fracS_\pm)=\calh_\pm(E)$ et on choisira l'indexation par $\pm$ dans la décomposition $S=S_+\oplus S_-$ de manière à avoir $\fracS_\pm \subset S_\pm$; $S_+$, (resp. $S_-$), est irréductible sous l'action de Spin$(E)$, ce qui implique que $S_+$, (resp. $S_-$), est engendré par $\fracS_+$, (resp. $\fracS_-$).\\

L'application $\Phi$ de 9.5 induit les isomorphismes
\[
\begin{array}{ll}
\Phi_J:(\oplus_{p\geq 0}\wedge^{0,2p}E^\ast_J)\otimes L_J\simeq S_+, & \forall J\in \calh_+(E)\\
\Phi_J:(\oplus_{p\geq 0}\wedge^{0,2p+1}E^\ast_J)\otimes L_J\simeq S_-, & \forall J\in \calh_+(E)\\
\Phi_J:(\oplus_{p\geq 0}\wedge^{0,2p}E^\ast_J)\otimes L_J\simeq S_-, & \forall J\in \calh_-(E)\\
\Phi_J:(\oplus_{p\geq 0}\wedge^{0,2p+1}E^\ast_J)\otimes L_J\simeq S_+, & \forall J\in \calh_-(E)
\end{array}
\]

\subsection{Exemples en basse dimension}

Soit $J\in \calh_+(E)$, alors $\Phi^{-1}_J$ composé avec la projection 
$$(\oplus_{p\geq0} \wedge^{0,2p}E^\ast_J)\otimes L_J \rightarrow (\mathbb C\oplus \wedge^{0,2}E^\ast_J)\otimes L_J$$
est surjective de $\fracS_+$ sur $(\mathbb C\oplus \wedge^{0,2}E^\ast_J)\otimes L_J\backslash \{0\}$ comme on le voit en utilisant le fait que l'algèbre de Lie de $\Spin(E)$ est engendré par les $[\gamma(\omega_1), \gamma(\omega_2)]$ avec $\omega_1,\omega_2\in E^\ast$ et que $\fracS_+$ est invariant dans l'action de $\Spin(E)$. Il en résulte que pour $\ell\leq 3$ on a $\fracS_+=S_+$, i.e. les spineurs simples sont tous les semi-spineurs non nuls. Ceci permet de retrouver
\[
\begin{array}{lllllll}
\calh_+(\mathbb R^2) & = & P(\fracS_+) & = & P(S_+) & \simeq &P_0(\mathbb C)\\
\calh_+(\mathbb R^4) & = & P(\fracS_+) & = & P(S_+) & \simeq &P_1(\mathbb C)\\
\calh_+(\mathbb R^6) & = & P(\fracS_+) & = & P(S_+) & \simeq &P_3(\mathbb C)
\end{array}
\]
Mais pour $\ell\geq 4$ l'inclusion $\fracS_+ \subset S_+$ est stricte car $\dim(S_+)=2^{\ell-1}$ tandis que $\dim(\fracS_+)=\dim(\calh_+(\mathbb R^{2\ell}))+1=\frac{\ell(\ell-1)}{2}+1$ et que l'on a $2^{\ell-1}>\frac{\ell(\ell-1)}{2}+1$, $\forall\ell \geq 4$.\\
Revenons à l'étude du fibré $L\rightarrow \calh(E)$.

\subsection{Proposition}

{\sl Pour $\ell\geq 2$, $-c_1(L\restriction \calh_+(\mathbb R^{2\ell})), $ (resp. $-c_1(L \restriction\calh_-(\mathbb R^{2\ell}))$), est le générateur de $H^2(\calh_+(\mathbb R^{2\ell})\simeq \mathbb Z$ (resp. $H^2(\calh_-(\mathbb R^{2\ell})$); $c_1(\bullet)$ désignant la première classe de Chern}.\\

\noindent\underbar{Démonstration}. La démonstration du corollaire 8.6 montre que dans la fibration $\calh_+(\mathbb R^{2\ell+2})\rightarrow S^{2\ell}$, l'injection d'une fibre $\calh_+(\mathbb R^{2\ell})\subset \calh_+(\mathbb R^{2\ell+2})$ induit des isomorphismes
$H_n(\calh_+(\mathbb R^{2\ell}))\simeq H_n(\calh_+(\mathbb R^{2\ell+2}))$ pour $n<2\ell$. La restriction du fibré $L$ de $\calh_+(\mathbb R^{2\ell+2})$ à la fibre $\calh_+(\mathbb R^{2\ell})$ étant le fibré $L$ naturel sur $\calh_+(\mathbb R^{2\ell})$, on a le résultat de la proposition par récurrence sur $\ell$ à partir de $\ell=2$ où $L$ est le fibré tautologique de $P_1(\mathbb C)$ ($\simeq \calh_+(\mathbb R^4)$).~$\square$\\

Le cas $\ell=1$ est trivial car $\calh_\pm(\mathbb R^2)=\{I_\pm\}$. Il résulte de la proposition précédente que tout fibré complexe de rang 1 sur $\calh_\pm(E)$ (on suppose $E$ orienté) est isomorphe, comme fibré $C^\infty$ ou $C^0$, à une puissance tensorielle du fibré naturel ou de son dual sur $\calh_\pm(E)$. Dans le cas de $\wedge^{\ell,0}E^\ast$ et du fibré canonique, nous voulons aller plus loin et les identifier en tant que fibrés holomorphes hermitiens. Nous nous servirons du lemme suivant.

\subsection{Lemme}

{\sl Soit $\calo$ un ouvert contractile de $\calh(E)$ sur lequel on a une section $\varphi$ du fibré naturel de rang 1 telle que $\lVert\varphi\rVert=1$ et une section ($\omega^1,\dots,\omega^\ell$) du fibré des repères orthonormés de $\wedge^{1,0}E^\ast$. Alors il existe une fonction réelle $\theta$ sur $\calo$ tel que l'on ait 
$$d\varphi=1/2(\sum_{m<n}(\omega^m,d\omega^n)\gamma(\bar\omega^m)\gamma(\bar\omega^n)+\sum_k(\bar\omega^k,d\omega^k)+id\theta)\varphi$$
où $d\varphi$ (resp. $d\omega^k)$ est la différentielle de $\varphi$ comme fonction à valeurs dans $S$ (resp. dans $E^\ast_c)$.}\\

\noindent \underbar{Démonstration}. On peut écrire $d\varphi=\frac{1}{4}\sum_{m,n} \alpha_{m,n} \gamma(\bar\omega^m)\gamma(\bar\omega^n)\varphi+\frac{1}{2}\alpha\varphi$, où $\alpha_{m,n}$ et $\alpha$ sont des 1-formes sur $\calo$ car un spineur simple se déduit d'un autre par une transformation de $\Spin(E)$ combiné à une homotétie, $\lVert\varphi\rVert=1$ entraine $\alpha+\bar\alpha=0$. En prenant la différentielle de $\gamma(\omega^m)\varphi=0$ et en tenant compte de $\gamma(\omega^n)\varphi=0$, il vient : $\alpha_{m,n}=(\omega^m,d\omega^n)$. En écrivant alors $d^2\varphi=0$, on obtient $d\alpha+\sum_{m,n}(\omega^m,d\omega^n)\wedge (\bar\omega^m,d\bar\omega^n)=0$, ce qui peut s'écrire\linebreak[4] $d(\alpha-\sum_k(\bar\omega^k,d\omega^k))=0$. Il en résulte $\alpha=\sum(\bar\omega^k,d\omega^k)+id\theta$ ce qui achève la démonstration.~$\square$\\

\subsection{Théorème}

{\sl Comme fibrés vectoriels complexes holomorphes hermitiens de rang 1 sur $\calh(E)$, $\wedge^{\ell,0}E^\ast$ est isomorphe à $L\otimes L$ et le fibré canonique $\wedge^{\frac{\ell(\ell-1)}{2},0}T^\ast(\calh(E))$ de $\calh(E)$ est isomorphe à $L^{\otimes2(\ell-1)}$. Si $L'$ est un fibré holomorphe hermitien de rang 1 sur $\calh'E)$ tel que $L^{\prime \otimes p}$ est isomorphe à $L^{\otimes p}$ pour un entier $p\geq 1$, alors $L'$ est isomorphe à $L$.}\\

\noindent\underbar{Démonstration}. Un isomorphisme de fibrés holomorphes hermitiens est un isomorphisme de fibrés différentiables hermitiens interchangeant leurs connexions canoniques de fibrés holomorphes hermitiens, i.e. leurs connexions hermitiennes de type (1,0), (c'est un cas particulier des théorèmes 4.6 et 4.8).\\

Reprenons les notations du lemme 9.10. La connexion canonique de $L$ dans la section $\varphi$ est $\langle\varphi \vert d\varphi\rangle=\frac{1}{2}(\sum_k(\bar\omega^k,d\omega^k)+id\theta)$. Il en résulte que la connexion de $L\otimes L$ dans la section unitaire $e^{-i\theta}\varphi\otimes \varphi$ est $\sum_k(\bar\omega^k,d\omega^k)$ ce qui est la connexion de $\wedge^{\ell,0}E^\ast$ dans la section $\omega^1\wedge\dots\wedge\omega^\ell$. Soit $\calo_\alpha$ un recouvrement de $\calh(E)$ par des ouverts contractiles et soit $\lambda_\alpha$, (resp. $\mu_\alpha$), des sections unitaires de $\wedge^{\ell,0}E^\ast$, (resp. $L\otimes L$) au-dessus des $\calo_\alpha$ telles que pour tout $\alpha$, la forme de connexion $\omega_\alpha$ de $\wedge^{\ell,0}E^\ast$ dans $\lambda_\alpha$ soit égale à celle de $L^{\otimes 2}$ dans $\mu_\alpha$. Soient $g_{\alpha\beta}$ et $\ell_{\alpha\beta}=g_{\alpha\beta} u_{\alpha\beta}$ les fonctions de transitions correspondantes ($g_{\alpha\beta}=\langle\lambda_\alpha\vert\lambda_\beta\rangle$ sur $\calo_\alpha\cap \calo_\beta$); on a $du_{\alpha\beta}=0$ puisque les formes de connexion sont les mêmes dans les 2 systèmes de section. $u_{\alpha\beta}$ définit donc un cocycle à valeurs dans $U(1)$ et on peut l'écrire $u_{\alpha\beta}=u^{-1}_\alpha u_\beta$ où $u_\alpha\in U(1)$ puisque $H_1(\calh(E))=0$ (corollaire 8.6). Les sections $u^{-1}_{\alpha}\mu_\alpha$ de $L\otimes L$ ont alors les mêmes fonctions de transition que les sections $\lambda_\alpha$ de $\wedge^{\ell,0}E^\ast$ et les formes de connexions sont les mêmes dans ces systèmes de sections, ce qui prouve que $\wedge^{\ell,0}E^\ast\simeq L^{\otimes 2}$comme fibrés holomorphes hermitiens.\\

Le théorème 7.5 implique que le fibré canonique est isomorphe à\linebreak[4] $(\wedge^{\ell,0}E^\ast)^{\otimes\ell-1}$ donc à $L^{\otimes 2(\ell-1)}$, comme fibré holomorphe hermitien.\\
Si $L^{\prime \otimes p}\simeq L^{\otimes p}$, ($p\geq 1$), on peut trouver des sections unitaires $\varphi_\alpha$ de $L$ et $\varphi'_\alpha$ de $L'$ au-dessus des $\calo_\alpha$ telles que les formes des connexions canoniques de $L$ et $L'$ soient les mêmes dans ces 2 systèmes de sections et telles que les fonctions de transitions $g_{\alpha\beta}$ et $g'_{\alpha\beta}$ respectives satisfassent $(g_{\alpha\beta})^p=(g'_{\alpha\beta})^p$. On a donc $g'_{\alpha\beta}=g_{\alpha\beta}z_{\alpha\beta}$ où $z_{\alpha\beta}$ est dans le groupe des racines p-ième de l'unité. Comme $H_1(\calh(E))=0$, $z_{\alpha\beta}=z^{-1}_{\alpha}z_\beta$ où les $z_\alpha$ sont des racines p-ième de 1. En remplaçant les $\varphi'_\alpha$ par $z^{-1}_\alpha\varphi'_\alpha$, on a 2 systèmes de sections unitaires sur $L$ et $L'$ ayant les mêmes fonctions de transitions et dans lesquelles les formes de connexions sont identiques. Ceci prouve que $L'\simeq L$ comme fibrés holomorphes hermitiens.~$\square$\\

\subsection{Remarque}

Au départ, on a défini le fibré naturel de rang 1 sur $\calh(E)$, $L$, à partir de la théorie des spineurs; il résulte de 9.11 que, à un isomorphisme près, on peut le définir comme étant la ``racine carrée" de $\wedge^{\ell,0}E^\ast$ par exemple (comme fibré holomorphe hermitien).

\newpage

\section{Structures spinorielles, opérateurs associés de Dirac et transformation de Penrose}

Dans ce paragraphe, $\rho:\Spin(n)\rightarrow SO(n)$ désigne le revêtement connexe (simplement connexe pour $n\geq 3$) d'ordre deux du groupe $SO(n)$, ($1\rightarrow \mathbb Z/2\rightarrow \Spin(n)\stackrel{\rho}{\rightarrow} SO(n)\rightarrow 1$). Bien qu'il soit possible d'être plus général, nous supposerons que $M$ est une variété riemannienne {\sl orientée} de dimension $n$ pour définir les structures spinorielles.

\subsection{Structures spinorielles} 
{\sl Une structure spinorielle sur $M$} est une paire, ($P,\tilde \rho$) où $P=P(M,\Spin(n))$ est un fibré principal sur $M$ de groupe structural $\Spin(n)$ et où $\tilde\rho$ est un homomorphisme de fibrés principaux sur $M$ de $P(M,\Spin(n))$ sur le fibré $R_+(M)$ des repères orthonormés d'orientation positive sur $M$ induisant le revêtement $\rho:\Spin(n)\rightarrow SO(n)$ pour les groupes de structures.\\
Par la suite, nous supposons $M$ de dimension $2\ell$, $E$ désigne $\mathbb R^{2\ell}$ muni de son produit scalaire et de son orientation usuels, $S$ est un espace de spineurs associé à $E=\mathbb R^{2\ell}$ et nous utilisons les notations du paragraphe 9 pour les spineurs simples, les semi-spineurs, etc. En particulier, $\Spin(E)=\Spin(2\ell)\subset SU(2^\ell)=SU(S)$.\\

Soit $(P,\tilde\rho)$ une structure spinorielle sur $M$; on a un fibré vectoriel hermitien $S_P$ sur $M$ de fibre type $S$ associé à $P$ et à l'action de $\Spin(E)$ sur $S$. $S_P$ sera appelé le {\sl fibré des spineurs sur $M$} correspondant à $(P,\tilde\rho)$. De même, on a un {\sl fibré $\fracS_P$ des spineurs simples} sur $M$ de fibre type $\fracS$ qui est un sous-fibré de $S_P$ ($\fracS_P\subset S_P)$. Le fibré $P(\fracS_P)$ des directions de spineurs simples dans $S_P$ ne dépend pas de $(P,\tilde\rho)$; on peut l'identifier (d'après 9.3) au fibré $\calh(T(M))$ qui est un fibré associé à $R_+(M)$ (de fibre type $\calh(E)=\calh(\mathbb R^{2\ell})$ pour l'action de $SO(2\ell)$ sur $\calh(\mathbb R^{2\ell})$). $\fracS_P$ s'identifie donc au fibré des repères (ou des directions) d'un fibré vectoriel complexe de rang 1, $L_P$, sur $\calh(T(M))$ correspondant au fibré naturel de rang 1,  $L\rightarrow \calh(\mathbb R^{2\ell})$. $L_P$ est un fibré hermitien de rang 1
sur $\calh(T(M))$ qui dépend de la structure spinorielle $(P,\tilde\rho)$ sur $M$; nous dirons que $L_P\rightarrow \calh(T(M))$ est {\sl le fibré naturel de rang 1 sur $\calh(T(M))$ correspondnat à la structure spinorielle} $(P,\tilde \rho)$, (ou $P$ lorsque cela n'entraîne pas de confusion), sur $M$.\\
Le théorème 9.11 implique que $\forall x\in M$, les restrictions des fibrés $L_P\otimes L_P$ et $\wedge^{\ell,0}_H \calh(T(M))$ à $\calh(T_(M))$ coïncident (comme fibrés hermitiens de rang 1); il n'est pas difficile de voir que la définition des structures spinorielles implique que $L_P\otimes L_P$ coïncide avec $\wedge^{\ell,0}_H\calh(T(M))$ sur $\calh(T(M))$ tout entier et que l'on a :

\subsection{Théorème}

{\sl La correspondance $P\mapsto L_P(=\eta)$ est une bijection de l'ensemble des structures spinorielles $P$ sur $M$ sur l'ensemble des fibrés hermitiens $\eta$ de rang 1 sur $\calh(T(M))$ satisfaisant $\eta\otimes \eta=\wedge^{\ell,0}_H\calh(T(M))$}.\\

Dans cet énoncé, il est entendu que l'on a muni $\wedge^{\ell,0}_H(T(M))$ de sa structure hermitienne associée à la métrique de $M$ et que, pour écrire l'égalité, on a choisi une fois pour toutes une identification de $L\otimes L$ avec $\wedge^{\ell,0}E^\ast$, (comme on a choisi une fois pour toutes un espace $S$ de spineurs associé à $E=\mathbb R^{2\ell}$).\\

Un corollaire immédiat de 10.2 qui ne se réfère pas à une structure riemannienne particulière sur $M$ est le résultat suivant.

\subsection{Corollaire}
{\sl Soit $M$ une variété de dimension 2$\ell$. Pour que la seconde classe de Stiefel-Whitney de $M$ soit nulle, il faut et il suffit que le fibré $\wedge^{\ell,0}_H\cali(T(M))$ admette une racine carrée. Les racines carrées, (i.e. les classes de fibrés $\eta$ de rang 1 sur $\cali(T(M))$ tels que $\eta\otimes \eta \simeq \wedge^{\ell,0}_H\cali(T(M))$), sont alors classifiées par la cohomologie $H^1(M,\mathbb Z/2)$.}\\

Pour passer de 10.2 à 10.3, nous avons utilisé la remarque 8.13 permettant de remplacer $\cali(T(M))$ par $\calh(T(M))$ (pour une structure riemannienne arbitraire sur $M$) et nous avons utilisé un résultat classique concernant les structures spinorielles.\\
Supposons que $M$ admette une structure presque complexe isométrique, i.e. on a une section  $J$ de $\calh(T(M))\rightarrow M$. Alors, 9.11 implique que dans 10.2 on peut se restreindre à cette section.

\subsection{Corollaire}

{\sl Soit $M$ une variété presque hermitienne de dimension réelle $2\ell$. $M$ admet une structure spinorielle si et seulement si son fibré canonique $\wedge^{\ell,0}T^\ast(M)$ (= restriction de $\wedge^{\ell,0}_H\calh(T(M))$ à la section $J$) admet une racine carrée. Les structures spinorielles sur $M$ sont indexées par ces racines carrées (aux isomorphismes près)}.

\subsection{Exemples}

1) On sait que le fibré canonique de $P_n(\mathbb C)$ est la $(n+1)$-ième puissance tensorielle du fibré tautologique, lequel correspond à la classe $
-1\in H^2(P_n(\mathbb C))=\mathbb Z$; il en résulte que le fibré canonique n'admet pas de racine carrée si $n+1$ est impair et qu'il en admet une unique si $n+1$ est pair. 10.4 implique donc, (le résultat classique), que les espaces projectifs complexes de dimensions paires, $P_{2\ell}(\mathbb C)$, n'admettent pas de structure spinorielle tandis que, $P_{2\ell+1}(\mathbb C)$ en admettent une.\\

2) Dans 8.4, nous avons vu que $\calh(T(S^{2\ell}))$ est identique à $\calh(\mathbb R^{2\ell+2})$; pour le faire, nous avons considéré les plongements successifs $S^{2\ell}\subset \mathbb R^{2\ell+1}\subset \mathbb R^{2\ell+1}\times \mathbb R=\mathbb R^{2\ell+2}$. La manière dont nous avons raisonné dans 8.4 montre que l'on a 
\[
\wedge^{\ell+1,0}(\mathbb R^{2\ell+2})^\ast=\wedge^{\ell,0}_H\calh(T(S^{2\ell}))\otimes \mathbb C = \wedge^{\ell,0}_H\calh(T(S^{2\ell}))
\]
Il en résulte, (en appliquant 9.11), que $S^{2\ell}$ admet une structure spinorielle unique, (ce qui est bien connu), et que le fibré naturel de rang 1 sur $\calh(T(S^{2\ell}))$ correspondant est identique au fibré naturel de rang 1 sur $\calh(\mathbb R^{2\ell+2})$, (défini en 9.4).

\subsection{Spineurs et formes différentielles}

Soit $P$ une structure spinorielle sur $M$ et $S_P$ le fibré des spineurs correspondant. Il résulte des définitions et de la proposition 9.5 que l'image inverse $\pi^\ast(S_P)$ par $\pi:\calh(T(M))\rightarrow M$ de $S_P$ est isomorphe au fibré 
$$(\oplus_{k\geq 0} \wedge^{0,k}_H\calh(T(M)))\otimes L_P$$ où $L_P$ est le fibré naturel de rang 1 sur $\calh(T(M))$ correspondant à $P$. En identifiant ces deux fibrés sur $\calh(T(M))$, on peut considérer que {\sl les champs spinoriels sur $M$ se relèvent en des formes différentielles inhomogènes à valeurs dans $L_P$ sur $\calh(T(M))$} n'ayant que des composantes pures de type $(0,k)$ et certaines propriétés, supplémentaires.\\
Pour chaque $x\in M$ et chaque élément $\omega$ de l'espace cotangent complexifié en $x$ a $M$, on a une {\sl multiplication de Clifford} $\gamma_P(\omega)$ qui est un endomorphisme de la fibre de $S_P$ en $x$ et correspond à l'opération $\gamma$ de 9.1; $\gamma_P:T^\ast(M)\otimes \mathbb C\rightarrow \End(S_P)$ est un homomorphisme de fibrés vectoriels complexes sur $M$. $\gamma_P(\bar\omega)$ est l'adjoint de $\gamma_P(\omega)$ pour la structure hermitienne de $S_P$.\\

$\gamma_P$ se relève à $\calh(T(M))$ en un homorphisme 
\[
\tilde\gamma_P:\wedge^1_H\calh(T(M))\rightarrow \End(\pi^\ast(S_P))=\End((\oplus_k\wedge^{0,k}_H\calh(T(M))\otimes L_P))
\]
 et il résulte de la proposition 9.5 que la restriction de $\tilde\gamma_P$ à $\wedge^{0,1}_H\calh(T(M))$ {\sl n'est autre que} 1/$\sqrt{2}$ fois {\sl le produit extérieur} de $\wedge^{0,1}_H\calh(T(M))$ par 
 $$(\oplus_k\wedge^{0,k}_H\calh(T(M)))\otimes L_P$$
  i.e. on a :
\[
\sqrt{2}\tilde \gamma_P(\omega)\varphi=\omega\wedge\varphi,
\]
$\forall \omega\in \wedge^{0,1}_H\calh(T(M))$ et $\varphi\in \oplus_k\wedge^{0,k}_H\calh(T(M))\otimes L_P$ (au même point de $\calh(T(M))$

\subsection{Opérateur de Dirac}

$P$ et $M$ étant comme dans 10.6 la connexion riemannienne de $M$ se relève en une connexion unique sur $P$. Nous désignerons par $\nabla_P$ la différentielle covariante correspondante
des sections de $S_P$; $\nabla_P:\Gamma(S_P)\rightarrow(S_P\otimes T^\ast(M))$. La multiplication de Clifford $\gamma_P$ définit canoniquement une application de $\Gamma(S_P\otimes T^\ast(M))$ dans $\Gamma(S_P)$ que nous noterons encore $\gamma_P:\Gamma(S_P\otimes T^\ast(M))\rightarrow \Gamma(S_P)$ ($\gamma_P(\psi\otimes\omega)=\gamma_P(\omega)\psi$). Par composition, on a un opérateur différentiel de premier ordre, $D_P=\gamma_P\circ \nabla_P(=``\gamma_P(\nabla_P)")$, sur $\Gamma(S_P)$ qui est appelé {\sl opérateur de Dirac}.\\

L'opérateur de Dirac $D_P$ se relève à $\calh(T(M))$ en un opérateur $\tilde D_P$. L'analyse précédente (10.6) suggère que, dans les ``bons cas" $\tilde D_P$ est proportionnel à la partie formellement hermitienne d'un opérateur de cohomologie $\bar\partial$ sur les formes à valeurs dans $L_P$. A priori, les ``bons cas" sont ceux pour lesquels la structure presque complexe de $\calh(T(M))$ (resp. de $\calh_+(T(M))$ ou $\calh_-(T(M))$) est intégrable; en effet, $L_P$ est alors canoniquement un fibré holomorphe sur $\calh(T(M))$ (resp. $\calh_+(T(M))$ ou $\calh_-(T(M))$)et on a un opérateur $\bar\partial$ bien défini sur les formes à valeurs dans $L_P$.\\

Examinons d'abord le cas où $M$ est plate.

\subsection{Proposition}

{\sl Soit $M$ une variété riemannienne plate de dimension $2\ell$ et $P$ une structure spinorielle sur $M$. On a pour tout champ spinoriel $\psi\in \Gamma(S_P)$
\[
\sqrt{2}\pi^\ast(D_P\psi)=(\bar\partial + \bar\partial^+)\pi^\ast(\psi)
\]
où $\pi^\ast(\psi)$ est (comme dans 10.6) la forme différentielle inhomogène sur $\calh(T(M))$ à valeurs dans $L_P$ correspondant à $\psi$ et $\bar\partial^+$ désigne l'adjoint formel de l'opérateur de cohomologie $\bar\partial$ sur les formes à valeurs dans le fibré holomorphe $L_P$ sur $\calh(T(M))$.}\\

\noindent\underbar{Démonstration}. Les formes horizontales du type $\pi^\ast(\psi)$ sont stables par $\bar\partial$ et $\bar\partial^+$. On peut se restreindre à une section horizontale, $i$, de $\pi:\calh(T(M))\rightarrow M$. L'énoncé étant de nature locale sur $M$, on peut se ramener au cas $M=\mathbb R^{2\ell}$. Le fibré $L_P$ est alors trivial sur $i$ et tout revient à identifier isométriquement $\mathbb R^{2\ell}$ à $\mathbb C^\ell$ (via la section $i$) et à identifier $S$ (ou plutôt $\pi^\ast(S)\restriction i$) ) à l'algèbre extérieure construite sur $\overline{\mathbb C^{\ell\ast}}$. Les champs spinoriels sont alors les formes différentielles sur $\mathbb C^\ell$ n'ayant que des composantes pures de types $(0,k)$, ($k=0,1,\dots,\ell$); si $\omega$ est une forme de type $(0,1)$, $\sqrt{2}\gamma(\omega)$ est la multiplication extérieure par $\omega$ et $\gamma(\bar\omega)=\gamma(\omega)^+$. On a donc $\sqrt{2}\gamma(d)\psi=(\bar\partial+\bar\partial^+)\psi$ ce qui est le résultat cherché.~$\square$

\subsection{Corollaire}

{\sl Soit $M$ une variété riemannienne compacte conforme plate de dimension $2\ell$ et soit $P$ une structure spinorielle sur $M$. On a une injection canonique de l'espace des solutions de $D_P\psi=0$ $(\psi\in\Gamma(P_P))$ dans la cohomologie $\oplus^{k=\ell}_{k=0}H^{0,k}(\calh(T(M)),L_P)$.}\\

En effet, l'équation $D_P\psi=0$ est invariante conforme (on peut rendre $D_P$ invariant conforme en considérant les spineurs de poids $2\ell-1/2$), ``on peut donc se ramener au cas plat et utiliser 10.8". La compacité de $M$ implique que l'on a séparément $\bar\partial\pi^\ast(\psi)=0$ et $\bar\partial^+\pi^\ast(\psi)=0$ et le résultat est une conséquence de la théorie harmonique sur les variétés complexes compactes.

\subsection{Remarques}

a) Ce résultat suggère qu'il est possible de généraliser aux dimensions paires la version riemannienne de la transformation de Penrose des équations de champs linéaires invariantes conformes sur les variétés conformes plates. De ce point de vue 10.9 est incomplet car, en dimension 4, on a essentiellement une bijection.\\

b) Si $M$ est orientée, on peut se restreindre à $\calh_+(T(M))$ et les solutions de $D_P\psi=0$ pour les semi-spineurs dans $S_{P^+}$ (resp. $S_{P^-}$) correspondent à la somme sur les valeurs paires (resp. impaires) de $k$ des $H^{0,k}(\calh_+(T(M)),L_P)$.

\newpage
\enlargethispage{2cm}
\begin{center}
{\large Références pour les chapitres 5, 6, 7, 8, 9 et 10}
\end{center}
\begin{itemize}
\item
Atiyah M.F., Hitchin N.J., Singer I.M. : Self-duality in four-dimesioanl Riemannian geometry. {\sl Proc. R. Soc. London\/ \bf A362} (1978), 425-461.\\

\item
Budinich P., Trautman A. : Fock space description of simple spinors. {\sl J. Math. Phys.\/ \bf 30} (1989), 2125-2131.\\

\item
Cartan E. : Leçons sur la théorie des spineurs, I, I. Hermann et Cie, Paris 1938.\\

\item
Chevalley C. : The algebraic theory of spinor, Columbia University Press, 1954.\\

\item
Dubois-Violette M. : Structures complexes au-dessus des variétés, appli\-cations $^\ast$, dans ``Mathématique et Physique", Séminaire à l'E.N.S.,\linebreak[4] Boutet de Monvel L., Douady A. et Verdier J.L., eds. {\sl Progress in Mathematics\/ \bf 37}, Birkhaüser 1983, pp. 1-42.\\

\item
Dubois-Violette M. :  Complex structures and the Elie Cartan approach to the theory of spinors in ``Spinors, Twistors, Clifford Algebras and Quantum Deformations", Oziewicz, Z. ed., Kluwer Academic Press 1993, pp. 17-23. \\

\item
Fegan H.D. : Conformally invariant first order differential operators. {\sl Quart. J. Math. Oxford (2), \/ \bf 27} (1976), 371-378.\\

\item
Hitchin N.J. : Linear field equations on self-dual spaces. {\sl Proc. R. Soc. Lond.\/ \bf A370} (1980), 173-191.\\

\item
Madore J., Richard J.L., Stora R. : An introduction to the twistor programme. {\sl Phys. Rep.\/ \bf 49} (1979), 113-130.\\

\item
Les articles de Penrose, R. sur les ``twisteurs".
\end{itemize}

\vspace{1cm}

\noindent $^\ast$ Les sections 5, 6, 7, 8, 9 et 10 sont essentiellement une reproduction de cette référence qui est le texte d'un exposé au Séminaire de Mathématique de l'E.N.S. en 1981.

\newpage

\section{Variétés quaternioniques}
\subsection{Définitions : $\GL(k,\mathbb H)$ et $\GL(k,\mathbb H)\bullet Sp(1)$}

Soit $e_0, e_1, e_2, e_3$ la base canonique de $\mathbb R^4$ (euclidien orienté); on définit sur $\mathbb R^4$ une structure d'algèbre associative en posant $e_0e_k=e_ke_0=e_k, e^2_0=e_0$ $e_ke_\ell=-\delta_{k\ell}+\sum^3_{m=1}\epsilon_{k\ell m}e_m$ ($k,\ell=1,2,3$). $e_0$ est l'unité de cette algèbre. Si $q=x^0e_0+x^ke_k$, on définit son conjugué $\bar q$ par $\bar q=x^0e_0-x^ke_k$; on a $q\bar q=\bar q q=\lVert q\rVert^2 e_0=((x^0)^2+\sum(x^k)^2)e_0$. Tout élément $q\not =0$ a un inverse, $\frac{\bar q}{\lVert q\rVert^2}$, donc l'algèbre ainsi définie est un corps appelé {\sl corps des quaternions}; nous le désignerons par $\mathbb H$. On considérera $\mathbb H^n$ comme un espace vectoriel quaternionique {\sl à droite}; les applications linéaires quaternioniques correspondent alors à la multiplicatio {\sl à gauche} par les matrices à coefficients quaternions. Le {\sl groupe linéaire quaternionique} $\GL(k,\mathbb H)$ sera considéré comme un sous-groupe de $\GL(4k,\mathbb R)$ en identifiant $\mathbb H^k$ avec $\mathbb R^{4k}$ et nous poserons 
$Sp(k)=\GL (k,\mathbb H)\cap O(4k)$ (= groupe des isométries quaternioniques). On a canoniquement
$Sp(1)=\{q\in \mathbb H\vert \lVert q\rVert^2=1\} \simeq SU(2)$ et $\GL(1,\mathbb H)=\mathbb H\backslash \{0\}=\mathbb R^+\bullet Sp(1)$. On remarquera que les éléments de $\GL(k,\mathbb H)$ {\sl préservent l'orientation} de $\mathbb R^{4k}$. On définit un homomorphisme $\GL(k,\mathbb H)\times Sp(1)\rightarrow \GL(4k,\mathbb R)$ noté $(L,u)\mapsto L\bullet u$ en posant $L\bullet u(V)=L(V\bar u)=L(V)\bar u$, $\forall V\in \mathbb H^k$, $L\in \GL(k,\mathbb H)$ et $u\in Sp(1)$ ; le noyau de cet homomorphisme est $\{(1,1),(-1,-1)\}\simeq \mathbb Z_2$, {\sl son image sera notée} $\GL(k,\mathbb H)\bullet Sp(1)$ $(\subset \GL(4k,\mathbb R))$.

\subsection{Définition de $Z_+(\mathbb R^{4n})$}

Dans le cas $k=1$, on a : $Sp(1)\bullet Sp(1)=SO(4)$, $\GL(1,\mathbb H)\bullet Sp(1)=CO^+(4)$, $\GL(1,\mathbb H)$ est l'ensemble des transformations conformes de $\mathbb R^4$ commutant avec les matrices antisymmétriques anti-auto-duales (donc avec les éléments de $\calh_+(\mathbb R^4)$ (voir 7)).\\
De manière générale, la multiplication à droite $J_u$, par $\sum^3_1 u^ke_k$ dans $\mathbb H^n$ avec $\vec u^2=1$ est une {\sl structure complexe isométrique d'orientation positive sur} $\mathbb R^{4n}$ ($J_u\in \calh_+(\mathbb R^{4n})$; {\sl nous désignerons l'ensemble de ces structures complexes sur} $\mathbb R^{4n}$ par $Z_+(\mathbb R^{4n})$, c'est une sous-variété algébrique complexe de $\calh_+(\mathbb R^{4n})$ qui est isomorphe à $P_1(\mathbb C)=\calh_+(\mathbb R^4)$. $\GL(n,\mathbb H)$ {\sl est le sous-groupe de $\GL(4n,\mathbb R)$ laissant $Z_+(\mathbb R^{4n})$ invariant point par point} tandis que {\sl $\GL(n,\mathbb H)\bullet Sp(1)$ est le sous-groupe de $\GL(4n,\mathbb R)$ laissant globalement $Z_+(\mathbb R^{4n})$ invariant}; il est évident que $\GL(n,\mathbb H)\bullet Sp(1)$, (et même $1\bullet Sp(1))$, opère transitivement sur $Z_+(\mathbb R^{4n})$.

\subsection{Définition des espaces projectifs quaternioniques}

L'ensemble des (directions) sous-espaces quaternioniques de dimension 1  de $\mathbb H^{n+1}$ est appelé {\sl espace projectif quaternionique de dimension quaternionique} $n$; nous le noterons $P_n(\mathbb H)$. $P_n(\mathbb H)$ est le quotient de $\mathbb H^{n+1}\backslash \{0\}$ par la relation d'équivalence $v_1\sim v_2$ si $v_1=v_2 q$, ($q\in \mathbb H$, $q\not=0$) nous munirons $P_n(\mathbb H)$ de la topologie quotient. Soit $U_k$ ($k=0,\dots,n)$ l'ouvert de $P_n(\mathbb H)$ image de $\{(q_0,q_1,\dots q_n)\in \mathbb H^{n+1}\vert q_k\not= 0\}$; posons pour $\xi\in U_k$ 
\[
\varphi_k(\xi)=(q_0q^{-1}_k,\dots,q_{k-1} q^{-1}_k, q_{k+1}q^{-1}_k,\dots,q_nq^{-1}_k)\in \mathbb H^n=\mathbb R^{4n},
\]
où $(q_0,\dots,q_n)$ est au-dessus de $\xi$. Si $\xi\in U_k\cap U_\ell$, on a évidemment\linebreak[4] $d(\varphi_\ell\circ \varphi_k^{-1})_x\in \GL(n,\mathbb H)\bullet Sp(1)$,
$\forall x\in \varphi_k(U_k\cap U_2)$ et $\varphi_\ell\circ \varphi_k^{-1}$ est un difféomorphisme de $\varphi_k(U_k\cap U_\ell)$ sur $\varphi_\ell(U_k\cap U_\ell)$. On a ainsi un atlas compatible avec $\Gamma_{\GL(n,\mathbb H)\bullet Sp(1)}(\mathbb R^{4n})$ sur $P_n(\mathbb H)$ et $P_n(\mathbb H)$ devient ainsi une variété différentiable de type $\GL(n,\mathbb H)\bullet Sp(1)$, (voir dans 1.3), ou ce qui revient au même, une variété différentiable munie d'une $GL(n,\mathbb H)\bullet Sp(1)$-structure intégrable. On notera que cette structure n'est pas pas une $GL(n,\mathbb H)$-structure et, en fait, {\sl il est impossible de munir $P_n(\mathbb H)$ d'une $GL(n,\mathbb H)$-structure}, (voir l'exemple ci-dessous).

\subsection{Exemple de $P_1(\mathbb H)=S^4$}

L'application de $\mathbb H^2\backslash \{0\}$ sur $S^4$ définie par 
\[
(q_0,q_1)\mapsto \xi=\left(\frac{2q_0\bar q_1}{\lVert q_0\rVert^2+\lVert q_1\rVert^2}, \frac{\lVert q_0 \rVert^2-\lVert q_1\rVert^2}{\lVert q_0\rVert^2+\lVert q_1\rVert^2}\right)
\]
 passe au quotient et définit un difféomorphisme de $P_1(\mathbb H)$ sur $S^4$. On a donc $P_1(\mathbb H)=S^4$.\\
Il est connu que, pour $\ell\not= 1$ et $\not=3$, $S^{2\ell}$ n'admet aucune structure presque complexe; comme, d'après 11.2, $GL(1,\mathbb H)\subset GL(\mathbb R^4_J)\simeq GL(2,\mathbb C)$ pour tout $J\in Z_+(\mathbb R^4)=\calh_+(\mathbb R^4)$, il est clair que $S^4=P_1(\mathbb H)$ n'admet pas de\linebreak[4] $GL(1,\mathbb H)$-structure.\\
Nous avons vu dans 11.2 que $GL(1,\mathbb H)\bullet Sp(1)$ est le groupe $CO^+(4)$ des transformations linéaires conformes préservant l'orientation de $\mathbb R^4$; la\linebreak[4] $GL(1,\mathbb H)\bullet Sp(1)$-structure intégrable de $P_1(\mathbb H)=S^4$ est la structure conforme plate orientée usuelle de $S^4$\\
Ces préliminaires sont la motivation de la définition suivante.

\subsection{Variétés presque quaternioniques}

Une variété différentiable de dimension $4k$ munie d'une $GL(k,\mathbb H)\bullet Sp(1)$-structure sera appelée une {\sl variété presque quaternionique}. Il résulte de 11.2 qu'une variété presque quaternionique est orientée et qu'une variété presque quaternionique de dimension 4 n'est autre qu'une {\sl variété orientée de dimension 4 munie d'une structure conforme}.\\

Soit $\fracg\subset \fracgl(4k,\mathbb R)$ l'algèbre de Lie de $\GL(k,\mathbb H)\bullet Sp(1)$; on peut montrer que $H^{r,2}(\fracg)=0$ si $r\geq 2$, de sorte que l'on n'a que deux obstructions $c_0(P)$ et $c_1(P)$ à l'intégrabilité d'une $\GL(k,\mathbb H)\bullet Sp(1)$-structure $P$ sur $V$. L'intégrabilité de $P$, i.e. $c_0(P)=0=c_1(P)$, entraîne que $V$ est localement isomorphe à $P_k(\mathbb H)$. Cela conduit à ``relaxer" la condition $c_1(P)=0$ et à adapter la définition suivante.

\subsection{Variétés quaternioniques}

Une {\sl variété quaternionique} est une variété (de dimension $4k$) munie d'une $\GL(k,\mathbb H)\bullet Sp(1)$-structure admettant une connexion de torsion nulle $(c_0(P)=0)$.\\
Une variété de dimension $4k$ munie d'une $\GL(k,\mathbb H)$-structure sera appelée une {\sl variété presque hyperquaternionique}; elle sera appelée {\sl hyperquaternionique} dans le cas où cette $\GL(k,\mathbb H)$-structure admet une connexion de torsion nulle.

\subsection{Le fibré $Z_+(T(V))$}

Soit $V$ une variété presque quaternionique de dimension $4k$; on a un fibré de fibre type $Z_+(\mathbb R^{4k})$ correspondant à l'action de $\GL(k,\mathbb H)\bullet Sp(1)$ sur $Z_+(\mathbb R^{4k})$ associé à la $\GL(k,\mathbb H)\bullet Sp(1)$-structure sur $V$. {\sl Nous noterons ce fibré} $Z_+(T(V))$; $Z_+(T(V))$ est parfois appelé {\sl espace des twisteurs de $V$}. La structure complexe de $Z_+(\mathbb R^{4k})\simeq P_1(\mathbb C)$ est invariante par\linebreak[4] $\GL(k,\mathbb H)\bullet Sp(1)$ et $Z_+(\mathbb R^{4k})$ est une sous-variété complexe de $\cali_+(\mathbb R^{4k})$, (voir Chapitre 5). Il en résulte que $Z_+(T(V))$ est un sous-fibré de $\cali_+(T(V))$ dont les fibres (isomorphes à $P_1(\mathbb C)=\calh_+(\mathbb R^4)$) sont des sous-variétés complexes des fibres de $\cali_+(T(V))$.\\
En particulier, pour toute connexion sur la $GL(k,\mathbb H)\bullet Sp(1)$-structure de $V$, $Z_+(T(V))$ est une sous-variété presque complexe de $\cali_+(T(V))$ pour les stuctures presque complexes associées (voir Chapitre 6).

\subsection{Théorème}

{\sl Soit $V$ une variété quaternionique. Les structures presque complexes sur $Z_+(T(V))$ associées à deux connexions sans torsions sur la $GL(k,\mathbb H)\bullet Sp(1)$-structure de $V$ coincident $(\dim V=4k)$. Si $k$ est supérieur ou égal à 2, (i.e. $\dim V\geq 8$), la structure presque complexe ainsi définie sur $Z_+(T(V))$ est intégrable}.\\

On remarquera que pour $k=1$ ($\dim V=4$) ce résultat est, en vertu de 11.2, un cas particulier du théorème 8.2.\\

Soit $P$ la $\GL(k,\mathbb H)\bullet Sp(1)$-structure de $V$ et soit $\fracg^{(1)}$ le prolongement d'ordre 1 de l'algèbre de Lie $\fracg$ de $\GL(k,\mathbb H)\bullet Sp(1)\subset \GL(4k,\mathbb R)$ on a un fibré associé \underbar{$\fracg^{(1)}$} à  $P$ sur $V$ (voir 1) et la différence $K$ de deux connexions sans torsions sur $P$ est une section de \underbar{$\fracg^{(1)}$}. La première partie du théorème résulte de 6.4 et du fait que si $t\in \fracg^{(1)}$ on a :
$t(JX,JY)-Jt(JX,Y)-Jt(X,JY)-t(X,Y)=0$, $\forall X,Y\in \mathbb R^{4k}$ et $\forall J\in Z_+(\mathbb R^{4k})$.\\

La deuxième partie s'obtient en utilisant 6.5 et la première identité de Bianchi (=torsion nulle). On trouve en prenant des traces de combinaisons appropriées que $(k-1)N=0$ où $N$ est la torsion de la structure presque complexe de $Z_+(T(V))$ (associée à une quelconque des connexions sans torsion sur $P$).

\subsection{Remarque sur la dimension 4}

Dans 11.8, on retrouve la situation très particulière de la dimension 4. En dimension 4, une $\GL(1,\mathbb H)\bullet Sp(1)$-structure est une structure conforme orientée donc elle admet toujours une connexion de torsion nulle ($c_0=0$); autrement dit, une structure presque quaternionique est quaternionique et la structure presque complexe de $Z_+=\calh_+$ n'est intégrable que si la structure conforme orientée est auto-duale (voir 8). C'est pourquoi, certains auteurs définissent les variétés quaternioniques de dimension 4 comme étant des variétés de dimension 4 munies de structures conformes orientées auto-duales. Pour d'autres auteurs, les variétés quaternioniques sont de dimensions $4k\geq 8$.

\subsection{Remarque sur le cas hyperquaternionique}

$\GL(k,\mathbb H)=\GL(k,\mathbb H)\bullet 1\subset \GL(4k,\mathbb R)$ opère trivialement sur $Z_+(\mathbb R^{4k})$ comme on l'a vu dans 11.2; il en résulte que, si $V$ est une variété presque hyperquaternionique de dimension $4k$, le fibré $Z_+(T(V))$ {\sl est un fibré trivial} et les sections correspondantes, qui sont indexées par $Z_+(\mathbb R^{4k})\simeq P_1(\mathbb C)$, sont des structures presque complexes sur $V$. Toute connexion sur la $\GL(k,\mathbb H)$-structure de $V$ est presque complexe pour chacune de ces structures presque complexes sur $V$. Le théorème 3.4 implique donc que {\sl si $V$ est hyperquaternionique toutes ces structures presques complexes indexées par $P_1(\mathbb C)$ sont intégrables}.

\subsection{Définitions supplémentaires}

La multiplication à droite dans $\mathbb H^k (\simeq \mathbb R^{4k})$ par les éléments de $\mathbb H$ détermine un isomorphisme de $\mathbb H$ sur une sous-algèbre $\mathbb H_+(\mathbb R^{4k})$ de l'algèbre $\End(\mathbb R^{4k})$ des endomorphismes de $\mathbb R^{4k}$ (contenant $\bbbone \in \End(\mathbb R^{4k})$). Si $V$ est une variété presque quaternionique de dimension $4k$ on a un fibré $\mathbb H_+(T(V))$ associé à sa $\GL(k,\mathbb H)\bullet Sp(1)$-structure correspondant à l'action (adjointe) de\linebreak[4] $\GL(k,\mathbb H)\bullet Sp(1)$ sur $\mathbb H_+(\mathbb R^{4k})$, (l'action de $\GL(k,\mathbb H)$ est triviale, on a en fait une action adjointe de $Sp(1)\simeq SU(2)$ i.e. une action de $SO(3)$). $\mathbb H_+(T(V))$ est un fibré en algèbres de quaternions; ce fibré est trivial si et seulement si $V$ est presque hyperquaternionique. De manière générale, nous appellerons {\sl structure presque quaternionique} sur une variété $V$ un sous-fibré en algèbres $K$ du fibré des endomorphismes de $T(V)$ contenant l'identité et dont les fibres sont isomorphes à $\mathbb H$; si $K$ est trivial,  nous parlerons de {\sl structure presque hyperquaternionique}. Si $V$ est munie d'une structure presque quaternionique $K$, on peut munir ses espaces tangents $T_x(V)$ de structures d'espaces vectoriels quaternioniques (à droite) mais $T(V)$ n'est un fibré vectoriel quaternionique (i.e. fonctions de transitions $\mathbb H$-linéaires) que si $K$ est trivial (i.e presque hyperquaternionique). Toute structure presque quaternionique $K$ sur $V$ détermine donc une $\GL(k,\mathbb H)\bullet Sp(1)$-structure sur $V$ pour un certain $k$ ($\Rightarrow \dim V=4k$) unique telle que $K=\mathbb H_+(T(V))$; $Z_+(T(V))=\{q\in \mathbb H_+(T(V))\vert q^2=-\bbbone\}$, et $q\in \mathbb H_+(T(V))$ s'écrit de manière unique $q=r\bbbone +\rho u$ avec $u\in Z_+(T(V))$, $r\in \mathbb R$ et $\rho\in \mathbb R_+$ (i.e.$\rho\geq 0$). En particulier, au voisinage de chaque point de $V$ on peut trouver des sections locales de $Z_+(T(V))$ $e_1,e_2,e_3$ satisfaisant $e_i e_j=-\delta_{ij}+\sum_k \epsilon_{ijk} e_k$ dans $\mathbb H_+(T(V))$ et $\mathbb H_+(T(V))=\mathbb R\bbbone \oplus \mathbb R e_1\oplus \mathbb R e_2 \oplus \mathbb R e_3$ dans le domaine de ces sections.\\

Si $V$ est une variété presque quaternionique de dimension $4k$, les connex\-ions linéaires $\nabla$ telles que $\forall X\in \Gamma(T(V))$, $\nabla_X\Gamma(\mathbb H_+(T(V))\subset \Gamma(\mathbb H_+(T(V))$) seront appelées {\sl presque quaternioniques}; ce sont les connexions sur la\linebreak[4] $\GL(4k,\mathbb H)\bullet Sp(1)$ structure de $V$ (ce sont aussi les connexions laissant\linebreak[4] $Z_+(T(V))$ globalement invariant).

\subsection{Variétés hyperkähleriennes}

Une variété presque quaternionique $V$ munie d'une métrique riemannienne $g$ telle que $g(uX,uY)=g(X,Y)$, $\forall X,Y\in T_x(V)$, $\forall u\in Z_+(T_x(V))$ et $\forall x\in V$ (i.e. telle que $Z_+(T(V))\subset \calh_+(T(V))$) est appelée {\sl variété presque hermitienne quaternionique}. Nous dirons que $V$ est {\sl kählerienne quaternionique} si la connexion riemannienne de $V$ est presque quaternionique.\\

Une variété kählerienne quaternionique est une variété quaternionique puisque la connexion riemannienne est sans torsion.\\
Une variété kählerienne quaternionique $V$ sera appelée {\sl hyperkählerienne} si sa connexion riemannienne induit une connexion de courbure nulle sur $\mathbb H_+(T(V))$ (ou, ce qui revient au même sur $Z_+(T(V))$).\\

Si $V$ est hyperkählerienne, les sections horizontales de $Z_+(T(V))$ sont des structures presque complexes isométriques et covariantes constantes; il en résulte (3.16) que $V$ {\sl est kählerienne pour chacune de ces structures (presque) complexes} indexées par $Z_+(\mathbb R^{4k})=P_1(\mathbb C)$.

\subsection {Cas de la dimension 4}
{\sl Une variété riemannienne $V$ de dimension 4 orienté est hyperkählerienne si et seulement si sa courbure riemannnienne est auto-duale}. En effet\linebreak[4] $Z_+(T(V))=\calh_+(T(V))$ est constitué par les structures complexes (sur les $T_x(V)$) correspondant aux 2-formes de norme 1 anti-autoduales; le ``commutant de
$Z_+$ dans les endomorphismes anti-symétriques" est constitué par ceux qui sont auto-duals $\Rightarrow$ l'énoncé en utilisant les symétries du tenseur de courbure. (Pour le cas anti-autodual) il suffit de changer l'orientation pour se ramener au cas auto-dual).

\subsection{Théorème}
{\sl La courbure de Ricci d'une variété hyperkählerienne est nulle. Toute variété kählerienne quaternionique de dimension $4k\geq 8$ est une variété d'Einstein}, (i.e. la courbure de Ricci est proportionnelle à la métrique).\\

Revenons à la dimension 4 pour laquelle on a des propriétés intéressantes particulières.

\subsection{Proposition}
{\sl Soit $V$ une variété riemannienne orientée de dimension 4; les propriétés} (a), (b) {\sl et} (c) {\sl suivantes sont équivalentes}\\

(a) $V$ {\sl est kählerienne et sa courbure de Ricci est nulle.}\\

(b) {\sl La courbure riemannienne de $V$ est auto-duale (comme 2-forme)}.\\

(c) $V$ {\sl est hyperkählerienne}.\\

\noindent\underbar{Démonstration}. Montrons $(a)\Rightarrow (b)\Rightarrow (c) \Rightarrow (a)$.\\

(a) $\Rightarrow$ (b). Soit $J\in \Gamma(\calh_+(T(V))$) la structure presque complexe de la variété kählerienne $V$. Comme $\nabla J=0$ pour la connexion riemannienne, on a $[\Omega, J]=0$, où $\Omega$ est la 2-forme de courbure riemannienne. Il en résulte $\Omega=\Omega^++J\otimes \lambda$ où $\Omega^+$ est (à valeurs dans les matrices anti-symmétriques auto-duales) dans $\Gamma(\calh_-(T(V))\otimes \wedge^2T^\ast(V))$ et où $\lambda$ est une 2-forme scalaire. L'annulation du tenseur de Ricci implique $\lambda=0$. On a donc $\Omega=\Omega^+$ et l'auto-dualité de $\Omega$ résulte de la symétrie du tenseur de courbure d'une variété riemannienne $(\Omega=\Omega^+ \Leftrightarrow \Omega =\ast \Omega)$\\

(b) $\Rightarrow$ (c).  $\Omega=\ast \Omega\Leftrightarrow \Omega=\Omega^+\in \Gamma(\calh_-(T(V))\otimes \wedge^2 T^\ast (V))\Leftrightarrow [\Omega,J]=0$ $\forall J\in \Gamma(\calh_+(T(V))$); autrement dit, la connexion riemannienne induit une connexion de courbure nulle sur $\calh_+(T(V))$. $V$ est donc hyperkählerienne avec $Z_+(T(V))=\calh_+(T(V))$.\\

(c) $\Rightarrow$ (a). En choisissant une section horizontale de $Z_+(T(V))=\calh_+(T(V))$, on a une structure presque complexe pour laquelle $V$ est kählerienne et l'annulation du tenseur de Ricci résulte de $\Omega=\ast\Omega$.~$\square$

\subsection{Exemples : Solutions de Hawking}

Soit $V$ une variété riemannienne orientée de dimension 4 qui est hyperkählerienne. La proposition 11.15 (b) implique qu'au voisinage de tout point de $V$  existe un champ de repères ($e_0,e_1,e_2,e_3$) orthonormés d'orientation positive, relativement auxquels la partie anti auto-duale de la forme de la connexion riemannienne s'annule, i.e dans lequel la connexion $\omega^a_b$ est à valeurs matrice antisymétrique auto-duale. Soit $(\theta^0,\theta^1,\theta^2,\theta^3)$ le corepère dual et considérons les 2-formes anti auto-duales
\[
j_k=\theta^0\wedge \theta^k- \theta^\ell \wedge \theta^m
\]
pour $k \in \{1,2,3\}$, ($k,\ell, m$) étant la permutation cyclique de (1,2,3) correspondante. Ces 2-formes constituent un champ de bases orthonormées des 2-formes anti auto-duales.\\
Les équations de structure et l'auto-dualité des $\omega^a_b$ implique
\[
dj_k=0
\]
pour tout $k\in \{1,2,3\}$. Inversement l'existence d'un tel champ local de bases orthonormées de 2-formes anti auto-duales constituées par des formes fermées sur une variété riemannienne orientée $V$ implique que $V$ est localement hyperkählerienne.\\
Cherchons à construire de tels $j_k$ en posant
\[
\theta^0=\frac{1}{\sqrt{A_0}}(dx^0+A_n dx^n)
\]
\[
\theta^k=\sqrt{A_0} dx^k
\]
où $A_0$ et les $A_k$ ($k\in \{1,2,3\}$) ne dépendent que des $x^k$ avec $k\in \{1,2,3\}$ et $A_0>0$. On a 
\[
j_k=(dx^0+A_n dx^n)\wedge dx^k-A_0 dx^\ell \wedge dx^m
\]
où ($k,\ell,m$) est une permutation cyclique de (1,2,3). Les équations
\[
dj_k=0
\]
s'écrivent alors
\[
\partial_k A_0=\partial_\ell A_m-\partial_m A_\ell
\]
pour toute permutation cyclique $(k,\ell,m)$ de (1,2,3). Elles impliquent que $A_0$ est harmonique
\[
\Delta A_0=0
\]
Inversement si $A_0$ est harmonique, on peut résoudre les équations ci-dessus en $A_k$ $(\grad (A_0)=\vec{\rot}(\vec A))$ et la solution est unique au gradiant d'une fonction des $x^k$ près (qui peut être absorbée dans une redéfinition de $x^0$). La métrique riemannienne correspondante
\[
ds^2=\frac{(dx^0+A_ndx^n)^2}{A_0}+A_0 \sum^3_{k=1} (dx^k)^2
\]
est par construction hyperkählerienne.\\

Remarquons qu'en posant
\[
A=A_0 dx^0+A_n dx^n
\]
et en tenant compte de $\frac{\partial A_\mu}{\partial x^0}=0$, les équations précédentes pour les $A_\mu$ sont équivalentes à l'anti auto-dualité de la 2-forme $F=dA$ sur $\mathbb R^4$, (i.e. $F+\ast F=0$).

\subsection{Analogue non linéaire de l'équation de Laplace}

Supposons que $V$ satisfasse les conditions équivalentes de la proposition 11.5. Alors, d'après 3.29 (3), la condition (a) de 11.5 implique que localement, on peut trouver des coordonnées complexes $z_1$ et $z_2$ de $V$ et un potentiel de Kähler $F$ tels que l'on ait
\[
\det (\partial\bar\partial F)=1
\]
Cette équation est un analogue non linéaire de l'équation de Laplace. En effet, la solution plate triviale est
\[
F_0=\vert z_1\vert^2 + \vert z_2\vert^2
\]
et pour une petite perturbation
\[
F=F_0+\varepsilon \phi
\]
l'équation précédente donne au premier ordre en $\varepsilon$
\[
\Tr (\partial\bar\partial \phi)=\Delta \phi=0
\]
ce qui est l'équation de Laplace pour le champ scalaire réel $\phi$.\\

 L'importance de cette remarque est que l'auto-dualité de la courbure correspond à la version ``euclidienne" d'un des 2 degrés d'hélicité du champ de gravitation et que en faisant ce choix de coordonnées (de jauge) nous avons, ainsi isolé un degré local des 2 degrés de liberté de la version riemannienne du champ de gravitation ainsi que la version non linéaire de l'équation de Laplace satisfaite par l'approximation linéaire. Signalons par ailleurs que cette équation non linéaire est du type Monge-Ampère complexe qui a été très étudiée et qui est une ``bonne version" non linéaire de l'équation de Laplace.\\

\newpage

\begin{center}
{\large Références pour le chapitre 11}
\end{center}
\vspace{1cm}
\begin{itemize}
\item
Alekseevskii, D.V. : Compact quaternion spaces. {\sl Funct. Anal. Appl.\/ \bf 2} (1968) 106-114.\\

\item
Bonan, E. : Sur les $G$-structures de type quaternionien. {\sl Cahiers Topologie Géom. Diff.\/ \bf 9} (1967) 389-463.\\

\item
Gibbons, G.W., Hawking, S.W. : Gravitational multi-instantons. {\sl Phys. Lett.\/ \bf 78B} (1978), 430-432.\\

\item
Hawking, S.W. : Gravitational instantons. {\sl Phys. Lett.\/ \bf 60A} (1977), 81-83.\\

\item
Hitchin, N. : Polygons and gravitons. {\sl Math. Proc. Camb. Phil. Soc.\/ \bf 85} (1979), 465-476.\\

\item
Ishihara, S. : Quaternionic Kählerian manifolds. {\sl J. Diff. Geom.\/ \bf 4} (1974) 483-500.\\

\item
Salamon, S. : Quaternionic manifold. Thesis, Oxford (1980).\\

\item
Salamon, S. : Quaternionic Kähler manifolds. {\sl Invent. math.\/ \bf 67} (1982), 143-171.\\

\end{itemize}

\newpage

\section{Les équations de Yang et Mills}

Dans ce chapitre nous rappelons différentes formulations des équations de Yang et Mills et nous décrivons quelques exemples classiques de solutions.
Pour se rapprocher des notations utilisées en physique, les formes de connections seront notées $A$ et leur courbure $F$.

\subsection{Equations de Yang et Mills}

Soit $M$ une variété pseudo-riemannienne orientée de métrique $g$ et soit $G$ un groupe de Lie d'algèbre de Lie $\fracg$. {\sl Un champ de jauge sur $M$ de groupe de structure} $G$ correspond à la donnée d'un $G$-fibré principal $P$ sur $M$ muni d'une connexion $A$ à un ismorphisme près de $G$-fibrés principaux sur $M$ interchangeant les connexions.\\

Soit $P$ un $G$-fibré principal sur $M$ et soit $A$ une forme de connexion sur $P$; des équations pour $A$ correspondent à des équations pour le champ de jauge sous-jacent si elles sont {\sl invariantes de jauge}, i.e. invariantes par les isomorphismes de $G$-fibrés principaux sur $M$ interchangeant les connexions. $A$ satifait aux {\sl équations de Yang et Mills} (sans source) si l'on a :
\begin{equation}
\nabla \ast F =0
\end{equation}
où $F=dA+1/2[A,A]$ est la courbure de $A$, $\ast$ désigne l'opération de Hodge sur les formes sur $M$ relevée aux formes horizontales sur $P$ et $\nabla$ est la différentielle extérieure covariante associée à $A$. Ces équations sont invariantes de jauge.

\subsection{Remarque}

(12.1) est invariant par changement d'orientation de $M$ (i.e. par $\ast \mapsto -\ast$) et est local, il en résulte que, à ce stade, il n'est pas nécessaire de supposer $M$ orientable.

\subsection{Action de Yang et Mills}
Reprenons les notations de 12.1 et supposons que $\fracg$ soit muni d'une forme bilinéaire symétrique non dégénérée, $k$, invariante par l'action adjointe de $G$. A l'aide de $k$ et de $g$, on peut définir, en chaque point de $P$, une forme bilinéaire symétrique non dégénérée $(\bullet, \bullet)$ sur les $p$-formes horizontales à valeurs dans $\fracg$ en ce point. En coordonnées locales $\{x^\mu\}$ sur $M$ et dans une base $\{E_A\}$ pour $\fracg$ on a :
\[
(H,H')=1/p!\sum_{A,A',\mu_\ell,\mu'_\ell} k_{AA'} g^{\mu_1\mu'_1}\dots g^{\mu_p\mu'_p}H^A_{\mu_1\dots\mu_p}H'^{A'}_{\mu'_1\dots \mu'_p}
\]

$F$ étant la courbure d'une forme de connexion $A$ sur $P$, la fonction $(F,F)$ est constante sur les fibres de $P$; c'est donc, canoniquement, une fonction sur $M$. On définit alors {\sl l'action de Yang et Mills} correspondante $S_{\text{Y.M.}}(A)$ par
\begin{equation}
S_{\text{Y.M.}}(A)=-1/2 \int(F,F)\ \text{vol}_M
\end{equation}
$\text{vol}_M$ étant la forme élément de volume de $M$ et l'intégrale étant à prendre au sens indéfini, i.e. on a une action $S^\calo_{\text{YM}}(A)$ pour chaque région $\calo$ de $M$ telle que l'on puisse intégrer.\\

Les équations de Lagrange correspondant à $S_{\text{Y.M.}}$ (le lagrangien est la forme de degré maximum $-1/2(F,F)\ \text{vol}_M$ sur $M$) sont précisément les équations de Yang et Mills qui sont, par conséquent des équations variationnelles.

\subsection {Interprétations de type Kaluza-Klein}
On reprend les notations et les hypothèses de 12.3. On a sur $G$ une structure invariante de variété pseudo-riemannienne correspondant à $k$ sur $\fracg$. Nous désignerons par $\Lambda$ la courbure scalaire de cette variété pseudo-riemannienne; c'est une constante. Si $\xi$ est un point de $P$, on munit la fibre $F_\xi$ passant par $\xi$ de la métrique pseudo-riemannienne rendant isométrique l'application $g\mapsto \xi g$ de $G$ sur $F_\xi$; cette structure ne dépend pas du point $\xi$ choisi dans la fibre.\\

Si $A$ est une forme de connexion sur $P$, on peut relever la métrique de $M$ sur les sous-espaces horizontaux à l'aide de la projection $p:P\mapsto M$. On définit alors une métrique pseudo-riemannienne $\gamma$ sur $P$ muni de $A$ en imposant que $\gamma$ coîncide sur les sous-espaces verticaux et horizontaux avec les structures précédentes et que les sous-espaces verticaux soient orthogonaux aux sous-espaces horizontaux. $\gamma$ est invariante par l'action de $G$ sur $P$ et est donnée par
\begin{equation}
\gamma(X,Y)=k(A(X),A(Y))+g(p_\ast(X),p_\ast(Y))
\end{equation}

La courbure scalaire $\tilde R$ de la variété pseudo-riemannienne $(P,\gamma)$ est donnée par :
\begin{equation}
\tilde R = -1/2(F,F)+R+\Lambda
\end{equation}
où $R$ est la courbure scalaire de $M$ et $F$ est, comme précédemment, la courbure de la connexion $A$ sur $P$. $\tilde R$ est canoniquement une fonction sur $M$ et 
\begin{equation}
S=\int \tilde R\ \  \text{vol}_M
\end{equation}
considérée comme fonction de $A$ est équivalente à l'action de Yang et Mills (2); si on considère $S$ comme fonction de $A$ et de la métrique $g$ de $M$, les équations de Lagrange correspondantes, (``$\delta S(A,g)=0$"), pour $A$ et $g$ sont les équations couplées de Yang et Mills et d'Einstein avec constante cosmologique $\Lambda$ (voir plus loin).\\

Ces remarques sont le point de départ de certains types d'unification des théories de jauge avec la gravitation généralisant la théorie de Kaluza-Klein (i.e. le cas où $G=U(1)$).

\subsection{Laplacien de la connexion $A$}

$P,A,\gamma$, etc. étant comme au-dessus, on munit $P$ d'une orientation en choisissant une orientation de $\fracg$. On peut introduire l'opération de Hodge pour les formes extérieures sur la variété pseudo-riemannienne $(P,\gamma)$ ainsi que les opérateurs $\delta$ et $\Delta=d\delta + \delta d$ correspondants. Pour la 1-forme $A$, on a alors :
\begin{equation}
\delta A =0
\end{equation}
\begin{equation}
\Delta A + QA=\varepsilon \ast \nabla \ast F
\end{equation}
où $\varepsilon=\pm1$ suivant les dimensions de $M$ et $G$ et la signature de $\gamma$; $\ast, \nabla, F$ sont comme dans  12.1 et $Q\in \End(\fracg)$ est défini par
\[
1/2 \Tr (ad(E).ad(F))=k(QE,F); \ \ \forall E, F\in \fracg.
\]
(En particulier, si $G$ est semi-simple de dimension $N$, $Q$ est la multiplication par $-\frac{2}{N}\Lambda$).\\
Il en résulte que les équations de Yang et Mills pour $A$ sont équivalentes à
\begin{equation}
\Delta A+QA=0
\end{equation}
sur $(P,\gamma)$.

\subsection{Remarques}

a) L'existence de $k$ sur $\fracg$ (comme dans 12.2, 12.3, 12.4 et 12.5) est une hypothèse restrictive sur $G$ qui implique notamment $\Tr (ad(E))=0$, $\forall E \in \fracg$, (cette relation est d'ailleurs l'origine de l'identité (12.6), $\delta A=0$).\\

b) $-1/2Q$ est l'opérateur correspondant à la courbure de Ricci de $G$.\\

c) Avec nos conventions le laplacien $\Delta$ est un opérateur positif dans le cas où $\gamma$ est positive.\\

d) Si $A$ satisfait aux équations de Yang et Mills, $\Delta A$ est (d'après (12.8)) vertical; il est facile de voir, (12.7), que réciproquement, si $\Delta A$ est vertical, $A$ est solution des équations de Yang et Mills.

\subsection{Exemples}

Il est remarquable que dès que l'on a une situation ``bien symétrique", les équations de Yang et Mills sont automatiquement satisfaites; de ce point de vue, il y a une certaine similarité entre les équations de Yang et Mills pour les fibrés avec connexions et les équations d'Einstein pour les variétés pseudo-riemanniennes, (on notera que 12.4 va un peu dans ce sens). Explicitement, si $M=K/H$ est un espace (pseudo)-riemannien symétrique, ($K$ est le groupe des isométries), alors $K\rightarrow K/H$ est un $H$-fibré principal sur $M$ et sa connexion canonique $K$-invariante (obtenue par projection de la forme Maurer-Cartan de $K$) satisfait les équations de Yang et Mills. Plus généralement, si $P$ est un $G$-fibré principal sur l'espace (pseudo)-riemannien symétrique $M=K/H$ tel que l'action de $K$ sur $M$ se relève en une action sur $P$ par des isomorphismes de $G$-fibrés principaux, alors on a une connexion canonique $K$-invariante sur $P$ et cette connexion satisfait les équations de Yang et Mills.\\
On a une démonstration simple de ces résultats en utilisant les résulats de 12.5 ou 12.6 d).

\subsection{Instantons $SU(2)$ sur $S^4=P_1(\mathbb H)$}

Un cas particulier de la situation symétrique décrite dans 12.7 est la fibration de Hopf $S^7 \xrightarrow[SU(2)]{}S^4$ où $S^4\simeq P_1(\mathbb H)$ est l'ensemble des sous-espaces vectoriels quaternioniques de dimension 1 de $\mathbb H^2$ et $S^7$ est l'ensemble des vecteurs unitaires de $\mathbb H^2$; $SU(2)=Sp(1)$ est le groupe des transformations $\mathbb H$-linéaires isométriques de $\mathbb H$ (voir dans Chapitre 11). Dans ce cas, la courbure $F$ de la connextion canonique invariante par $SO(5)$ ($Sp(2)$ localement) satisfait $F=-\ast F$, (ce qui implique (12.1) en vertu de l'identité de Bianchi); c'est l'anti-instanton invariant par $SO(5)$.\\

Bien qu'étant une solution symétrique, la solution précédente est un cas particulier de toute une classe de solutions des équations de Yang et Mills qui ne sont pas de type symétrique, mais qui sont reliées à l'invariance conforme des équations de Yang et Mills (12.1) lorsque $M$ est de dimension 4.\\

 Si $M$ est de dimension 4, l'opération de Hodge, $\ast$, sur les 2-formes est invariante conforme; il en résulte que les équations (12.1) sont invariantes conformes. Si $M$ est, en plus, proprement riemannienne, les connexions de courbure $F$ telles que
\begin{equation}
F=\pm \ast F
\end{equation}
sont solutions de (12.1) en vertu de l'identité de Bianchi $\nabla F=0$. Les solutions de (12.9) se relèvent en condition d'holomorphie via la transformée de Penrose décrite au chapitre 8, théorème 8.10. \\

Dans le cas où $M=S^4$ et $G=SU(2)$, les solutions de (12.9) sont connues.
Elles sont obtenues en utilisant le théorème 8.10 et la proposition 8.4 pour $\ell=2$.

\subsection{Instantons généraux}

Supposons que $M$ est une variété riemannienne compacte de dimension 4 et que $G$ est un sous-groupe du groupe unitaire $U(p)$. Dans ce cas, nous prendrons pour forme bilinéaire $k$ sur $\fracg$ la forme bilinéaire définie par
\[
k(E,F)=-\Tr(E\ F)\ \ \ \forall E,F\in \fracg,
\]
$E$ et $F$ sont des matrices complexes $p\times p$ anti-hermitiennes).\\

On peut intégrer l'action de Yang et Mills sur $M$ et, l'on a, (inégalité de Schwarz):
\[
-2S^M_{\text{Y.M.}}(A)=\int_M \parallel F\parallel^2\ \text{vol}_M\geq \vert \int_M(F,\ast F)\ \text{vol}_M\vert=\vert \int_M \Tr (F\wedge F)\vert
\]
On a l'égalité si et seulement si $F$ vérifie (12.9).\\
$\int_M\ \Tr (F\wedge F)$est une classe caractéristique ne dépendant que du fibré. On en déduit que, pour un fibré donné, les solutions de (12.9) sont les extrema absolus de $S^M_{\text{Y.M.}}(A)$ et que toute solution de (12.9) avec $F\not= 0$ correspond à un fibré non trivial sur $M$.\\

Si l'on veut généraliser ce qui précède à des dimensions différentes de 4, il y a lieu de généraliser les équations de Yang et Mills; si $M$ est de dimensions $4k$, on peut remplacer l'action de Yang et Mills par
\[
S_{4k}(A)=-1/2\int (F^{\wedge k},F^{\wedge k})\ \text{vol}_M.
\]
Dans le cas où la dimension de $M$ est seulement un multiple de 2 ($2\ell$), c'est une généralisation de certains ``modèles $\sigma$" à 2 dimensions qu'il faut introduire.

\subsection{Sources}

Les équations (12.1) sont les équations de Yang et Mills sans sources. Il est naturel de leur adjoindre un second membre correspondant à l'interaction avec d'autres champs, etc. et d'écrire
\begin{equation}
\nabla\ast F=J
\end{equation}
où $J$ est comme $\nabla\ast F$ une ($n-1$)-forme tensorielle de type $\ad$ satisfaisant
\begin{equation}
\nabla J=0
\end{equation}
i.e. $J$ est à valeurs dans $\fracg$, horizontal satisfaisant $dJ+[A,J]=0$

\subsection{Les équations de Yang et Mills comme lois de conservation sur le fibré principal $P$}

Supposons que $J$ est donné comme $(n-1)$-forme à valeurs dans $\fracg$ horizontale satisfaisant (12.11). Alors (12.10) implique que l'on a sur $P$
\begin{equation}
d(-[A,\ast F]+J)=0
\end{equation}
ce qui peut s'interpréter comme une loi de conservation. Inversement, si (12.12) est vérifié, on en déduit
\[
[A,\nabla \ast F-J]=0
\]
ce qui implique (12.10) dans le cas où $\fracg$ est semi-simple (i.e lorsque la représentation adjointe de $\fracg$ est injective).\\

Dans une section, $J$ correspond à la ``densité de charge de Yang-Mills" de la matière, notion pas intrinsèque car dépendant de la section (sauf dans le cas abélien, pour $U(1)$ cela correspond à la charge électrique). Sur $P$ par contre, c'est bien défini et $-[A,\ast F]$ est l'analogue de $J$ pour le champ Yang et Mills. Dans le cas $\fracg$ semi-simple on a l'équivalence
\[
d(-[A,\ast F])= 0 \Leftrightarrow \nabla \ast F=0
\]
sur le fibré $P$.

\subsection{Les équations de Yang et Mills comme conditions d'intégrabilité}

Considérons l'opérateur de Yang et Mills
\[
A\mapsto \nabla \ast F
\]
sur les connexions $A$ sur $P$. Soit $B$ une 1-forme tensorielle de type $\ad$, i.e. une 1-forme horizontale à valeurs dans $\fracg$ sur $P$ se transformant par la représentation adjointe dans l'action de $G$ sur $P$. Une variation infinitésimale $\delta A=\varepsilon B$ conduit à
\[
\delta (\nabla\ast F)=\varepsilon (\nabla\ast\nabla B+[B,\ast F])
\]
c'est-à-dire à un opérateur différentiel linéaire pour $B$ dont les coefficients dépendent localement de A
\begin{equation}
\cald_A B=\nabla\ast \nabla B+ [B,\ast F]
\end{equation}
Etudions les conditions d'intégrabilité du système linéaire
\begin{equation}
\cald_A B=0
\end{equation}
pour $B$. On a 
\begin{equation}
\nabla \cald_A B = -[B,\nabla\ast F]
\end{equation}
et par conséquent on a 
\begin{equation}
\ad (\nabla\ast F)=0
\end{equation}
comme conditions d'intégrabilité de (12.14). On peut montrer (voir référence D.-V. à la fin  du chapitre 13) que ce sont les seules conditions d'intégrabilité de (12.14). Par conséquent, si $\fracg$ est semi-simple, les équations de Yang et Mills (12.1) sont les conditions d'intégrabilité du système linéaire (12.14).

\newpage
\section{Les équations d'Einstein}
\subsection{Première approche}

Soit $M$ une variété différentiable orientée de dimension $n$ et $g$ une métrique pseudo-riemannienne sur $M$. Soit ($\calo,\varphi$) une carte locale, notons $x^\mu=\varphi^\mu(x)$ ($\mu\in \{1,\dots,n\}$) les coordonnées locales correspondantes de $x\in \calo$. Le champ de repères naturels tangents correspondants sur $\calo$ sera noté $(\partial_\mu)=(\partial_1,\dots,\partial_n)$ et le champ des corepères duaux sera noté($dx^\mu$). La métrique en $x\in \calo$ peut s'écrire
\[
g(x)=g_{\mu\nu}(x)dx^\mu \otimes dx^\nu
\]
où les $g_{\mu\nu}(x)$ sont les composantes de $g$ correspondantes avec $g_{\mu\nu}=g_{\nu\mu}$ pour $\mu,\nu\in \{1,\dots,n\}$. La courbure riemannienne de $g$ a pour composantes
\[
R^\lambda_{\rho\mu\nu}=\partial_\mu \Gamma^\lambda_{\rho\nu}-\partial_\nu \Gamma^\lambda_{\rho\mu} +\Gamma^\lambda_{\sigma\mu} \Gamma^\sigma_{\rho\nu}-\Gamma^\lambda_{\sigma\nu} \Gamma^\sigma_{\rho\mu}
\]
où les $\Gamma^\lambda_{\mu\nu}$ donnés par
\begin{equation}
\Gamma^\lambda_{\mu\nu}=\frac{1}{2}g^{\lambda\rho}(\partial_\mu g_{\nu\rho}+\partial_\nu g_{\mu\rho}-\partial_\rho g_{\mu\nu})
\end{equation}
sont les composantes correspondantes de la connexion riemannienne.\\
Les composantes du tenseur de Ricci sont données par
\[
R_{\mu\nu}=R^\lambda_{\mu\lambda\nu}
\]
et la courbure scalaire est
\[
R=g^{\mu\nu} R_{\mu\nu}
\]
où les $g^{\mu\nu}$ sont définis par
\[
g^{\mu\rho}g_{\rho\nu}=\delta^\mu_\nu
\]
pour $\mu,\nu\in \{1,\dots,n\}$. Avec ces notations, dans le repère naturel de la carte, le tenseur d'Einstein $G$ est donné par 
\begin{equation}
G_{\mu\nu}=R_{\mu\nu}-\frac{1}{2}R g_{\mu\nu}
\end{equation}
$\mu, \nu\in \{1,\dots,n\}$. L'annulation de la torsion de la connexion riemannienne
implique la symétrie $R_{\mu\nu}=R_{\nu\mu}$ du tenseur de Ricci et par conséquent aussi la symétrie $G_{\mu\nu}=G_{\nu\mu}$ du tenseur d'Einstein. Les identités de Bianchi impliquent alors, par contraction
\begin{equation}
g^{\lambda\mu}\nabla_\lambda G_{\mu\nu}=0
\end{equation}
i.e. l'annulation de la divergence covariante du tenseur d'Einstein.\\

La métrique $g$ satisfait aux {\sl équations d'Einstein (sans source)} si l'on a
\begin{equation}
G_{\mu\nu}=0
\end{equation}
pour $\mu,\nu\in \{1,\dots,n\}$. Ces équations en $x\in M$ ne dépendent pas de la carte locale ($\calo,\varphi$) choisie avec $x\in \calo$.

\subsection{L'action d'Einstein-Hilbert}

Les équations d'Einstein (13.3) dérivent de {\sl l'action d'Einstein-Hilbert}
\[
S_{\text{E.H.}}(g)=\frac{1}{4} \int R\  \vol_M
\]
i.e. on a
\begin{equation}
\delta S_{\text{E.H.}}=\frac{1}{4}\int G_{\mu\nu}\ \delta g^{\mu\nu}\ \vol_M
\end{equation}
modulo un ``terme de surface" $\int d\omega$ où $\omega$ est une ($n-1$)-forme.\\

Notons que les équations (13.3) sont identiquement vérifiées en dimension $n=2$ et que pour $n\geq 3$ elles sont équivalentes à $R_{\mu\nu}=0$, i.e. à l'annulation de la courbure de Ricci de la variété pseudo-riemannienne $M$.

\subsection{Repères orthonormés}

Soit $\GL(M)$ le $GL(n,\mathbb R)$-fibré principal sur $M$ des repères tangents. En chaque $x\in M$, on peut trouver un repère tangent $(e_1,\dots,e_n)$ tel que
\begin{equation}
g(e_i,e_j)=\eta_{ij}
\end{equation}
$(\forall i,j\in \{1,\dots,n\})$ où $\eta_{ij}=\eta_{ii} \delta_{ij}$ est diagonale avec ses $s$ premiers éléments diagonaux égaux à 1 et les autres $n-s$ éléments diagonaux égaux à $-1$. L'entier $s\in \{0,\dots,n\}$ est un invariant de la variété pseudo-riemannienne $(M,g)$ indépendant du point $x\in M$. Un repère tangent $(e_1,\dots,e_n)$ satisfaisant (13.6) sera appelé {\sl repère orthonormé}. L'ensemble $\O(M,g)$ de ces repères orthonormés est un sous-fibré de $\GL(M)$ qui est un $O(s,n-s)$-fibré principal sur $M$. Par définition $\O(M,g)$ est la $O(s,n-s)$-structure sur $M$ correspondant à la pseudo-métrique $g$ (voir Chapitre 1).\\

Rappelons que sur $\GL(M)$ on a une 1-forme à valeurs dans $\mathbb R^n$ notée $\theta=(\theta^i)$, la {\sl forme de soudure},  définie de la manière suivante. Pour tout vecteur tangent $X$ en $r$ à $\GL(M)$, $\theta^k(X)$ est la $k$-ième composante du projeté $\pi_\ast(X)$ dans le repère $r$, $\pi:\GL(M)\rightarrow M$ étant la projection du fibré $GL(M)$ sur $M$.\\

Soit $\omega=(\omega^i_j)$ une 1-forme de connexion sur $\GL(M)$. Une telle connexion existe toujours et les $n+n^2=n(n+1)$ formes ($\theta^k, \omega^i_j$) définissent un parallélisme absolu de la variété $\GL(M)$, ($\dim\  \GL(M)=n(n+1)$). {\sl Les équations de structures}
\begin{equation}
d\theta^i+\omega^i_j\wedge \theta^j=\Theta^i
\end{equation}
\begin{equation}
d\omega^i_j+\omega^i_k \wedge \omega^k_j=\Omega^i_j
\end{equation}
définissent les 2-formes de torsion $\Theta^k$ et de courbure $\Omega^i_j$. Par différentiation on obtient {\sl les identités de Bianchi}
\begin{equation}
d\Theta^i+\omega^i_j \wedge \Theta^j=\Omega^i_j \wedge \theta^j
\end{equation}
\begin{equation}
d\Omega^i_j+\omega^i_k\wedge \Omega^k_j-\Omega^i_k\wedge \omega^k_j=0
\end{equation}
pour la torsion et la courbure de la connexion.\\

La forme de soudure $\theta$ se restreint à $\O(M,g)$, ainsi qu'à toute $G$-structure sur $M$, et nous désignerons encore par $\theta^k$ les restrictions à $\O(M,g)$ de ses composantes. Avec ces conventions, la métrique s'écrit
\begin{equation}
g=\eta_{ij} \theta^i\otimes \theta^j
\end{equation}
expression qui est indépendante du repère orthonormé $r\in \O(M,g)$ choisi au-dessus de $x\in M$. Supposons que la connexion $\omega$ préserve la métrique $g$. Alors les $\omega^i_j$ restreints à $\O(M,g)$ définissent une connexion sur $\O(M,g)$ et la préservation de $g$ s'écrit
\[
d\eta_{ij}=0=\eta_{ik}\omega^k_j+\eta_{kj} \omega^k_i
\]
autrement dit en posant
\begin{equation}
\omega_{ij}=\eta_{ik}\omega^k_j
\end{equation}
on a 
\begin{equation}
\omega_{ij}+\omega_{ji} =0
\end{equation}
$\forall i,j\in \{1,\dots,n\}$. Par la suite nous descendrons et nous montrerons les indices $i\in \{1,\dots,n\}$ sur $\O(M,g)$ en contractant avec les $\eta_{ik}$ et les $\eta^{ik}$ où $\eta^{ik}\eta_{kj}=\delta^i_{j}$ par exemple on posera $\omega^{ij}=\omega^i_k \eta^{kj}$.\\
On notera qu'en vertu de (13.13) on a seulement $\frac{n(n-1)}{2}$ 1-formes indépendantes sur $\O(M,g)$ données par les $\omega^i_j$ et que les $\frac{n(n+1)}{2}$ 1-formes indépendantes
$\theta^k$, $\omega^{ij} (i<j)$ définissent un parallélisme absolu de la variété $\O(M,g)$.\\

Dans les hypothèses précédentes, les équations de structures et les identités de Bianchi se restreignent à $\O(M,g)$ et nous désignerons par les mêmes symboles $\Theta^k$, $\Omega^i_j$ les 2-formes de torsion et de courbure sur $\O(M,g)$. On a 
\begin{equation}
\Omega_{ij}+\Omega_{ji}=0
\end{equation}
sur $\O(M,g)$. Les relations (13.13) et (13.14) expriment le fait que les formes $\omega$ et $\Omega$ sur $\O(M,g)$ sont à valeurs dans l'algèbre de Lie de $O(s,n-s)$.\\

Il existe une connexion et une seule de torsion nulle sur $\O(M,g)$ i.e. une connexion et une seule préservant la métrique de torsion nulle, c'est la connexion riemannienne $\mathring\omega=(\mathring\omega^i_j)$. Dans le repère naturel associé à une carte locale (section locale de $\GL(M))$ cette connexion est donnée par
\begin{equation}
\mathring\omega^\lambda_\mu=\Gamma^\lambda_{\mu\nu} dx^\nu
\end{equation}
où les $\Gamma^\lambda_{\mu\nu}$ sont donnés par (13.1).

\subsection{Actions d'Einstein-Hilbert et d'Einstein-Cartan}

Définissons les ($n-q$)-formes sur $\O(M,g)$
\begin{equation}
\theta^\ast_{i_1\dots i_q}=\frac{1}{(n-q)!}\varepsilon_{i_1\dots i_qi_{q+1}\dots i_n}\theta^{i_{q+1}}\wedge \dots \wedge \theta^{i_n}
\end{equation}
où $\varepsilon_{i_1\dots i_n}$ est complètement antisymétrique avec $\varepsilon_{1\dots n}=1$.\\

Considérons la $n$-forme sur $\O(M,g)$
\begin{equation}
L=\frac{1}{4}\Omega^{ij}\wedge \theta^\ast_{ij}
\end{equation}
où $\Omega^i_j$ est la courbure d'une connexion $\omega^i_j$ sur $\O (M,g)$, (i.e. $\omega$ est une connexion métrique satisfaisant par conséquent (13.13)). C'est une forme horizontale invariante par l'action de $O(s,n-s)$ donc basique. Autrement dit $L$ est l'image inverse par $\pi:\O(M,g)\rightarrow M$ d'une $n$-forme $\call$ sur $M$, $L=\pi^\ast(\call)$.\\

 \'Etant donnée une section de $\O(M,g)$ sur l'ouvert $\calo\subset M$, l'intégrale de $L$ sur la section ne dépend pas de la section choisie et coïncide avec l'intégrale de $\call$ sur $\calo$.\\
 
 Nous noterons par
 \[
 S=\int L =\frac{1}{4}\int \Omega^{ij}\wedge \theta^\ast_{ij}
 \]
 l'intégrale (indéfinie $\calo\mapsto S^\calo$) correspondante. $S$ est a priori une fonctionnelle locale de la métrique $g$ et de la connexion métrique $\omega$.\\
 
 Dans le cas où on prend pour $\omega$ la connexion riemannienne $\mathring{\omega}$ on obtient une fonctionnelle de $g$ qui n'est autre que l'action d'Einstein-Hilbert, 
 \[
 S_{E.H.}(g)=S(g,\mathring{\omega})
 \]
 car $\call \lvert_{\omega=\mathring{\omega}}=\frac{1}{4} R\ \vol_M$.\\
 
 Dans le cas général, il est plus commode de prendre comme variable locale une section du fibré des repères $\GL(M)$ ou plutôt de manière équivalente les restrictions $\theta^k(x)$ des composantes de la forme de soudure $\theta$ à cette section de $\GL(M)$, la métrique $g$ étant donnée par
 \[
 g(x)=\eta_{ij}\theta^i(x) \theta^j(x)
 \]
 et la connexion métrique dans la section étant donnée par des 1-formes $\omega^i_j(x)$ telles que l'on ait (13.13) avec (13.12), i.e. $\omega_{ij}(x)=-\omega_{ji}(x)$. L'action correspondante pour les 1-formes $x\mapsto \theta^k(x), x\mapsto\omega_{ij}(x)$ sur $M$ est {\sl l'action d'Einstein-Cartan}
 \begin{equation}
S_{\text{E.C.}}(\theta,\omega)=\frac{1}{4}\int \Omega^{ij}\wedge \theta^\ast_{ij}
\end{equation}
où $\Omega^i_j(x)$ et $\theta^\ast_{ij}(x)$ sont toujours donnés par
\[
\Omega^i_j(x)=d\omega^i_j(x)+\omega^i_k(x)\wedge \omega^k_j(x)
\]
\[
\theta^\ast_{ij}(x)=\frac{1}{(n-2)!}\varepsilon_{ijk_1\dots k_{n-2}}\theta^1(x)\wedge\dots\wedge \theta^{n-2}(x),
\]
$\frac{1}{n!}\varepsilon_{i_1\dots i_n}\theta^{i_1}(x)\wedge\dots\wedge \theta^{i_n}(x)=\theta^1(x)\wedge\dots\wedge\theta^n(x)$
étant la forme élément de volume $\vol_M(x)$ de la variété pseudo-riemannienne orientée $M$ de métrique $g(x)$ en $x$.

\subsection{Equations sans source}

On suppose que $\dim(M)=n\geq 3$. On a pour la variation de l'action d'Einstein-Cartan
\[
\delta S_{\text{E.C.}}=\frac{1}{4}\int d(\delta\omega^{ij}\wedge \theta^\ast_{ij})+\frac{1}{4}\int (\delta\omega^{ij}\wedge\Theta^k+\Omega^{ij}\wedge\delta\theta^k)\wedge\theta^\ast_{ijk}
\]
et par conséquent les équations qui dérivent de cette action sont
\begin{equation}
\Theta^k\wedge\theta^\ast_{ijk}=0
\end{equation}
et
\begin{equation}
\Omega^{ij}\wedge\theta^\ast_{ijk}=0
\end{equation}
Les premières (13.19) sont équivalentes à l'annulation de la torsion $\Theta^k=0$ et impliquent que $\omega$ est la connexion riemannienne $\mathring{\omega}$; les dernières (13.20) se réduisent alors aux équations d'Einstein (13.4). Cela découle également de l'identité
\[
-\frac{1}{2}\Omega^{ij}\wedge \theta^\ast_{ijk}=(\bar R^i_k-\frac{1}{2}\bar R\delta^i_k)\theta^\ast_i
\]
où $\bar R^i_k=\bar R^{mi}_{mk}$avec $\Omega^{ij}=\frac{1}{2}\bar R^{ij}_{km}\theta^k\wedge \theta^m$ et $\bar R=\bar R^m_m$ et du fait que pour $\omega^i_j=\mathring{\omega}^i_j$, $\bar R^i_k-\frac{1}{2}\bar R\delta^i_k=R^i_k-\frac{1}{2}R\delta^i_k$ est le tenseur d'Einstein $G^i_k$.\\

On voit donc que sans source ou couplage avec d'autres champs, les équations qui dérivent de l'action d'Einstein-Cartan sont les mêmes que celles qui dérivent de l'action d'Einstein-Hilbert c'est-à-dire les équations d'Einstein (sans source) décrites dans 13.1. Il faut souligner à ce sujet qu'il existe une infinité d'actions locales naturelles (i.e. covariantes, etc.) dont dérivent les équations d'Einstein sans source. C'est le cas en particulier pour l'action de Holst utilisée pour la gravité quantique en boucles (pour les valeurs génériques du paramètre d'Immirzi).

 \subsection{Les équations d'Einstein comme conditions d'inté\-grabilité}
 
 Considérons l'opérateur d'Einstein $g\mapsto G_{\mu\nu}(g)$ sur les métriques $g$ sur $M$. Soit $h$ un champ de tenseurs (0,2) symétrique $x\mapsto h_{\mu\nu}(x)$ sur $M$. Une variation infinitésimale $\delta g=\varepsilon h$ de la métrique conduit à
 \[
 \delta(G_{\mu\nu})=\varepsilon G'_{\mu\nu} h
 \]
 où $G'_{\mu\nu}$ est un opérateur différentiel linéaire pour $h$ dont les coefficients dépendent localement de $g$.\\

\'Etudions les conditions d'intégrabilité du système linéaire
\begin{equation}
G'_{\mu\nu} h=0
\end{equation}
pour $h$. La variation de l'identité (13.6), $\nabla^\mu G_{\mu\nu}=0$, conduit à 
\[
\nabla^\mu(G'_{\mu\nu}h)-h_{\lambda\mu} \nabla^\lambda G^\mu_\nu-(g^{\lambda\rho} G^\mu_\nu+g^{\mu\rho} G^\lambda_\nu + G^{\lambda\mu} \delta^\rho_\nu - G^\rho_\nu g^{\lambda\mu})\frac{1}{2} \nabla_\rho h_{\lambda\mu}=0
\]
pour $\nu\in \{1,\dots,n\}$.\\

Il en résulte que l'intégrabilité de (13.21) implique $(\nabla^\mu(G'_{\mu\nu}h)=0)$
\[
\nabla^\lambda G^\mu_\nu=0
\]
et 
\begin{equation}
g^{\lambda\rho}G^\mu_\nu+g^{\mu\rho} G^\lambda_\nu+G^{\lambda\mu}\delta^\rho_\nu-G^\rho_\nu g^{\lambda\mu}=0
\end{equation}
$\forall \lambda, \mu, \nu, \rho\in \{1,\dots,n\}$. Les équations (13.22) impliquent
\[
G_{\mu\nu}=0
\]
c'est-à-dire les équations d'Einstein (sans source) pour $g$. On peut montrer (voir référence D.-V. à la fin) que ce sont les seules conditions d'intégrabilité de (13.21).

\subsection{Remarques}

\noindent a) Les équations d'Einstein $G_{\mu\nu}=0$ sont également les conditions nécessaires et suffisantes pour que les équations linéaires (13.21) pour $h$ soient invariantes par
\[
h_{\mu\nu} \mapsto h_{\mu\nu}+\nabla_\mu X_\nu + \nabla_\nu X_\mu
\]
pour tout champ de vecteurs $X$ sur $M$.\\

\noindent b) Les équations (13.21) sont les équations pour un champ $h$ de spin 2 et de masse nulle dans le champ de gravitation extérieur $g$.\\

\noindent c) On a une analyse similaire pour les équations (13.19), (13.20) en repère orthonormé.

\subsection{Généralisation des paires de Lax}

En utilisant les résultats de 13.6 pour les équations d'Einstein et de 12.12 pour les équations de Yang et Mills, on peut montrer (référence D.-V. à la fin de ce chapitre) que ces équations sont équivalentes à des conditions de courbure nulle du type
\begin{equation}
\partial_\mu\cala_\nu-\partial_\nu \cala_\mu+[\cala_\mu,\cala_\nu]=0
\end{equation}
$\mu, \nu\in \{1,\dots,n\}$, où les $\cala_\mu(x)$ sont les composantes en $x\in M$ d'une connexion sur un fibré vectoriel sur $M$ qui ne dépendent que des champs et de leurs dérivés première en $x$, c'est-à-dire que l'on a par exemple
\[
\cala_\mu(x)=F_\mu(g_{\alpha\beta}(x),\partial_\rho g_{\tau\sigma}(x))
\]
pour les équations d'Einstein.\\

Les $\cala_\mu$ satisfaisant (13.23) constituent une généralisation en dimension $n$ pour les équations d'Einstein ou pour les équations de Yang et Mills des paires de Lax pour les systèmes intégrables en dimension 2. Ce qu'il manque ici par rapport aux systèmes intégrables en dimension 2 est l'analogue des équations de Yang-Baxter classiques.

\subsection{Les équations d'Einstein comme lois de conservation sur le fibré de repères orthonormés}

Soit $\omega=(\omega^i_j)$ une 1-forme de connexion sur le fibré $\O(M,g)$ des repères orthonormés de la variété pseudo-riemannienne $(M,g)$.\\

Définissons sur $\O(M,g)$ la $(n+1)$-forme $\tau=(\tau_i)$ à valeurs dans $\mathbb R^n$ par
\begin{equation}
\tau_i=-\frac{1}{2}(\omega^j_i \wedge \omega^{k\ell}\wedge \theta^\ast_{jk\ell}+\omega^j_\ell\wedge \omega^{\ell k} \wedge \theta^\ast_{ijk})
\end{equation}
et la $(n-2)$-forme $\sigma=(\sigma_i)$ à valeurs dans $\mathbb R^n$ par
\begin{equation}
\sigma_i=-\frac{1}{2}\omega^{jk}\wedge \theta^\ast_{ijk}
\end{equation}
$\forall i\in \{1,\dots,n\}$. On a les identités suivantes 
\begin{equation}
d\sigma_i=\tau_i+\frac{1}{2}(\Theta^j \wedge \omega^{k\ell}\wedge \theta^\ast_{ijk\ell}-\Omega^{jk}\wedge \theta^\ast_{ijk})
\end{equation}
et 
\begin{eqnarray}
d\tau_i &=&-\Theta^m\wedge \frac{1}{2}(\omega^\ell_i\wedge \omega^{jk}\wedge \theta^\ast_{jk\ell m}\nonumber +\omega^j_\ell \wedge \omega^{\ell k} \wedge \theta^\ast_{ijkm}) \nonumber \\
&+&\frac{1}{2} (\omega^{jk}\wedge \Omega^\ell_n \wedge \theta^n \wedge \theta^\ast_{ijk\ell}+
\omega^j_i \wedge \Omega^{k\ell} \wedge \theta^\ast_{jk\ell})
\end{eqnarray}
où les $\Theta^k$ et les $\Omega^i_j$ sont les 2-formes de torsion et de courbure de $\omega$.\\

En annulant les composantes des différents degrés verticaux, on vérifie que $d\tau_i=0$ est équivalent à la torsion nulle plus les équations d'Einstein sans source pour $g$. Plus précisément, on a le théorème suivant.

\subsection {Théorème}

{\sl Les conditions $\mathrm{a),\  b),\ c)}$ suivantes sont équivalentes :\\
$\mathrm{a)}$ $d\tau_i=0$\\
$\mathrm{b)}$ $\tau_i=d\sigma_i$\\
$\mathrm{c)}$ $\omega$ est la connexion riemannienne (i.e. $\Theta=0$) et $g$ satisfait aux équations d'Einstein (sans source)}.\\

\noindent\underbar{Démonstration}. Supposons $d\tau_i=0$. En égalant à 0 les différents degrés verticaux du second membre de (13.27), on obtient
\[
(\omega^n_i\wedge \omega^{jk}\wedge \theta^\ast_{njk\ell} + \omega^j_m \wedge \omega^{mk}\wedge \theta^\ast_{ijk\ell})\wedge \Theta^\ell=0
\]
et 
\[
\omega^{jk}\wedge \Omega^\ell_n\wedge \theta^n\wedge \theta^\ast_{ijk\ell} + \omega^n_i\wedge \Omega^{jk}\wedge \theta^\ast_{njk}=0
\]
ce qui est équivalent à 
\[
\Theta^k=0
\]
et 
\[
\Omega^{ij}\theta^\ast_{ijk}=0
\]
$\forall k\in \{1,\dots,n\}$, où on a utilisé l'identité de Bianchi (13.9) et l'annulation de $\Theta$ qui impliquent $\Omega^\ell_n\wedge \theta^n=0$.\\

Ces dernières équations sont équivalentes aux équations d'Einstein, comme on l'a vu en 13.5, elles impliquent d'autre part en vertu de (13.26)
\[
\tau_i=d\sigma_i
\]
$\forall i\in \{1,\dots,n\}$, ce qui implique évidemment $d\tau_i=0$. On a donc montré $\mathrm{a)} \Rightarrow \mathrm{c)} \Rightarrow \mathrm{b)} \Rightarrow \mathrm{a)}$. $\square$

\subsection{Sources}

Pour décrire les interactions du champ de gravitation qui s'identifie à la métrique $g$ on adjoint un second membre aux équations sans source et on écrit (avec nos conventions $c=1, 4\pi G=1$)
\begin{equation}
G_{\mu\nu}=2T_{\mu\nu}
\end{equation}
où $G_{\mu\nu}$ est le tenseur d'Einstein de $g$ et $T_{\mu\nu}$ est le tenseur d'impulsion-énergie représentant la densité d'impulsion-énergie de ``la matière" (particules, champ électromagnétique, etc.). Ce tenseur satisfait à 
\begin{equation}
g^{\lambda\mu} \nabla_\lambda T_{\mu\nu}=0
\end{equation}
ce qui la trace de ce qui reste de la conservation de l'impulsion énergie de la matière en présence du champ de gravitation. Ces équations (13.29) sont aussi nécessaires pour (13.28) en vertu des identités (13.3). Notons aussi que $T$ est symétrique; en fait $T_{\mu\nu}$ est la dérivée variationnelle de l'action $S_{Mat}$ de la matière par rapport à $g^{\mu\nu}$, $T_{\mu\nu}=\delta S_{Mat}/\delta g^{\mu\nu}$.\\

En repère orthonormé, on introduit la $(n-1)$-forme horizontale à valeurs dans $\mathbb R^n$
\begin{equation}
t_i=T^j_i\theta^\ast_j
\end{equation}
qui est tensorielle de type représentation fondamentale dans $\mathbb R^n$. Les équations (13.29) s'écrivent alors
\begin{equation}
dt_i-\omega^j_i \wedge t_j=0
\end{equation}
$\forall i\in \{1,\dots,n\}$, pour $\omega =\mathring\omega$ (i.e. $\Theta=0$). Les équations d'Einstein avec source dans $\O(M,g)$ deviennent
\begin{equation}
\left\{
\begin{array}{l}
\Theta^i=0\\
\Omega^{jk}\wedge \theta^\ast_{ijk}+4t_i=0
\end{array}
\right.
\end{equation}
pour $i\in \{1,\dots,n\}$.

\subsection{Théorème}

{\sl Soit $\omega=(\omega^i_j)$ une forme de connexion sur $\O(M,g)$ et soit $t=(t_i)$ une $(n-1)$-forme à valeurs dans $\mathbb R^n$ sur  \O$(M,g)$ telles que l'on ait $dt_i=\omega^j_i\wedge t_j=0$. Alors les conditions $\mathrm{a), b), c)}$ suivantes sont équivalentes :\\

$\mathrm{a)}$ : $d(\tau_i+t_i)=0$\\

$\mathrm{b)}$ : $\tau_i+t_i=d\sigma_i$\\

$\mathrm{c)}$ : $\Theta^i=0$ et $\Omega^{jk}\wedge \theta^\ast_{ijk}+ 4 t_i=0$\\

\noindent où $\Theta$ et $\Omega$ sont les formes de torsion et de courbure de $\omega$}.

La démonstration est la même que pour le théorème (13.10) en utilisant $dt_i=\omega^j_i\wedge t_j$ et en identifiant les différents degrés verticaux.\\

La condition c) signifie que l'on a les équations d'Einstein avec source $t$ pour $g$ (et que $\omega$ est la connexion riemannienne).

\subsection{Interprétation physique}

Dans le cas où $(M,g)$ est l'espace-temps et $t=(t_i)$ est construit avec le tenseur d'impulsion-énergie de la matière, $t=(t_i)$ représente dans n'importe quelle section locale de $\O(M,g)$ la densité d'impulsion-énergie de la matière (i.e. de tout ce qui n'est pas le champ de gravitation $g$). Les équations $d(\tau_i+t_i)=0$ peuvent donc s'interpréter dans n'importe quelle section locale comme la conservation locale de l'impulsion énergie totale en attribuant $\tau_i$ à la partie gravitationnelle.\\

En fait, les images de $\tau_i$ dans les différentes sections correspondent à des pseudo-densités associées à différents pseudo-tenseurs d'impulsion-énergie pour la gravitation.\\

Il est possible d'étendre les définitions des formes $\tau_i$ et $\sigma_i$ à valeurs dans $\mathbb R^n$ sur le fibré $\GL(M)$ de tous les repères (en introduisant les composantes de la métrique dans les différents repères) de manière à avoir les théorèmes 13.10 et 13.12 sur le fibré $\GL(M)$, $\omega$ étant toujours une connexion métrique. Les pseudo-densités associées aux différents pseudo-tenseurs d'impulsion-énergie pour la gravitation correspondent aux images de $\tau_i$ dans les différentes sections.\\

On peut donc dire que bien qu'il ne puisse pas exister sur $M$ de notion d'impulsion-énergie locale, il existe un objet canonique sur les fibré des repères dont la conservation locale est équivalente aux équations d'Einstein.

\subsection{Théorie d'Einstein et d'Einstein-Cartan}

Nous avons vu dans 13.5 que les équations variationnelles déduites de l'action d'Einstein-Hilbert et celles déduites de l'action d'Einstein-Cartan coïncident et se réduisent aux équations d'Einstein sans source.\\

Si on ajoute à ces actions des actions pour d'autres champs qui ne dépendent pas de la connexion $\omega$ comme l'action de Maxwell ou de Yang-Mills ou même d'un champ scalaire sur $(M,g)$ elles restent équivalentes et conduisent aux mêmes équations couplées pour ces champs et le champ de gravitation.\\

Il n'en est pas de même pour les actions faisant intervenir des champs de spineurs par exemple car dans ce cas il est naturel dans le cas Einstein-Cartan de faire intervenir les dérivées covariantes relatives à la connexion $\omega$. Dans un tel cas les équations variationnelles relativement à $\omega$ n'entraînent plus $\Theta^i=0$ mais déterminent la torsion en fonction du ``tenseur de spin" du nouveau champ. Cette torsion ne se propage pas de sorte que les équations se réduisent toujours à des équations couplées entre la métrique et les nouveaux champs, mais elles diffèrent de celles obtenues en partant de l'action d'Einstein-Hilbert et des nouveaux champs dans l'action desquels seule la métrique apparaît (à travers la connexion riemannienne en particulier).\\

Il faut cependant souligner deux points importants :\\

1. Le principe de couplage minimal est ambigu dans la théorie d'Einstein-Cartan,\\

2. Les écarts entre les 2 théories sont très petits et indécelables aux mesures actuelles.\\

Ces deux versions de la théorie sont donc parfaitement viables dans l'état de nos connaissances actuelles; il en est de même des théories obtenues en rajoutant à l'action d'Einstein-Cartan le terme de Holst ou d'autres termes similaires donnant les mêmes équations sans sources.

\newpage
\begin{center}
{\large Références pour les chapitres 12 et 13}
\end{center}
\vspace{1cm}
\begin{itemize}
\item
Atiyah M.F., Hitchin N.J., Drinfeld V.G., Manin Yu.I : Construction of instantons, {\sl Phys. Lett.\/ \bf A 65} (1978) 185-187.\\

\item
Atiyah M.F., Hitchin N.J, Singer I.M. : Self-duality in four-dimensional Riemannian geometry, {\sl Proc. R. Soc. Lond.\/ \bf A 362} (1978) 425-461.\\

\item
Cho Y.M. : Higher-dimensional unifications of gravitation and gauge theories, {\sl J. Math. Phys.\/ \bf 16} (1975) 2029-2035.\\

\item  
Dubois-Violette M. :
The theory of linear over-determined systems and its applications to
  non-linear field equations,
{\sl J. Geom. Phys.\/ \bf 1} (1984) 139-172.\\

\item 
Dubois-Violette M., Madore J. :
Conservation laws and integrability conditions for gravitational and {Y}ang-{M}ills field equations,
{\sl Commun. Math. Phys.\/ \bf 108} (1987) 213-223.\\

\item 
Dubois-Violette M., Lagraa M.: Abundance of local actions for the vacuum Einstein equation, {\sl Lett. Math. Phys.\/ \bf 91} (2010) 83-91.\\

\item
Harnad J, Tafel J., Shnider S. : Canonical connections on Riemannian symmetric spaces and solutions to the Einstein-Yang-Mills equations, {\sl J. Math. Phys.\/ \bf 21} (1980) 2236-2240.\\

\item
Kaluza T. : Zum unitätsproblem der Physik, {\sl Sitz. ber. Preuss. Akad. Wiss. Berl. \/ \bf 54} (1921) 966-972.\\

\item
Kerner R. : Generalization of the Kaluza-Klein theory for an arbitrary non-abelian gauge group. {\sl Annales de l'I.H.P., section A\/ \bf 9} (1968) 143-152.\\

\item
Klein O. : Quantentheorie und fünfdimensionale Relativitätstheorie, {\sl Z. Physik\/ \bf A37} (1926) 895-906.\\

\item
Lichnerowicz A. : Théories relativistes de la Gravitation et de l'\'Electro\-magnétisme, Masson, Paris 1955.\\

\item
Nowakowski J., Trautman A. : Natural connections on Siefel bundles are sourceless gauge fields, {\sl J. Math. Phys.\/ \bf 19} (1978) 1100-1103.\\

\item
Sparling G.A.J. : Twisters, spinors and the Einstein vacuum equations, University of Pittsburg preprint (1984).\\

\item
Stora R. : Continuum gauge theories, in New Developments in Quantum Field Theory and Statistical Mechanics, Cargèse 1976, Ed. M. Lévy and P. Mitter, Plenum Press 1977, pp. 201-224.\\

\item
Trautman A. : Fiber bundles, Gauge fields and Gravitation, in General Relativity and Gravitation, Vol. I, Ed. A. Held, Plenum Press 1980, pp. 287-308.\\

\item
Trautman A. : Einstein-Cartan theory, in Encylopedia of Mathematical Physics, Vol. 2, eds J.-P. Françoise et al., Elsevier 2006, pp. 189-195.

\end{itemize}
\vspace{2cm}
\noindent Voir aussi dans ce volume les cours de Gérard Clément, Jean-Pierre Provost et de Simone Speziale.
\end{document}